\documentclass[10pt]{article}

\usepackage{amssymb}
\usepackage{amsmath}
\usepackage{amscd}
\usepackage[dvips]{graphicx}
\usepackage{mathrsfs}
\usepackage{wrapfig} 
\usepackage{picins}
\usepackage{picinpar}
\usepackage{psfig}
\usepackage{setspace}
\numberwithin{equation}{section}

\newtheorem{theorem}{Theorem}
\newtheorem{lemma}{Lemma}

\begin{document}

\title{Some Asymptotic Results for the Transient Distribution of the Halfin-Whitt Diffusion Process}

\author{
Qiang Zhen\thanks{
Department of Mathematics and Statistics,
University of North Florida, 1 UNF DR, Jacksonville, FL 32224, USA.
{\em Email:} q.zhen@unf.edu.}
\and
and
\and
Charles Knessl\thanks{
Department of Mathematics, Statistics, and Computer Science,
University of Illinois at Chicago, 851 S. Morgan St.(M/C 249),
Chicago, IL 60607, USA.
{\em Email:} knessl@uic.edu.\
\newline\indent\indent{\bf Acknowledgement:} C. Knessl was partly supported by NSA grants H 98230-08-1-0102 and H 98230-11-1-0184.
}}
\date{ }
\maketitle

\begin{abstract}
\noindent We consider the Halfin-Whitt diffusion process $X_d(t)$, which is used, for example, as an approximation to the $m$-server $M/M/m$ queue. We use recently obtained integral representations for the transient density $p(x,t)$ of this diffusion process, and obtain various asymptotic results for the density. The asymptotic limit assumes that a drift parameter $\beta$ in the model is large, and the state variable $x$ and the initial condition $x_0$ (with $X_d(0)=x_0>0$) are also large. We obtain some alternate representations for the density, which involve sums and/or contour integrals, and expand these using a combination of the saddle point method, Laplace method and singularity analysis. The results give some insight into how steady state is achieved, and how if $x_0>0$ the probability mass migrates from $X_d(t)>0$ to the range $X_d(t)<0$, which is where it concentrates as $t\to\infty$, in the limit we consider. We also discuss an alternate approach to the asymptotics, based on geometrical optics and singular perturbation techniques.  
\\

\end{abstract}

\section{Introduction}

We consider a diffusion process $X_d(t)$ whose state space is $(-\infty, \infty)$, the process has constant drift $=-\beta$ in the range $X_d>0$, and a linear drift 
$=-X_d-\beta$ in the range $X_d<0$, which acts as a ``restoring force". Thus the drift function is continuous at $X_d=0$ but not smooth there. This process arises naturally as a heavy traffic limit of the discrete $m$-server Markovian $M/M/m$ queue, as well as the $GI/M/m$ model. When the number of servers $m$ is large, Halfin and Whitt \cite{hal} showed that the number of customers $N(t)$ in the system may be approximated by $N(t)\approx m+\sqrt{m}\,X_d(t)$. Here we must also scale $\lambda/(\mu\,m)=1-\beta/\sqrt{m}=1-O(m^{-1/2})$, where $\lambda$ and $\mu$ are the arrival and service rates in the model. Thus the drift parameter $\beta$ is defined as the limit of $\sqrt{m}\big[1-\lambda/(\mu\,m)\big]$ as $m\to \infty$. Note that time $t$ is not scaled in this diffusion limit. The scaling $\lambda/(\mu\,m)=1+O(m^{-1/2})$ was also used earlier by Erlang \cite{bro} and Pollaczek \cite{pol} for related models. Note that the diffusion model is stable for $\beta>0$. Assuming that $X_d(0)=x_0$ with probability one, we let $p(x,t)=p(x,t;x_0,\beta)$ be the density of this diffusion process, which satisfies $p(x,0)=\delta(x-x_0)$, and the forward equation in (\ref{2_pde}). 

For $\beta>0$ it is relatively easy to obtain the steady state limit $p(x,\infty)$ (see (\ref{2_steady}) below), but the transient density is much more complicated. In \cite{lee} integral representations are given for $p(x,t)$, where it proves necessary to consider separately the four cases $x\lessgtr 0$ and $x_0\lessgtr 0$. The integrals involve inverse Laplace transforms and parabolic cylinder functions. In view of the complexity of the transient density, it is useful to evaluate it in certain asymptotic limits, as the asymptotic results will presumably be simpler and yield more insight. One asymptotic limit worthy of consideration is that of large time, where we estimate the approach to equilibrium. In \cite{lee} and \cite{gam} estimates are given of the form $p(x,t)-p(x,\infty)\approx e^{-r(\beta)t}$ for $t\to\infty$ where $r(\beta)$ is the ``relaxation rate", which depends on the drift parameter. In particular it was shown in \cite{gam} that $r(\beta)=\beta^2/4$ if $\beta\le \beta_*$ and $r(\beta)=r_0(\beta)$ if $\beta\ge \beta_*$, where $\beta_*\approx 1.85722$ is the minimal positive root of $D'_{\beta^2/4}(-\beta)$. Here $D$ is the parabolic cylinder function $D_p(z)$ and $D'_p(z)=\frac{d}{dz}D_p(z)$ is its derivative. Then $r_0(\beta)$ is the minimal positive solution to 
$$\sqrt{\frac{\beta^2}{4}-r_0(\beta)}\,D_{r_0(\beta)}(-\beta)=D'_{r_0(\beta)}(-\beta),\quad 0< r_0(\beta)< \frac{\beta^2}{4},$$
which can be shown to exist for $\beta\ge \beta_*$. In \cite{gam}, $r_0(\beta)$ was characterized in a somewhat different form. If we let $\widehat{p}(x;\theta)=\int_0^{\infty}e^{-\theta t}p(x,t)\,dt$ be the Laplace transform of the density, then $\theta=-\beta^2/4$ corresponds to a branch point of $\widehat{p}$ and $\theta=-r_0(\beta)$ is a pole of $\widehat{p}$. In \cite{lee} we examined further the spectral properties of the Halfin-Whitt diffusion, showing that there is always a continuous spectrum, corresponding to $\Re(\theta)\le -\beta^2/4$, and that the point spectrum is empty for $\beta<0$, and may contain any number of points for $\beta>0$. For $0<\beta<\beta_*$, $\theta=0$ is the only pole of $\widehat{p}$, but as $\beta$ increases, other poles appear in the range where $\theta$ is real with $\theta\in (-\beta^2/4,\,0)$. In fact an additional pole appears precisely when $\beta$ increases past a root of $D'_{\beta^2/4}(-\beta)$, and there are infinitely many of these. In the limit $\beta\to +\infty$ the point spectrum approaches $-\theta=\{0,1,2,\cdots\}$, which is that of the standard free space Ornstein-Uhlenbeck process, where the linear drift $-X_d-\beta$ would apply for all $X_d\in(-\infty,\infty)$. 

The relaxation rate asymptotics in \cite{lee}, \cite{gam} assume that $x$, $x_0$ and $\beta$ are fixed, and that only $t\to\infty$. Here we shall consider another asymptotic limit (see (\ref{2_limit})) where $\beta$ is assumed large and positive, and $(x,x_0)$ are also scaled to be $O(\beta)$. Then $x_0>0$ and $x_0=O(\beta)$ means that the process $X_d(t)$ starts at a large positive value, but as $t$ increases it will most likely be in the range $X_d(t)<0$, due to the strong negative drift. For very large times we expect that $X_d(t)\approx -\beta$, since this is the equilibrium point in the range $X_d<0$. We shall obtain a variety of asymptotic results for $p(x,t)$ in this limit, which will show how the migration of mass from $X_d>0$ to $X_d<0$ occurs. The present asymptotic results are very different from the relaxation rate results, as for large $\beta$ the continuous spectrum will not play much of a role, and there will be many points in the point spectrum, which will correspond to poles of $\widehat{p}$ that are close to $-\theta=0,1,2,\cdots$. Getting an accurate estimate of these poles for $\beta\to+\infty$ will be crucial to the asymptotic analysis. 

We mention that the Halfin-Whitt diffusion limit is currently of interest in the modeling of call centers (see \cite{bor}, \cite{gan}, \cite{jan}). Related limit theorems for other more general and more complicated models appear in \cite{man}, \cite{gar}, \cite{jel}, \cite{ree}, and the mean hitting times for the Halfin-Whitt diffusion were analyzed in \cite{mag}. 

Our analysis is based on complex variable techniques, the asymptotic evaluation of integrals and sums, and properties of the parabolic cylinder functions $D_p(z)$. Good references on asymptotic methods for evaluating integrals/sums, such as the Laplace method, saddle point method, Euler-Maclaurin formula, singularity analysis, etc., may be found in the books of Wong \cite{won}, Bleistein and Handelsman \cite{ble}, Bender and Orszag \cite{ben}, Flajolet and Sedgewick \cite{fla}, and Szpankowski \cite{szp}. Basic properties of parabolic cylinder functions appear in the books of Abramowitz and Stegun \cite{abr}, Gradshteyn and Ryzhik \cite{gra}, Magnus, Oberhettinger and Soni \cite{magnus} and Temme \cite{tem}. 

The paper is organized as follows. In section 2 we summarize our main results (see Theorems \ref{theorem2} and \ref{theorem3}), and interpret them probabilistically. In section 3 we derive the results for the range $x>0$, and in section 4 for $x<0$. Finally in section 5 we re-derive the results using an alternate approach, which is based on a singular perturbation analysis of the PDE satisfied by $p(x,t)$. This more geometric approach is then compared and contrasted to the asymptotic evaluation of sums and integrals. As the analysis is fairly complicated, we include many graphs to illustrate specific points.

\section{Summary of results}

We let $p(x,t)$ be the transient probability density of the Halfin-Whitt diffusion $X_d(t)$, subject to the initial condition $X_d(0)=x_0$. The density satisfies the Kolmogorov forward equation
\begin{equation}\label{2_pde}
p_t=\left\{ \begin{array}{ll}
\beta p_x+p_{xx}, & x>0\\
\beta p_x+(xp)_x+p_{xx}, & x<0
\end{array} \right.
\end{equation}
with the initial condition $p(x,0)=\delta(x-x_0)$ and the ``interface" condition $p(0^-,t)=p(0^+,t)$ and $p_x(0^-,t)=p_x(0^+,t)$. This diffusion arises naturally in the multi-server $M/M/m$ queue, in the limit where $m\to\infty$ and traffic intensity $\rho=\lambda/(\mu\,m)\to 1$, with the scaling $\rho=1-\beta/\sqrt{m}=1+O(m^{-1/2})$. Here $\lambda$ and $\mu$ are the arrival and service rates in the discrete model. Note that the discrete model is stable for $\rho<1$ which corresponds to $\beta>0$. The queue length $N(t)$ in the $M/M/m$ queue can then be approximated by $N(t)\approx m+\sqrt{m} X_d(t)$. 

When $\beta>0$ the steady state distribution of the diffusion can be easily obtained as
\begin{equation}\label{2_steady}
p(x,\infty)=C(\beta)\left\{ \begin{array}{ll}
\beta e^{-\beta x}, & x>0\\
\beta e^{-\beta x}e^{-x^2/2}, &x<0
\end{array} \right.
\end{equation}
where the normalization constant is $1/C(\beta)=1+\beta e^{\beta^2/2}\int_{-\infty}^\beta e^{-u^2/2}du.$
In the limit of $\beta\to +\infty$ we have $1/C(\beta)\sim \sqrt{2\pi}\,\beta e^{\beta^2/2}$ and then 
\begin{equation}\label{2_steady_appro}
p(x,\infty)\sim \frac{1}{\sqrt{2\pi}}
\left\{ \begin{array}{ll}
e^{-\beta x}e^{-\beta^2/2}, & x>0\\
e^{-(x+\beta)^2/2}, &x<0
\end{array} \right.
\end{equation}
which is concentrated near $x=-\beta<0$. 

Now assume that the initial condition $x_0>0$. In terms of the discrete model this means that the initial number of customers $N(0)$ exceeds the number of servers $m$, so that some customers are in the queue. In \cite{lee} we obtained an explicit expression for the Laplace transform $\widehat{p}(x,\theta)=\int_0^{\infty}e^{-\theta t}p(x,t)\,dt$ of the transient density, which we summarize below.

\begin{theorem}\label{theorem1}
If $x_0>0$ we have the following.
\begin{enumerate}
\item For $x>0$:
	\begin{equation}\label{2_th1_1}
	\widehat{p}(x,\theta)=\frac{e^{\beta(x_0-x)/2}}{\sqrt{\beta^2+4\theta}}\bigg[e^{-|x-x_0|\sqrt{\theta+\beta^2/4}}-e^{-(x+x_0)\sqrt{\theta+\beta^2/4}}\bigg]+\frac{e^{\beta(x_0-x)/2}e^{-(x+x_0)\sqrt{\theta+\beta^2/4}}}{\sqrt{\theta+\beta^2/4}-R_\beta(\theta)}
	\end{equation}
with 
$$R_\beta(\theta)=\frac{D'_{-\theta}(-\beta)}{D_{-\theta}(-\beta)}=-\frac{d}{d\beta}\log\big[D_{-\theta}(-\beta)\big]$$
where $D_{-\theta}(-\beta)$ is the parabolic cylinder function.

\item For $x<0$:
\begin{equation}\label{2_th1_2}
\widehat{p}(x,\theta)=e^{-x^2/4}e^{\beta(x_0-x)/2}\frac{e^{-x_0\sqrt{\theta+\beta^2/4}}}{\sqrt{\theta+\beta^2/4}-R_\beta(\theta)}\frac{D_{-\theta}(-\beta-x)}{D_{-\theta}(-\beta)}.
\end{equation}
\end{enumerate}
\end{theorem}
When $x=0$, either (\ref{2_th1_1}) or (\ref{2_th1_2}) may be used to compute $\widehat{p}(0,\theta)$, as $\widehat{p}(0^+,\theta)=\widehat{p}(0^-,\theta)$. By inverting the Laplace transform we obtain integral representations for $p(x,t)$ in the form
\begin{equation}\label{2_invlap}
p(x,t)=\frac{1}{2\pi i}\int_{Br}e^{\theta t}\,\widehat{p}(x,\theta)\,d\theta,
\end{equation}
where $Br$ is a vertical Bromwich contour in the complex $\theta$-plane, on which $\Re(\theta)>0$. 

In \cite{lee} we studied some spectral properties of this problem. The expressions in (\ref{2_th1_1}) and (\ref{2_th1_2}) have a branch point at $\theta=-\beta^2/4$, and there is always a continuous spectrum in the range $\Re(\theta)\in(-\infty, -\beta^2/4)$. For $\beta<0$ the spectrum is purely continuous, but for $\beta>0$ there may be any number of discrete eigenvalues. For any $\beta>0$, $\theta=0$ is a simple pole of (\ref{2_th1_1}) and (\ref{2_th1_2}) and the residue at this pole corresponds to the steady state solution in (\ref{2_steady}). As $\beta$ increases past each root of $D'_{\beta^2/4}(-\beta)=0$ (there are infinitely many such roots in the range $\beta>0$), a new discrete eigenvalue enters the spectrum.

The first root occurs when $\beta=\beta_*\approx 1.85722\cdots$ and thus for $0<\beta<\beta_*$, $\theta=0$ is the only pole of (\ref{2_th1_1}) and (\ref{2_th1_2}), but for $\beta>\beta_*$ there exists at least one pole on the real axis with $\Re(\theta)\in(-\beta^2/4,0)$. In \cite{lee} we also studied how as $\beta\to +\infty$ the spectrum approaches that of the standard free space Ornstein-Uhlenbeck process, whose Laplace transformed density would have poles at $\theta=0,-1,-2,\cdots$. Understanding the spectrum is useful in estimating the approach to equilibrium, i.e., to obtain $p(x,t)-p(x,\infty)$ for $t\to\infty$ and fixed $x$ and $x_0$. 

The expressions in (\ref{2_th1_1})-(\ref{2_invlap}) do not yield much immediate insight into the qualitative structure of the density $p(x,t)$. Thus it is useful to evaluate these expressions in various asymptotic limits, to gain more insight, for example, as to how equilibrium is achieved when $\beta>0$. Here we shall consider the limit
\begin{equation}\label{2_limit}
\beta\to +\infty, \quad x=\beta X,\quad x_0=\beta X_0,\quad X_0>0.
\end{equation}
Thus we are assuming that the drift parameter $\beta$ is large and positive, and then scale $(x,x_0)$ by $\beta$. Furthermore, $X_0>0$ means that initially the number of customers $N(t)$ well exceeds the number of servers. But for $t\to\infty$ the density in the steady state becomes sharply concentrated near $X=-1$, as the second formula in (\ref{2_steady_appro}) becomes 
$p(x,\infty)\sim e^{-\beta^2(X+1)^2/2}/\sqrt{2\pi}$, for $X<0$.
Then $N(t)$ tends to be less than the number of servers $m$, and there is no need for customers to queue. Thus our analysis will provide some insight into how the $M/M/m$ model tends to process large queues, at least in this diffusion limit. Note that for $\beta\to +\infty$ there are many discrete eigenvalues, and in fact our analysis will precisely estimate these. However, our asymptotic results are much different than the relaxation time asymptotics that we previously obtained in \cite{lee}, which assume that $t\to\infty$ but with fixed $\beta$, $x$, $x_0$. 

To interpret some of our results it is useful to consider a fluid approximation to the diffusion process $X_d(t)$, where we neglect completely the ``diffusive" or variability effects. For $\beta\to +\infty$ we have $X_d(t)\approx \beta Y_{fl}(t)$ where, since $x_0=\beta X_0$ with $X_0>0$,
\begin{equation}\label{2_Yfl}
Y_{fl}(t)=\left\{ \begin{array}{ll}
X_0-t, & 0\le t\le X_0\\
e^{X_0-t}-1, & t\ge X_0.
\end{array} \right.
\end{equation}
Thus the fluid approximation starts at $X_0$, leaves the range $Y_{fl}>0$ when $t$ increases past $X_0$, and then remains in the range $Y_{fl}<0$ for $t>X_0$, ultimately approaching $Y_{fl}(\infty)=-1$. Then (\ref{2_steady_appro}) gives the Gaussian spread about this limiting value.

The fluid approximation indicates where most of the probability mass concentrates, but it does not, for example, predict how some small amount of probability mass would remain in the range $X>0$, which would lead to the first formula in (\ref{2_steady_appro}). Our goal is to obtain a more ``global" description of the density $p(x,t)=p(x,t;x_0,\beta)$ and to show how it behaves in various space-time ranges.

We introduce the small parameter $\varepsilon$ with 
$\varepsilon=\beta^{-2}$, so that $\beta=1/\sqrt{\varepsilon}$ and
$\varepsilon\to 0^+$, and summarize our main asymptotic results in the $(X,t)$ plane, first for $X>0$ and then $X<0$. 

\begin{theorem}\label{theorem2}
For $X(=x/\beta=x\sqrt{\varepsilon})>0$ we have the following asymptotic results:
\begin{enumerate}
	\item[(i)] $0<t<t_1(X;X_0)$
	\begin{equation}\label{2_th2_1}
	 p(x,t)-\frac{1}{2\sqrt{\pi t}}\exp\Big[-\frac{(t+X-X_0)^2}{4\varepsilon t}\Big]\sim  \frac{\varepsilon}{4\sqrt{\pi}}\frac{t^{3/2}(X+X_0+t)}{(X+X_0)^3}e^{-X/\varepsilon}\exp\Big[-\frac{(t-X-X_0)^2}{4\varepsilon t}\Big].
	 \end{equation}
\item[(ii)] $t_1(X;X_0)<t<t_2(X;X_0)$
	\begin{equation}\label{2_th2_2}
	 p(x,t)-\frac{1}{2\sqrt{\pi t}}\exp\Big[-\frac{(t+X-X_0)^2}{4\varepsilon t}\Big]\sim   \frac{1}{2\sqrt{2\pi}}\frac{(1-4z_*)^{1/4}\big[1+\sqrt{1-4z_*}\big]}{\sqrt{1-4z_*-2(X+X_0)z_*}}e^{f(z_*)/\varepsilon},
	 \end{equation}
where
\begin{equation}\label{2_th2_2f}
f(z_*)=-\frac{1}{2}(X+X_0+1)\sqrt{1-4z_*}+\frac{X_0-X}{2}-\frac{X+X_0}{\sqrt{1-4z_*}}z_*
\end{equation}
and $z_*=z_*(X,t)$ is defined implicitly, as the solution to 
\begin{equation}\label{2_th2_2z*}
t=\frac{X+X_0}{\sqrt{1-4z_*}}+2\log\bigg(\frac{1+\sqrt{1-4z_*}}{2\sqrt{z_*}}\bigg).
\end{equation}
Furthermore, $t_2$ is implicitly defined as the solution to 
\begin{equation}\label{2_th2_2t2}
t=\frac{(X-X_0)^2}{t}-\sqrt{\mathcal{A}}+2\log\bigg[\frac{2(X+X_0+2)+t+\frac{(X-X_0)^2}{t}+\sqrt{\mathcal{A}}}{2(X+X_0+2)-t-\frac{(X-X_0)^2}{t}-\sqrt{\mathcal{A}}}\bigg],
\end{equation}
$$\mathcal{A}=\bigg[t+\frac{(X-X_0)^2}{t}\bigg]^2-4(X+X_0)(X+X_0+2),$$
and $t_1$ is implicitly defined as the solution to 
\begin{equation}\label{2_th2_2t1}
t=\frac{(X+X_0)^2}{t}-\sqrt{\mathcal{B}}+2\log\bigg[\frac{2(X+X_0+2)+t+\frac{(X+X_0)^2}{t}+\sqrt{\mathcal{B}}}{2(X+X_0+2)-t-\frac{(X+X_0)^2}{t}-\sqrt{\mathcal{B}}}\bigg],
\end{equation}
$$\mathcal{B}=\bigg[t+\frac{(X+X_0)^2}{t}\bigg]^2-4(X+X_0)(X+X_0+2).$$
\item[(iii)] $t>t_2(X;X_0)$
\begin{equation}\label{2_th2_3}
p(x,t)\sim \frac{1}{2\sqrt{2\pi}}\frac{(1-4z_*)^{1/4}\big[1+\sqrt{1-4z_*}\big]}{\sqrt{1-4z_*-2(X+X_0)z_*}}e^{f(z_*)/\varepsilon},
\end{equation}
where $z_*$ and $f$ are as in (\ref{2_th2_2z*}) and (\ref{2_th2_2f}). 
\end{enumerate}
\end{theorem}
On the time scale $t=\log(1/\varepsilon)+O(1)$, with 
\begin{equation}\label{2_Tscale}
t=\log(1/\varepsilon)+T=2\log\beta+T,\quad T=O(1),
\end{equation}
we have the explicit expression
\begin{equation}\label{2_pT}
p(x,t)\sim\frac{1}{\sqrt{2\pi}}e^{-(X+1/2)/\varepsilon}\exp\big[e^{-(T-X-X_0)}\big].
\end{equation}

To interpret these results we note that the transient density for a free space Brownian motion with unit negative drift is given by
\begin{equation}\label{2_pbm}
p_{_{BM}}(X,t;X_0)=\frac{1}{2\sqrt{\pi t}}\exp\bigg[-\frac{(t+X-X_0)^2}{4\varepsilon t}\bigg].
\end{equation}
Here the process starts at $X_0$ when $t=0$ and has $\varepsilon$ as its diffusion coefficient.

If the first PDE in (\ref{2_pde}) were to apply for all $x$, rather than only $x>0$, then $p(x,t)$ would be given exactly by (\ref{2_pbm}). Theorem \ref{theorem2} shows that there are three ranges in the $(X,t)$ plane where the behavior of $p(x,t)$ is different. For $t<t_2(X;X_0)$ items (i) and (ii) (with (\ref{2_th2_1}) and (\ref{2_th2_2})) show that $p\sim p_{_{BM}}$ with an exponentially small error. The error term undergoes a ``phase transition" along the curve $t=t_1(X;X_0)$, and we gave the results in Theorem \ref{theorem2} so as to estimate the difference $p-p_{_{BM}}$ for those space-time ranges where $p\sim p_{_{BM}}$. We also note that along the first transition curve $t_1$ we have
$-X-(t-X-X_0)^2/(4t)=f(z_*)$
and along the second curve $t_2$ we have
\begin{equation}\label{2_t2eqn}
-\frac{1}{4t}(t+X-X_0)^2=f(z_*).
\end{equation}
Thus along the curve $t=t_1$ the right sides of (\ref{2_th2_1}) and (\ref{2_th2_2}) become roughly comparable, while along $t=t_2$ the right side of (\ref{2_th2_2}) becomes roughly comparable to $p_{_{BM}}$. In the range $t>t_2$, $p_{_{BM}}$ can no longer be used to approximate the density and we obtain the complicated expression in (\ref{2_th2_3}). From (\ref{2_th2_2z*}) we can show that 
\begin{equation}\label{2_z*sim}
z_*(X,t)\sim e^{-t}e^{X+X_0},\quad t\to\infty
\end{equation}
and then
\begin{equation}\label{2_fz*sim}
f(z_*)\sim -\frac{1}{2}-X+e^{-t}e^{X+X_0},\quad t\to\infty.
\end{equation}
Thus if we scale time as in (\ref{2_Tscale}) then on the $T=O(1)$ time scale $f(z_*)+1/2+X$ becomes $O(\varepsilon)$ and (\ref{2_th2_3}) reduces to (\ref{2_pT}) as a special case. The latter expression shows how the steady state in (\ref{2_steady_appro}) for $x>0$ is slowly achieved on the $T$-scale, as $p(x,t)$ decreases toward $p(x,\infty)$ as a double exponential function in $T$.

In section 5 we give more probabilistic and geometric interpretations of the three main results in Theorem \ref{theorem2}. Roughly speaking, the approximation $p\approx p_{_{BM}}$ corresponds to sample paths that go directly from $(X_0,0)$ to $(X,t)$, the right side of (\ref{2_th2_1}) corresponds to sample paths that go from $(X_0,0)$ to $X=0$ and then ``reflect" at the interface before returning to $(X,t)$, and the right side of (\ref{2_th2_3}) corresponds to paths that go from $(X_0,0)$, hit $X=0$, spend some time in the range $X<0$, and then return to $X=0$ and re-enter the range $X>0$. Note that the exponential part of the right side of (\ref{2_th2_1}) would be consistent with a reflection (in $X=0$) of the Brownian motion in (\ref{2_pbm}), but the algebraic factors indicate a more complicated reflection/transmission law, through the interface $X=0$. 

When $t\approx t_1$ the right sides of (\ref{2_th2_1}) and (\ref{2_th2_2}) become comparable, and then $p-p_{_{BM}}$ can be estimated by adding the two expressions. Similarly, when $t\approx t_2$ we can approximate $p$ by $p_{_{BM}}$ plus the right side of (\ref{2_th2_3}). Finally, we note that Theorem \ref{theorem2} applies not only for $\beta\to +\infty$ with $X=x/\beta>0$, but also on the original $x$ scale, for $x>0$ (then $X=O(1/\beta)$). In fact it can be used to compute $p(0,t)$ asymptotically by letting $X=0$. 

In Figure \ref{figure1} we plot the transition curves $t_1$ and $t_2$ for the three values $X_0=0.2,\, 1,\, 5$ to illustrate their dependence on the initial condition. Note that when $X=0$, $t_1(0;X_0)=t_2(0;X_0)$, as then (\ref{2_th2_2t2}) and (\ref{2_th2_2t1}) coincide.

We next consider the range $X<0$, where the asymptotics of $p(x,t)$ become much more complicated. It proves useful to divide $\{X<0\}$ as $\{X\le -1\}\cup\{-1<X<0\}$, and the strip $-1<X<0$ in the $(X,t)$ plane must be further subdivided into 7 regions. We first define
\begin{equation}\label{2_t*}
t_*=t_*(X;X_0)=\frac{X_0}{\sqrt{-2X-X^2}}+\log\bigg[\frac{1+X}{1-\sqrt{-2X-X^2}}\bigg],\quad -1<X<0,
\end{equation}
and note that $-2X-X^2=1-(1+X)^2\in (0,1)$ when $X\in (-1,0)$. We also have $t_*\to +\infty$ as either $X\to -1^+$ or $X\to 0^-$. More precisely, 
$t_*\sim{X_0}/{\sqrt{-2X}}$ as $X\to 0^-$
and
$t_*=-\log(1+X)+\log 2+X_0+o(1)$ as $X\to (-1)^+$.
We furthermore define the curves $t_c$ and $t_d$ by
\begin{equation}\label{2_tc}
t_c=t_c(X;X_0)=\alpha_c+\log\bigg[\frac{X+1+\sqrt{(X+1)^2+{X_0^2}/{\alpha_c^2}-1}}{1-{X_0}/{\alpha_c}}\bigg],
\end{equation}
where
\begin{equation}\label{2_alpha_c}
\alpha_c =\Big[-X(X_0+2)(X+2)+X_0+\sqrt{X_0}(X+1)\sqrt{(X^2+2X)(X_0+4)+X_0}\Big]^{1/2}\sqrt{\frac{X_0}{-2X(X+2)}},
\end{equation}
and
\begin{equation}\label{2_td}
t_d=t_d(X;X_0)=\alpha_d+\log\bigg[\frac{X+1+\sqrt{(X+1)^2+{X_0^2}/{\alpha_d^2}-1}}{1-{X_0}/{\alpha_d}}\bigg],
\end{equation}
where
\begin{equation}\label{2_alpha_d}
\alpha_d  =\Big[-X(X_0+2)(X+2)+X_0-\sqrt{X_0}(X+1)\sqrt{(X^2+2X)(X_0+4)+X_0}\Big]^{1/2}\sqrt{\frac{X_0}{-2X(X+2)}}.
\end{equation}
The curves $t_c$ and $t_d$ exist only in the range $X\in (X_{cusp},0)$ where
\begin{equation}\label{2_Xcusp}
X_{cusp}=-1+\frac{2}{\sqrt{X_0+4}}.
\end{equation}
Note that $\big[(X_{cusp}+1)^2-1\big](X_0+4)+X_0=0$. The curves $t_c$ and $t_d$ thus coincide at the point $(X_{cusp},t_{cusp})$, where
\begin{equation}\label{2_tcusp}
t_{cusp}=\sqrt{X_0(X_0+3)}+\log\bigg[\frac{X_0+2+\sqrt{X_0(X_0+3)}}{\sqrt{X_0+4}}\bigg].
\end{equation}
When $X=X_{cusp}$ we have, in view of (\ref{2_tcusp}) and (\ref{2_t*}), $t_{cusp}>t_*$. From (\ref{2_tc}) and (\ref{2_alpha_c}) we see that $t_c\to \infty$ as $X\to 0^-$, and comparing the small $X$ behaviors of $t_c$ and $t_*$ we can show that $t_c>t_*$ for $X\to 0^-$. On the other hand $t_d$ tends to the finite limit $t_d(0;X_0)$ as $X\to 0^-$, with
\begin{equation*}\label{2_td_0}
t_d(0;X_0)=\sqrt{X_0(X_0+2)}+\log\bigg[\frac{\sqrt{X_0+2}+\sqrt{X_0}}{\sqrt{X_0+2}-\sqrt{X_0}}\bigg],
\end{equation*}
and also $\alpha_d\to \alpha_d(0;X_0)=\sqrt{X_0(X_0+2)}$. We use the term ``cusp" to denote the intersection of $t_c$ and $t_d$ since these curves have the common slope
${dt}/{dX}=\sqrt{{(X_0+3)(X_0+4)}/{X_0}}$
at this point. We can also view the union of $t_c$ and $t_d$ as representing a single curve with a cusp point at $(X_{cusp},t_{cusp})$. In section 5 we give more geometric interpretation to $t_c$ and $t_d$, since they correspond to the envelope(s) of a certain family of curves, and thus $t_c$ and $t_d$ are called ``caustic" curves.

\begin{figure}[h]   
  \begin{minipage}[t]{0.5\linewidth}  
    \centering   
    \includegraphics[width=0.75\textwidth]{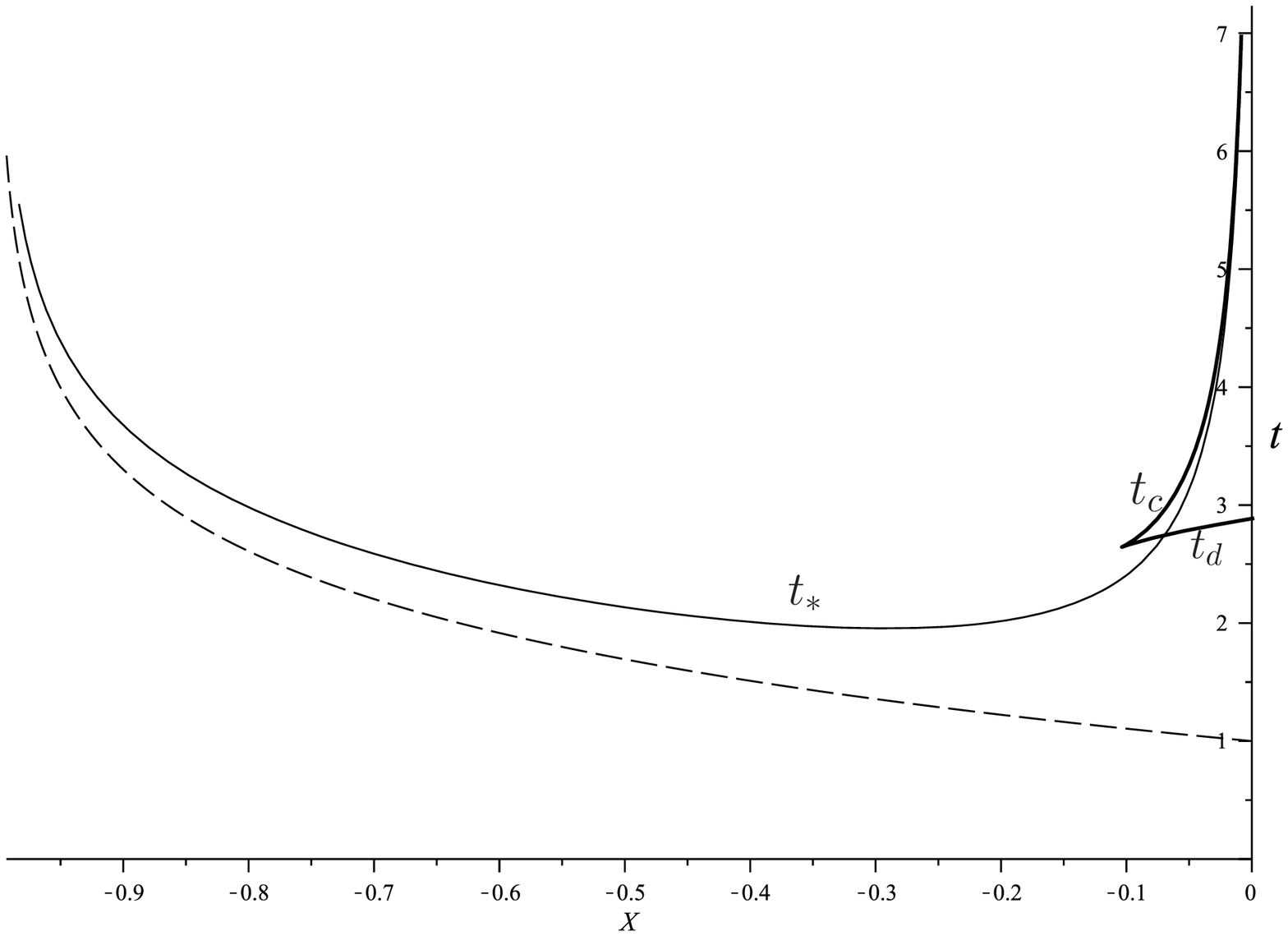}   
    \caption[]{The curves $t_*$, $t_c$, $t_d$ and the fluid\\approximation(the dashed curve) for $X_0=1$. }   
    \label{figure2}   
  \end{minipage}%
  \begin{minipage}[t]{0.5\linewidth}   
    \centering   
    \includegraphics[width=0.75\textwidth]{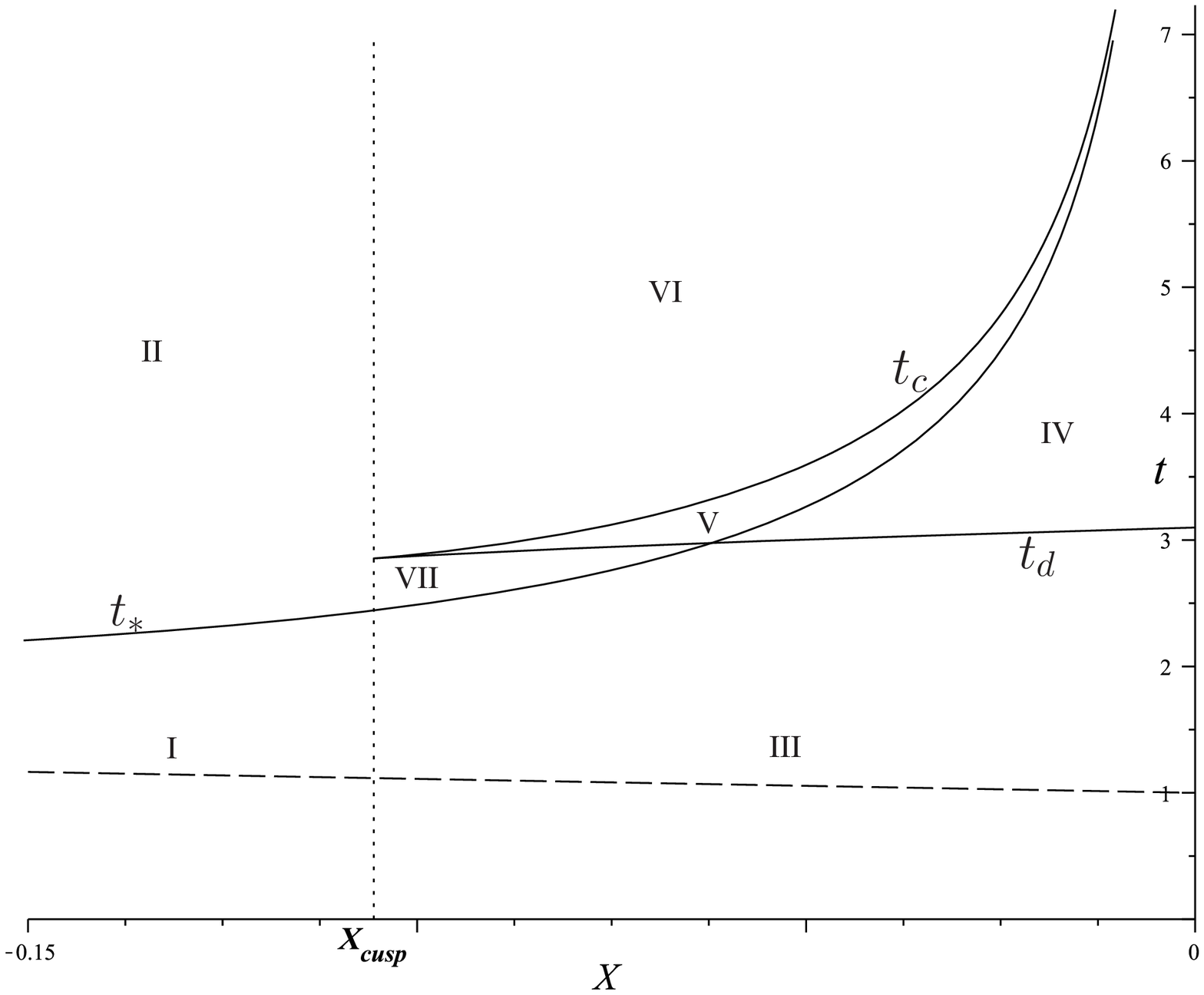}   
    \caption{The curves $t_*$, $t_c$ and $t_d$ split the strip $-1<X<0$ into seven regions (the dashed curve is the fluid approximation).}   
    \label{figure3}   
  \end{minipage}   
\end{figure}

In Figure \ref{figure2} we sketch the curves $t_*$, $t_c$ and $t_d$ in the range $X\in(-1,0)$ for $X_0=1$. This shows that $t_c$ (which exists only for $X>X_{cusp}$) lies above $t_*$, while $t_d$ intersects $t_*$ at a point, say $(X_*, t_{**})$, which may be obtained numerically by simultaneously solving (\ref{2_t*}) and (\ref{2_td}) with (\ref{2_alpha_d}). We also sketch (the dashed line) the fluid approximation in (\ref{2_Yfl}), which is the curve $X=e^{X_0-t}-1\in (-1,0)$, or $t=X_0-\log(1+X)$, and this lies below each of the curves $t_*$, $t_c$, $t_d$ (when the latter two exist). While the figure fixes $X_0=1$, the pattern of intersection is the same for any $X_0>0$. Hence the three curves $t_*$, $t_c$, $t_d$ split the strip $-1<X<0$ into the following 7 regions:
\begin{eqnarray}
&&\textrm{Region  I}=\{X\le -1\}\cup\{0<t<t_*,\;-1<X<X_{cusp}\}\label{2_region1}\\
&&\textrm{Region  II}=\{t>t_*,\;-1<X<X_{cusp}\}\label{2_region2}\\
&&\textrm{Region  III}=\{0<t<\min\{t_*,t_d\},\;X_{cusp}<X<0\}\label{2_region3}\\
&&\textrm{Region  IV}=\{t_d<t<t_*,\;X_*<X<0\}\label{2_region4}\\
&&\textrm{Region  V}=\{\max\{t_*,t_d\}<t<t_c,\;X_{cusp}<X<0\}\label{2_region5}\\
&&\textrm{Region  VI}=\{t>t_c,\;X_{cusp}<X<0\}\label{2_region6}\\
&&\textrm{Region  VII}=\{t_*<t<t_d,\;X_{cusp}<X<X_*\}.\label{2_region7}
\end{eqnarray}
It will prove convenient to include $\{X\le -1\}$ within Region I. Note that $\min\{t_*,t_d\}=t_*$ if $X<X_*$ and $\min\{t_*,t_d\}=t_d$ if $X>X_*$. The different regions are indicated in Figure \ref{figure3}. The fluid approximation lies within Region I if $X<X_{cusp}$ and within Region III if $X\in(X_{cusp},0)$. The different regions will correspond to different asymptotic formulas for $p(x,t)$. 

Our results will involve the solutions of certain transcendental equations, similar to the one in (\ref{2_th2_2z*}), so we state two lemmas about these equations. 

\begin{lemma}\label{lemma1}
For $t>0$ and $X<0$ the equation
\begin{equation}\label{2_lm1}
t-\frac{X_0}{\sqrt{1+4\phi_s}}+\log\bigg[\frac{\sqrt{1+4\phi_s}-1}{\sqrt{(1+X)^2+4\phi_s}-1-X}\bigg]=0,\quad \phi_s>-\frac{1}{4}(1+X)^2
\end{equation}
has a unique solution $\phi_s=\phi_s(X,t)$ for $(X,t)$ in Regions I, III and IV. Also, $\phi_s>0$ if $X\le -1$ or $X\in(-1,0)$ and $t<X_0-\log(1+X)$ (i.e., below the fluid approximation). When $t=X_0-\log(1+X)$ we have $\phi_s=0$. When $X_0-\log(1+X)<t<t_*(X;X_0)$ we have $\phi_s<0$, and $\phi_s\to-(1+X)^2/4$ as $t\uparrow t_*$. In Regions II and V-VII, (\ref{2_lm1}) has no solutions in the range $\phi_s>-(1+X)^2/4$.
\end{lemma} 
\begin{lemma}\label{lemma2}
For $t>0$ and $-1<X<0$ the equation
\begin{equation}\label{2_lm2}
t-\frac{X_0}{\sqrt{1-4z}}+\log\bigg[\frac{1-\sqrt{1-4z}}{1+X+\sqrt{(1+X)^2-4z}}\bigg]=0,\quad z>0
\end{equation}
has a unique solution for $(X,t)$ in Regions II and VI which we denote by $z_1=z_1(X,t)$. In Region V, (\ref{2_lm2}) has three solutions that we order as $z_1<z_2<z_3$. In Region VII there is again a unique solution that we now denote by $z_3=z_3(X,t)$. In Region IV there are two solutions that we order as $z_1<z_2$, and they lie in the range $z\in (0,-\phi_s)$, where $\phi_s$ is the solution to (\ref{2_lm1}). For Regions II and V-VII, these conclusions consider (\ref{2_lm2}) over the ``full" range $z\in(0,(1+X)^2/4)$. For Regions I and III, (\ref{2_lm2}) has no solutions in the range $z\in (0,-\phi_s)$ if $\phi_s<0$.
\end{lemma}

Note that if we consider the totality of the solutions to (\ref{2_lm1}) and (\ref{2_lm2}), this number is 1 for Regions I-III, VI, and VII, and 3 for Regions IV and V. Thus there are a total of 3 solutions precisely when $t_d<t<t_c$ (for $X\in(X_{cusp},0)$). We can now state our main results for $X<0$. 

\begin{theorem}\label{theorem3}
For $X(=x/\beta=x\sqrt{\varepsilon})<0$ we have the following asymptotic results:
\begin{enumerate}
\item[(i)] $(X,t)$ in Regions I and III:
\begin{equation}\label{2_th3_1}
p(x,t)\sim G(X,t)\exp\Big[\frac{1}{\varepsilon}F(X,t)\Big],
\end{equation}
\begin{equation}\label{2_th3_1F}
F(X,t)=-\frac{X^2}{4}-\frac{X}{2}+\frac{X_0}{2}+\frac{X_0\phi_s}{\sqrt{1+4\phi_s}}-\Big(\frac{X_0}{2}+\frac{1}{4}\Big)\sqrt{1+4\phi_s}+\frac{1+X}{4}\sqrt{(1+X)^2+4\phi_s},
\end{equation}
\begin{eqnarray}\label{2_th3_1G}
G(X,t)&=&\frac{1}{\sqrt{2\pi}}\frac{(1+4\phi_s)^{-1/4}}{\big[(1+X)^2+4\phi_s\big]^{1/4}}\bigg[\frac{\sqrt{1+4\phi_s}+1}{\sqrt{(1+X)^2+4\phi_s}+1+X}\bigg]^{1/2}\nonumber\\
&&\times\bigg|\frac{2X_0}{(1+4\phi_s)^{3/2}}+\frac{1}{2\phi_s}\bigg(\frac{1}{\sqrt{1+4\phi_s}}-\frac{1+X}{\sqrt{(1+X)^2+4\phi_s}}\bigg)\bigg|^{-1/2},
\end{eqnarray}
and $\phi_s=\phi_s(X,t)$ is the unique solution to (\ref{2_lm1}). For $(X,t)$ close to the fluid approximation, we let 
$X=-1+e^{X_0-t}+\sqrt{\varepsilon}\Delta$.
Then for $\varepsilon\to 0$ and $\Delta$ fixed, (\ref{2_th3_1}) simplifies to the Gaussian approximation
\begin{equation}\label{2_th3_1gau}
p(x,t)\sim\frac{1}{\sqrt{2\pi}}\frac{e^{t-X_0}}{\sqrt{2X_0+e^{2(t-X_0)}-1}}\exp\bigg[-\frac{e^{2(t-X_0)}}{e^{2(t-X_0)}-1+2X_0}\frac{\Delta^2}{2}\bigg].
\end{equation}
For $t$ large, with $t=\log(1/\varepsilon)+T$ as in (\ref{2_Tscale}), and $X<-1$, (\ref{2_th3_1}) simplifies to 
\begin{equation}\label{2_th3_1T}
p(x,t)\sim\frac{1}{\sqrt{2\pi}}\exp\bigg[-\frac{(1+X)^2}{2\varepsilon}\bigg]\exp\bigg[e^{X_0-T}(1+X)\bigg],
\end{equation}
which gives the approach to equilibrium as $T\to\infty$. For $t\approx t_*$ (the upper limit for Region I if $X\in(-1,X_{cusp})$ and Region III if $X\in(X_{cusp},X_*)$), (\ref{2_th3_1}) simplifies to 
\begin{equation}\label{2_th3_1t*}
p(x,t)\sim\frac{1}{2\sqrt{\pi}}\frac{\sqrt{1+\sqrt{1-(1+X)^2}}}{\big[1-(1+X)^2\big]^{1/4}}\exp\bigg[\frac{F_c(X)}{\varepsilon}-\frac{(1+X)^2}{4\,\varepsilon^{2/3}}\Lambda-\frac{(1+X)^2}{12}\Lambda^3\bigg]
\end{equation}
where
\begin{equation}\label{2_th3_1lambda}
t=t_*+\varepsilon^{1/3}\Lambda,\quad\Lambda=O(1)
\end{equation}
and
\begin{equation}\label{2_th3_1Fc}
F_c(X)=-\frac{X^2}{4}-\frac{X}{2}+\frac{X_0}{2}-\frac{X_0}{4\sqrt{-X^2-2X}}-\Big(\frac{X_0}{4}+\frac{1}{4}\Big)\sqrt{-X^2-2X}.
\end{equation}
\item[(ii)] $(X,t)$ in Regions II and VI:
\begin{eqnarray}\label{2_th3_2}
p(x,t)\sim \widetilde{G}(z_1)\exp\Big[\frac{1}{\varepsilon}\widetilde{f}(z_1)\Big],
\end{eqnarray}
\begin{equation*}
\widetilde{f}(z)=\widetilde{f}(z;X,t)=-\frac{X^2}{4}-\frac{X}{2}+\frac{X_0}{2}-\frac{zX_0}{\sqrt{1-4z}}-\Big(\frac{X_0}{2}+\frac{1}{4}\Big)\sqrt{1-4z}-\frac{1+X}{4}\sqrt{(1+X)^2-4z},
\end{equation*}
\begin{eqnarray*}
\widetilde{G}(z)&=&\widetilde{G}(z;X,t)=\frac{1}{\sqrt{2\pi}}\frac{(1-4z)^{-1/4}}{\big[(1+X)^2-4z\big]^{1/4}}\bigg[\frac{1+\sqrt{1-4z}}{1+X-\sqrt{(1+X)^2-4z}}\bigg]^{1/2}\nonumber\\
&&\times\bigg|-\frac{2X_0}{(1-4z)^{3/2}}+\frac{1}{2z}\bigg(\frac{1}{\sqrt{1-4z}}+\frac{1+X}{\sqrt{(1+X)^2-4z}}\bigg)\bigg|^{-1/2},
\end{eqnarray*}
and $z_1=z_1(X,t)$ is the unique solution to (\ref{2_lm2}) in the range $z\in(0,(1+X)^2/4)$. For $t\to\infty$, with $t=\log(1/\varepsilon)+T$, (\ref{2_th3_2}) simplifies to (\ref{2_th3_1T}), so that (\ref{2_th3_1T}) applies both to $X<-1$ and $X\in(-1,0)$.
\item[(iii)] $(X,t)$ in Region VII:
\begin{equation}\label{2_th3_3}
p(x,t)\sim\widetilde{G}(z_3)\exp\Big[\frac{1}{\varepsilon}\widetilde{f}(z_3)\Big],
\end{equation}
where now $z_3$ is the unique solution to (\ref{2_lm2}).
\item[(iv)] $(X,t)$ in Region V:
\begin{equation}\label{2_th3_4}
p(x,t)\sim\widetilde{G}(z_1)\exp\Big[\frac{1}{\varepsilon}\widetilde{f}(z_1)\Big]+\widetilde{G}(z_3)\exp\Big[\frac{1}{\varepsilon}\widetilde{f}(z_3)\Big],
\end{equation}
where $z_1$ and $z_3$ are the minimal and maximal of the three solutions to (\ref{2_lm2}). 

The equation $\widetilde{f}(z_1(X,t);X,t)=\widetilde{f}(z_3(X,t);X,t)$ define a curve $t_\Gamma$ in the $(X,t)$-plane within Region V, that satisfies $t_d<t_\Gamma<t_c$. For $t_d<t<t_\Gamma$, the term involving $z_3$ in (\ref{2_th3_4}) dominates, while for $t_\Gamma<t<t_c$ the term involving $z_1$ dominates. The three curves $t_c$, $t_d$ and $t_\Gamma$ all go through the cusp point $(X_{cusp},t_{cusp})$.
\item[(v)] $(X,t)$ in Region IV:
\begin{equation}\label{2_th3_5}
p(x,t)\sim G(X,t)\exp\Big[\frac{1}{\varepsilon}F(X,t)\Big]+\widetilde{G}(z_1)\exp\Big[\frac{1}{\varepsilon}\widetilde{f}(z_1)\Big],
\end{equation}
where $z_1$ is the unique solution of (\ref{2_lm2}). 

The equation $F(X,t)=\widetilde{f}(z_1(X,t);X,t)$ defines a curve $t_\Gamma$ within Region IV such that the term in (\ref{2_th3_5}) involving $F$ is the dominate one for $t<t_\Gamma$, while the term involving $\widetilde{f}(z_1)$ dominates for $t>t_\Gamma$. The curves $t_\Gamma$ in Regions V and IV are smooth continuations of one another, and as $X\to 0^-$ in Region IV, $t_\Gamma$ coincides with $t_1(0;X_0)=t_2(0;X_0)$ in Theorem \ref{theorem2}.
\item[(vi)] $(X,t)\approx (X_{cusp},t_{cusp})$:

Introducing the scaled variables $(\xi,\eta)$ with
\begin{equation}\label{2_th3_6XT}
X=X_{cusp}+\sqrt{\varepsilon}\,\xi=\Big(-1+\frac{2}{\sqrt{X_0+4}}\Big)+\sqrt{\varepsilon}\,\xi,\quad t=t_{cusp}+\sqrt{\frac{(X_0+4)(X_0+3)}{X_0}}\sqrt{\varepsilon}\,\xi+\varepsilon^{3/4}\eta,
\end{equation}
we have
\begin{eqnarray}\label{2_th3_6}
p(x,t)&\sim & \frac{\varepsilon^{-1/4}}{2\pi}\sqrt{\frac{X_0+3}{X_0}}\sqrt{\sqrt{X_0}\sqrt{X_0+3}+X_0+2}\times\exp\bigg[-{\xi^2}\frac{\sqrt{X_0}+2\sqrt{X_0+3}}{4\sqrt{X_0}}\,\bigg]\nonumber\\
&&\times\exp\bigg\{\frac{1}{\varepsilon}\frac{\sqrt{X_0}}{4(X_0+4)}\Big[-2(X_0+3)^{3/2}+\sqrt{X_0}(2X_0+9)\Big]\bigg\}\nonumber\\
&&\times\exp\bigg[-\frac{\xi}{\sqrt{\varepsilon}}\,\frac{\big(\sqrt{X_0}+2\sqrt{X_0+3}\big)^2}{4\sqrt{X_0(X_0+3)(X_0+4)}}-\frac{\eta}{\varepsilon^{1/4}}\frac{3}{4(X_0+3)}\bigg]\times \mathcal{J}(\xi,\eta),
\end{eqnarray}
where $\mathcal{J}$ has the integral representation 
\begin{eqnarray}\label{2_th3_6J}
\mathcal{J}(\xi,\eta)=\int_{-\infty}^{\infty}\exp\bigg\{-\eta w+\Big[\frac{(X_0+4)(X_0+3)}{X_0}\Big]^{3/2}\xi w^2-\frac{2(X_0+3)^{11/2}}{3X_0^{5/2}}w^4\bigg\}\,dw.
\end{eqnarray}

\end{enumerate}
\end{theorem}

The asymptotic structure of $p(x,t)$ is not quite as complicated as Theorem \ref{theorem3} suggests, for the following reasons. In (\ref{2_region1})-(\ref{2_region7}) we defined the various regions as open sets, which exclude the curves that separate pairs of regions. For example, the curve $t=t_*$ for $X<X_{cusp}$, which separates Regions I and II, is not contained in either region. This may suggest that another expansion is needed for $t\approx t_*$ but this is in fact not the case. If we expand (\ref{2_th3_1}) for $t\uparrow t_*$ (and $X<X_{cusp}$) we obtain the expression in (\ref{2_th3_1t*}) for $t_*-t=O(\varepsilon^{1/3})$. But expanding the Region II result in (\ref{2_th3_2}) for $t\downarrow t_*$ (and $X<X_{cusp}$) again leads to the expression in (\ref{2_th3_1t*}), which is a smooth function of the variables $\Lambda$ and $X$. Hence (\ref{2_th3_2}) is smooth through the curve $t=t_*$ and the asymptotics of $p(x,t)$ do not undergo a phase transition here. Similarly, going from Region III to Region VII through the curve $t=t_*$ (with $X_{cusp}<X<X_*$) does not lead to a phase transition, as (\ref{2_th3_3}) is the smooth continuation of (\ref{2_th3_1}). Also, going from Region IV to Region V through $t=t_*$ (with $X_*<X<0$) does not lead to a phase transition as $\widetilde{G}(z_3)e^{\widetilde{f}(z_3)/\varepsilon}$ is the natural continuation of $Ge^{F/\varepsilon}$. Thus there are no transitions near $t=t_*$ for any $-1<X<0$.

Now suppose we fix $X\in(X_{cusp},X_*)$ and increase $t$, thus going from Region III to VII to V and finally VI. As mentioned above going from III to VII is just a smooth continuation. Going from VII to V leads to the appearance of the term involving $z_1$ in (\ref{2_th3_4}). However along $t=t_d$ the term involving $z_3$ is the dominant one asymptotically, as $\widetilde{f}(z_3)>\widetilde{f}(z_1)$, so the transition involves only the exponentially small error term. Going from V to VI leads to the disappearance of the term involving $z_3$, but by then the term involving $z_1$ dominates, and this occurs precisely as $t$ increases past $t_\Gamma$, by definition. Thus along $t\approx t_c$ the phase transition again involves only the exponentially small correction term, which now arises from $z_3$. 

If we fix $X\in(X_*,0)$ and increase $t$ we go from Region III to IV to V and finally VI. Going from III to IV leads to the birth of the term involving $z_1$ in (\ref{2_th3_5}), but for $t\approx t_d$ this is exponentially smaller than the first term in (\ref{2_th3_5}). Going from IV to V involves crossing $t_*$ and this is a smooth continuation that we already discussed. Finally, going from V to VI again leads to the disappearance of the $z_3$ term, which is by then dominated by $z_1$, as was the case for $X\in(X_{cusp},X_*)$. 

This discussion shows that the only true transition in the leading term for $p(x,t)$ occurs along the curve $t_\Gamma$, as when $t$ increases past this curve $\widetilde{G}(z_1)e^{\widetilde{f}(z_1)/\varepsilon}$ begins to dominate $\widetilde{G}(z_3)e^{\widetilde{f}(z_3)/\varepsilon}$ or $Ge^{F/\varepsilon}$ (as the latter two are smooth continuations of one another). Then the only other non-uniformity or transition in the asymptotics of $p(x,t)$ occurs in a small neighborhood of the cusp point $(X_{cusp},t_{cusp})$, defined precisely via the scaling in (\ref{2_th3_6XT}), where the approximation for $p(x,t)$ involves the complicated integral in (\ref{2_th3_6J}), in terms of the local $(\xi,\eta)$ coordinate system.

In section 4 we shall establish Theorem \ref{theorem3} (and also Lemmas \ref{lemma1} and \ref{lemma2}) by using the Laplace transform in (\ref{2_th1_2}) and the integral in (\ref{2_invlap}), which we shall expand using the saddle point method and singularity analysis. Then the curves $t_*$, $t_c$ and $t_d$ will naturally come out of that analysis. However, in section 5 we shall briefly discuss an alternate, more geometrical approach. This will regain all of the results in Theorem \ref{theorem3}, albeit in a slightly different form. The alternate approach will make it clear that only the exterior and interior (defined by $t_d<t<t_c$ for $X_{cusp}<X<0$) of the cusped region need to be considered separately. Thus Regions I, II, III, VI and VII will lead to one result for $p(x,t)$, while regions IV and V will need a different analysis. 

In Figure \ref{figure4} we indicate the curves $t_c$, $t_d$ and $t_\Gamma$, for the three initial conditions $X_0=1,\,2,$ and $5$. Finally we comment that to obtain results for the scale $X=O(\sqrt{\varepsilon})$ (thus $x=O(1)$) we can simply expand the Regions III and IV results for small $X$. When $X=0$ the Region III expression agrees with the approximation $p_{_{BM}}(0,t;X_0)$ in (\ref{2_pbm}), and this remains true in Region IV as long as $t<t_\Gamma\big|_{X=0}$. For $t>t_\Gamma\big|_{X=0}$ setting $X=0^-$ in (\ref{2_th3_5}), where now $p(0,t)\sim \widetilde{G}(z_1)e^{\widetilde{f}(z_1)/\varepsilon}\big|_{X=0}$, we obtain (\ref{2_th2_3}) with $X=0^+$ (and $z_*=z_*(0^+,t)$). Thus all the asymptotic approximations are continuous along $X=0$.

Theorems \ref{theorem2} and \ref{theorem3} apply only for a fixed $X_0>0$, and break down if $X_0\to 0^+$. Note that the cusp point $(X_{cusp},t_{cusp})\to (0,0)$ as $X_0\to 0$ and then a completely separate asymptotic analysis is needed, which we do not attempt here. Theorem \ref{theorem1} holds for any $x_0\ge 0$ but the asymptotics of $p(x,t)$ will be very different if $X_0=x_0/\beta$ is small. When $x_0<0$ a different expression for the Laplace transform must be used (see \cite{lee}) and again the asymptotics of $p(x,t)$ will be completely different.

\section{Analysis for $x>0$}

Here we establish the asymptotic results in Theorem \ref{theorem2}, which apply for $X>0$ (thus $x>0$) and $X=0$. We shall use (\ref{2_th1_1}) in (\ref{2_invlap}) and expand the resulting contour integral by the saddle point method and singularity analysis. General reference on methods for asymptotically evaluating integrals and sums can be found in \cite{won}, \cite{ble}, \cite{ben}. We shall also need to use various properties of the parabolic cylinder functions $D_{-\theta}(\cdot)$, and they are summarized in \cite{abr}-\cite{tem}.

Since we are taking $\beta\to +\infty$, the function in (\ref{2_th1_1}) will have many poles in the range $\theta$ real and $\theta\in(-\beta^2/4,0)$, in addition to the branch cut along $\Im(\theta)=0$, $\Re(\theta)\in(-\infty,-\beta^2/4)$. Let us denote by $\theta=\theta_N$ the $N^{\textrm{th}}$ pole of (\ref{2_th1_1}), with $
\theta_0=0$ corresponding to the steady state limit $p(x,\infty)$, and $\theta_N<0$ with $0<|\theta_1|<|\theta_2|<\cdots$. The number of poles in the range $(-\beta^2/4,0)$ is certainly finite for finite $\beta$, and in view of (\ref{2_th1_1}) and the definition of $R_\beta(\theta)$, the poles are the solutions of
\begin{equation}\label{3_Rbeta}
\sqrt{\theta+\frac{\beta^2}{4}}=R_\beta(\theta)=\frac{D'_{-\theta}(-\beta)}{D_{-\theta}(-\beta)},\quad -\frac{\beta^2}{4}<\theta\le 0.
\end{equation}
Note that $D_{-\theta}(-\beta)$ is well known to be an entire function of $\theta$ (and also $\beta$) \cite{magnus} and that zeros of $D_{-\theta}(-\beta)=0$ are not poles of the $\widehat{p}(x,\theta)$, again by (\ref{2_th1_1}). We have $D_0(-\beta)=e^{-\beta^2/4}$ so that $\theta=0$ is a root of (\ref{3_Rbeta}) for any $\beta>0$. Note also that $D'_{-\theta}(z)$ is the derivative of the function with respect to its argument $z$, hence $D'_{-\theta}(-\beta)=-\frac{\partial}{\partial\beta}D_{-\theta}(-\beta)$.

In \cite{lee} we estimated the poles $\theta_N$ for $\beta\to +\infty$ and $N$ fixed, showing that $\theta_N\approx -N$ with an exponentially small error:
\begin{eqnarray}\label{3_p1p2}
\theta_N &=& -N+\frac{e^{-\beta^2/2}}{\sqrt{2\pi}\,(N-1)!}\beta^{2N-3}\big[1+o(1)\big],\quad \beta\to\infty,\nonumber\\
&=& -N+\frac{\varepsilon^{3/2-N}}{(N-1)!}\exp\Big(-\frac{1}{2\varepsilon}\Big)\big[1+o(1)\big],\quad \varepsilon\to 0.
\end{eqnarray}
This applies also to $N=0$ since $\theta_0=0$ exactly. The present analysis, however, will require computing the contribution(s) from poles $\theta_N$ for large values of $N$, and we shall see that a natural scaling will have $N=O(\varepsilon^{-1})=O(\beta^2)$.

Using (\ref{2_th1_1}) and (\ref{2_invlap}) we can decompose $p(x,t)$ into the two parts
$p(x,t)=p_1(x,t)+p_2(x,t)$,
with
\begin{equation}\label{3_p1}
p_1(x,t) = \frac{1}{2\pi i}\int_{Br} e^{\theta t}e^{\beta(x_0-x)/2}\,\frac{e^{-|x-x_0|\sqrt{\theta+\beta^2/4}}}{\sqrt{\beta^2+4\theta}}\,d\theta=\frac{1}{2\sqrt{\pi t}}\,\exp\bigg[-\frac{1}{4\varepsilon t}(t+X-X_0)^2\bigg],
\end{equation}
where we set $(x,x_0)=\beta(X,X_0)$ and $\beta^2=\varepsilon^{-1}$, and 
\begin{equation*}\label{3_p2}
p_2(x,t)=\frac{1}{2\pi i}\int_{Br} e^{\theta t}e^{\beta(x_0-x)/2}e^{-(x+x_0)\sqrt{\theta+\beta^2/4}}\,{H}(\theta)\,d\theta,
\end{equation*}
\begin{equation}\label{3_p2H}
{H}(\theta)=-\frac{1}{\sqrt{\beta^2+4\theta}}+\frac{1}{\sqrt{\theta+\beta^2/4}-R_\beta(\theta)}.
\end{equation}
We again have $\Re(\theta)>0$ on the $Br$ contours, so these lie to the right of all poles. The contour integral in (\ref{3_p1}) may be explicitly evaluated and the result is precisely the density of the free space Brownian motion in (\ref{2_pbm}). The term $p_2$ gives the derivation of $p$ from $p_{_{BM}}$, and represents the effects of the interface $X=0$. Note that the first term in the expression for $H(\theta)$ would correspond to absorption at $X=0$, but the second part of $H$ indicates a much more complicated reflection/transmission law, that involves the parabolic cylinder functions.

We can obtain an alternate representation of $p_2(x,t)$ by shifting the contour $Br$, on which $\Re(\theta)>0$, into another vertical contour $Br(\theta_{sa})$ on which $-\beta^2/4<\Re(\theta_{sa})<0$. But in shifting we must take into account the residues from the poles. We thus have:
\begin{lemma}\label{lemma3}
The difference $p(x,t)-p_{_{BM}}(x,t)=p_2(x,t)$ between the density of the present diffusion and a free space Brownian motion with drift has the alternate representation
\begin{eqnarray}\label{3_p2lemma3}
p_2(x,t)&=&\sum_{N=0}^{\lfloor -\theta_{sa}\rfloor}e^{\theta_N t}e^{\beta(x_0-x)/2}e^{-(x+x_0)\sqrt{\theta_N+\beta^2/4}}\,h_N+\frac{1}{2\pi i}\int_{Br(\theta_{sa})}e^{\theta t}e^{\beta(x_0-x)/2}e^{-(x+x_0)\sqrt{\theta+\beta^2/4}}{H}(\theta)d\theta\nonumber\\
&\equiv&\mathrm{SUM}+\mathrm{INT},
\end{eqnarray}
where $Br(\theta_{sa})$ is a vertical contour on which $\Re(\theta_{sa})\in(-\beta^2/4,0)$, $\lfloor\cdot\rfloor$ is the greatest integer function, and
\begin{equation}\label{3_hN}
h_N = \lim_{\theta\to \theta_N}\bigg\{\big(\theta-\theta_N\big)\bigg[-\frac{1}{\sqrt{\beta^2+4\theta}}+\frac{1}{\sqrt{\theta+\beta^2/4}-R_\beta(\theta)}\bigg]\bigg\}=\bigg[\frac{1}{\sqrt{\beta^2+4\theta}}-\frac{d}{d\theta}\,R_\beta(\theta)\bigg]^{-1}\bigg|_{\theta=\theta_N}
\end{equation}
is the residue of ${H}(\theta)$ at the pole $\theta_N$, where $0=\theta_0>\theta_1>\theta_2>\cdots$.
\end{lemma}

By choosing $\theta_{sa}$ to be the saddle point of the integrand, the integral in (\ref{3_p2lemma3}) may be easily evaluated for $\beta\to +\infty$, and then we will need to estimate $h_N$ for $\beta$ large and evaluate the residue sum by a Laplace type method. Getting the asymptotics of $p$ then involves determining whether the integral in (\ref{3_p2lemma3}) dominates the sum, and then whether $p_1$ dominates $p_2$. Note that for short times the approximation $p\sim p_{_{BM}}$ will certainly hold, so that $p_1\gg p_2$, but for $t\to\infty$, $p_1=p_{_{BM}}\to 0$ and then $p_2\gg p_1$, as $p_2$ must approach the steady state limit in (\ref{2_steady}) or (\ref{2_steady_appro}). In fact the steady state limit comes from the pole $\theta_0=0$, so that for sufficiently large time, the residue sum in (\ref{3_p2lemma3}) must dominate the integral, and the summand with $N=0$ must ultimately dominate the other terms in the sum. The details of these asymptotic transitions will ultimately lead to Theorem \ref{theorem2}. 

We first consider the integral in (\ref{3_p2lemma3}), and use the fact that $\beta\to +\infty$ to simplify the function ${H}(\theta)$. Using asymptotic properties of the parabolic cylinder functions $D_{-\theta}(-\beta)$ and scaling $\theta$ as
$\theta={\varepsilon}/{\phi}$ and $\theta_{sa}={\varepsilon}/{\phi_{**}}$,
we obtain the following result (see \cite{lee})
\begin{equation}\label{3_D'/D}
\frac{D'_{-\theta}(-\beta)}{D_{-\theta}(-\beta)}=-\frac{1}{\sqrt{\varepsilon}}\sqrt{\phi+\frac{1}{4}}+\frac{\sqrt{\varepsilon}\big[1+\sqrt{1+4\phi}\big]}{2(1+4\phi)}+O(\varepsilon^{3/2}).
\end{equation}
Thus we are taking $\beta\to +\infty$ with $\theta=O(\varepsilon^{-1})=O(\beta^2)$. The result in (\ref{3_D'/D}) applies for $\phi>0$, and also for $\phi\in(-1/4,0)$ as long as we remain bounded away from the zeros of $D_{-\theta}(-\beta)=D_{-\phi/\varepsilon}(-1/\sqrt{\varepsilon})$. For $\varepsilon\to 0^+$ these zeros will be very close to $\phi=-\varepsilon N$ (then $\theta=-N$) where $N$ is an integer. We shall show that the poles of ${H}(\theta)$ are also very close to negative integer values. Thus as long as $-\theta_{sa}$ in (\ref{3_p2lemma3}) is not close to a positive integer, we can use (\ref{3_D'/D}) to approximate $H(\theta)$, thus obtaining 
\begin{equation}\label{3_mathcalH}
H(\theta)=\frac{1+\sqrt{1+4\phi}}{2(1+4\phi)^2}\,\varepsilon^{3/2}+O(\varepsilon^{5/2}).
\end{equation}
Note that the second term in the expansion in (\ref{3_D'/D}) is needed to obtain the leading term in (\ref{3_mathcalH}), since $\sqrt{\theta+\beta^2/4}-R_\beta(\theta)\sim 2\sqrt{\theta+\beta^2/4}=\sqrt{4\theta+\beta^2}$ in this limit. With (\ref{3_mathcalH}) and the scaling $(x,x_0)=\beta(X,X_0)$ the integral in (\ref{3_p2lemma3}) becomes
\begin{equation}\label{3_INT}
\mathrm{INT}=\frac{\sqrt{\varepsilon}}{2\pi i}\int_{Br(\phi_{**})}\frac{1+\sqrt{1+4\phi}}{2(1+4\phi)^2}\,e^{(X_0-X)/(2\varepsilon)}\,\big[1+O(\varepsilon)\big]\exp\bigg\{\frac{1}{\varepsilon}\bigg[\phi\, t-\frac{X+X_0}{2}\sqrt{1+4\phi}\bigg]\bigg\}\,d\phi.
\end{equation}
The right side of (\ref{3_INT}) is in the standard form for applying the saddle point method for $\varepsilon\to 0^+$. There is a saddle at 
\begin{equation}\label{3_phistar}
\frac{d}{d\phi}\bigg[\phi\, t-\frac{X+X_0}{2}\sqrt{1+4\phi}\bigg]=0,\quad \textrm{i.e.,}\quad \phi=\frac{(X+X_0)^2}{4t^2}-\frac{1}{4}\equiv \phi_{**}(X,t).
\end{equation}
By integrating through the saddle in the imaginary direction we obtain the standard estimate 
\begin{equation}\label{3_INTapp}
\mathrm{INT}\sim \frac{\varepsilon}{4\sqrt{\pi}\sqrt{X+X_0}}\frac{\sqrt{1+4\phi_{**}}+1}{\big[1+4\phi_{**}\big]^{5/4}}\,e^{(X_0-X)/(2\varepsilon)}\exp\bigg\{\frac{1}{\varepsilon}\bigg[\phi_{**} t-\frac{X+X_0}{2}\sqrt{1+4\phi_{**}}\bigg]\bigg\}.
\end{equation}
But, in view of (\ref{3_phistar}), the expression in (\ref{3_INTapp}) is the same as the right side of (\ref{2_th2_1}), since $\sqrt{1+4\phi_{**}}=(X+X_0)/t$. Note that the saddle $\phi_{**}\lessgtr 0$ for $X+X_0\lessgtr t$. When $t<X+X_0$, $\theta_{sa}>0$ and the sum in (\ref{3_p2lemma3}) is absent, and then (\ref{3_INTapp}) gives the leading term for $p_2(x,t)$. For $t>X+X_0$ we must estimate the residue sum in (\ref{3_p2lemma3}). Note also that the right side of (\ref{2_th2_1}) (or (\ref{3_INTapp})) is always exponentially smaller than $p_{_{BM}}$. The two become comparable only along $X=0$ (more precisely for $X=O(\varepsilon)$), but then $p_{_{BM}}$ is still larger since (\ref{3_INTapp}) contains an additional $O(\varepsilon)$ algebraic factor that is absent from $p_{_{BM}}$.

To estimate the sum in (\ref{3_p2lemma3}) for $t>X+X_0$, we need to evaluate the residue(s) $h_N$. We first locate the poles $\theta_N$ for $\beta\to +\infty$ with $N=O(\beta^2)$. Note that for $N=O(1)$, (\ref{3_p1p2}) applies. Using the symmetry relation (see \cite{gra}, p. 1030)
\begin{equation}\label{3_Dsymm}
D_{-\theta}(-\beta)=e^{-\pi i\theta}D_{-\theta}(\beta)+\frac{\sqrt{2\pi}}{\Gamma(\theta)}e^{\pi i(1-\theta)/2}D_{\theta-1}(i\beta)
\end{equation} 
and noting that the poles satisfy
\begin{equation}\label{3_Dpoles}
\frac{d}{d\beta}D_{-\theta}(-\beta)+\sqrt{\theta+\frac{\beta^2}{4}}\,D_{-\theta}(-\beta)=0
\end{equation}
we must solve 
\begin{equation}\label{3_Deqn}
D'_{-\theta}(\beta)+\sqrt{\theta+\frac{\beta^2}{4}}\,D_{-\theta}(\beta)=\frac{\sqrt{2\pi}}{\Gamma(\theta)}e^{\pi i\theta/2}(-i)\bigg[\frac{d}{d\beta}D_{\theta-1}(i\beta)+\sqrt{\theta+\frac{\beta^2}{4}}\,D_{\theta-1}(i\beta)\bigg].
\end{equation}
In the limit $\beta\to +\infty$ with $\theta=O(\beta^2)$ we have the estimates (see \cite{abr}), for $\beta^2+4\theta>0$,
\begin{equation}\label{3_Dbapp}
D_{-\theta}(\beta)\sim \Big(1+\frac{\theta}{W^2}\Big)^{-1/2}\,W^{-\theta}\,\exp\Big[\frac{\theta}{2}-\frac{\beta}{4}\sqrt{\beta^2+4\theta}\Big],
\end{equation}
$$W=\frac{1}{2}\big[\beta+\sqrt{\beta^2+4\theta}\big]=\frac{1}{2\sqrt{\varepsilon}}\Big[1+\sqrt{1+4\phi}\Big]$$
and
\begin{equation}\label{3_Dibapp}
D_{\theta-1}(i\beta)\sim -ie^{\pi i\theta/2}\Big(1+\frac{\theta}{W^2}\Big)^{-1/2}\,W^{\theta-1}\,\exp\Big[-\frac{\theta}{2}+\frac{\beta}{4}\sqrt{\beta^2+4\theta}\Big].
\end{equation}
Also, analogous estimates of the derivatives show that
\begin{equation}\label{3_D'/Dapp}
\frac{D'_{-\theta}(\beta)}{D_{-\theta}(\beta)}=-\sqrt{\theta+\frac{\beta^2}{4}}-\frac{1}{2(\beta^2+4\theta)}\Big[\beta-\sqrt{\beta^2+4\theta}\Big]+O(\beta^{-3}),
\end{equation}
\begin{equation*}
\frac{\frac{d}{d\beta}D_{\theta-1}(i\beta)}{D_{\theta-1}(i\beta)}\sim\sqrt{\theta+\frac{\beta^2}{4}}.
\end{equation*}
Now, we know that the solutions of (\ref{3_Dpoles}) and (\ref{3_Deqn}) are on the real axis. But in view of (\ref{3_Dbapp}) and (\ref{3_Dibapp}) the ratio of the left side of (\ref{3_Deqn}) to the bracketed factor in the right side is roughly of the order of $O\big[\exp\big(\theta-\beta\sqrt{\theta+\beta^2/4}\big)\big]$, which is exponentially small in this limit. Thus if (\ref{3_Deqn}) is to hold then $1/\Gamma(\theta)$ must also be exponentially small, so that $\theta$ must be very close to a zero of $\Gamma^{-1}(\theta)$, or a pole of $\Gamma(\theta)$, and these occur at $\theta=0,-1,-2,\cdots$. We know that $\theta=0$ is a pole exactly, and near a pole we can approximate $\Gamma$ by 
\begin{equation}\label{3_1/Gamma}
\frac{1}{\Gamma(\theta)}=(-1)^NN!(\theta+N)+O((\theta+N)^2)
\end{equation}
which is just a Taylor expansion. Using (\ref{3_1/Gamma}) and (\ref{3_D'/Dapp}) we can replace (\ref{3_Deqn}) by the asymptotic relation
\begin{equation}\label{3_Dasym}
-\frac{1}{2}\frac{1}{\beta^2+4\theta}\Big[\beta-\sqrt{\beta^2+4\theta}\Big]\frac{D_{-\theta}(\beta)}{D_{\theta-1}(i\beta)}\sim\sqrt{2\pi}(-1)^{N+1}i\,N!\,(\theta+N)\,e^{\pi i\theta/2}\sqrt{\beta^2+4\theta}
\end{equation}
and this can be used to estimate $\theta_N+N$, which we summarize below.
\begin{lemma}\label{lemma4}
The roots $\theta_N$ of (\ref{3_Rbeta}) for $\beta= 1/\sqrt{\varepsilon}\to +\infty$ satisfy $\theta_N\sim -N$ with the correction term
\begin{equation}\label{3_lemma4}
\theta_N+N\sim\frac{e^{-N}}{\sqrt{2\pi}(N-1)!(\beta^2-4N)^{3/2}}\,\exp\bigg[-\frac{\beta}{2}\sqrt{\beta^2-4N}+2N\log\bigg(\frac{\beta+\sqrt{\beta^2-4N}}{2}\bigg)\bigg].
\end{equation}
\end{lemma}
Note that (\ref{3_lemma4}) remains valid for $N=O(1)$, as then it reduces to (\ref{3_p1p2}), since for $N\ll \beta^2$ we have $\beta\sqrt{\beta^2-4N}=\beta^2-2N+o(1)$. Expression (\ref{3_lemma4}) thus applies for all $N=O(1)$ and $N=O(\beta^2)$, as long as $N/\beta^2<1/4$. Lemma \ref{lemma4} follows immediately from (\ref{3_Dasym}) by using (\ref{3_Dbapp}) and (\ref{3_Dibapp}) to estimate the parabolic cylinder functions, and asymptotically replacing $\theta$ by $-N$ in all terms except $\theta+N$. Also, (\ref{3_lemma4}) is consistent with $\theta_0=0$, if we define $(-1)!=\Gamma(0)=\infty$. For $N\gg 1$, we can approximate $(N-1)!=\Gamma(N)$ by $N^Ne^{-N}\sqrt{2\pi/N}$, due to Stirling's formula.

Having estimated the location of the poles we now evaluate asymptotically the residues $h_N$ in (\ref{3_p2lemma3}). We use (\ref{3_p2H}) and (\ref{3_Dsymm}) and note that only the second term in (\ref{3_p2H}) will contribute to the residues. Thus,
\begin{eqnarray}\label{3_reshN}
&&\frac{1}{\sqrt{\theta+\frac{\beta^2}{4}}-R_\beta(\theta)}=\frac{D_{-\theta}(-\beta)}{\sqrt{\theta+\frac{\beta^2}{4}}\,D_{-\theta}(-\beta)+\frac{d}{d\beta}D_{-\theta}(-\beta)}\nonumber\\
&=&\frac{D_{-\theta}(\beta)+\frac{\sqrt{2\pi}}{\Gamma(\theta)}\,e^{\pi i(1+\theta)/2}\,D_{\theta-1}(i\beta)}{\sqrt{\theta+\frac{\beta^2}{4}}\,D_{-\theta}(\beta)+D'_{-\theta}(\beta)+i\frac{\sqrt{2\pi}}{\Gamma(\theta)}\,e^{\pi i\theta/2}\big[\sqrt{\theta+\frac{\beta^2}{4}}\,D_{\theta-1}(i\beta)+\frac{d}{d\beta}D_{\theta-1}(i\beta)\big]}
\end{eqnarray}
and the residue will be the numerator in (\ref{3_reshN}) evaluated at $\theta=\theta_N$, divided by the derivative of the denominator at $\theta=\theta_N$. Since $\Gamma^{-1}(\theta)\sim (-1)^NN!\,(\theta+N)$ we use (\ref{3_lemma4}), (\ref{3_Dbapp}) and (\ref{3_Dibapp}) to conclude that the numerator is asymptotically the same as $D_N(\beta)$, while the derivative of the denominator is asymptotically dominated by $\frac{d}{d\theta}\Gamma^{-1}(\theta)\sim (-1)^NN!,\; \theta\to\theta_N$. Hence the residue is asymptotically given by 
\begin{eqnarray}\label{3_reshNfinal}
h_N &\sim &\frac{(-1)^N\,e^{\pi i N/2}}{i\sqrt{2\pi}\,N!\,\sqrt{\beta^2-4N}}\frac{D_N(\beta)}{D_{-N-1}(i\beta)}\sim 2(\theta_N+N)\frac{\beta^2-4N}{\beta-\sqrt{\beta^2-4N}} \nonumber\\
&\sim&\frac{1}{N!}\frac{\beta+\sqrt{\beta^2-4N}}{2\sqrt{2\pi}\sqrt{\beta^2-4N}}\Big(\frac{\beta+\sqrt{\beta^2-4N}}{2}\Big)^{2N}\,e^{-N}\,\exp\Big[-\frac{\beta}{2}\sqrt{\beta^2-4N}\Big].
\end{eqnarray}
Here we again used (\ref{3_Dbapp}), (\ref{3_Dibapp}), (\ref{3_Dasym}) and (\ref{3_lemma4}).
Using (\ref{3_reshNfinal}) we have thus approximated the summand in (\ref{3_p2lemma3}) as follows.
\begin{lemma}\label{lemma5}
The summand in (\ref{3_p2lemma3}) has the asymptotic expansion
\begin{equation}\label{3_lemma5}
\sqrt{\varepsilon}\,g(\varepsilon N)\exp\Big[\frac{1}{\varepsilon}f(\varepsilon N)\Big],\quad 0<\varepsilon N<\frac{1}{4}
\end{equation}
where
$$f(z)=f(z;X,t)=-z\,t-\frac{X+X_0}{2}\sqrt{1-4z}+\frac{X_0-X}{2}-\frac{1}{2}\sqrt{1-4z}+2z\log\Big[\frac{1+\sqrt{1-4z}}{2\sqrt{z}}\Big],$$
$$g(z)=\frac{1+\sqrt{1-4z}}{4\pi\sqrt{1-4z}\sqrt{z}},$$
and for $N=O(1)$ the expansion is 
\begin{equation}\label{3_lemma5N}
\frac{1}{N!}\frac{1}{\sqrt{2\pi}}\varepsilon^{-N}\,\exp\Big[-\frac{1}{\varepsilon}\Big(X+\frac{1}{2}\Big)\Big]\,e^{-N(t-X-X_0)}.
\end{equation}
\end{lemma}

To obtain (\ref{3_lemma5}) we simply multiplied (\ref{3_reshNfinal}), after expanding $N!$ by Stirling's formula, by the factors that multiply $h_N$ in (\ref{3_p2lemma3}), using also $\theta_N\sim -N$ and $x\sqrt{\theta_N+\beta^2/4}\sim\frac{1}{2}\beta X\sqrt{\beta^2-4N}=\frac{1}{2}\varepsilon^{-1}X\sqrt{1-4\varepsilon N}$, etc. To obtain (\ref{3_lemma5N}) we also used $\sqrt{\beta^2-4N}\sim\beta=\varepsilon^{-1/2}$ for $N=O(1)$.

To evaluate asymptotically the sum in (\ref{3_p2lemma3}) we first note that for $t=O(1)$ the expression in (\ref{3_lemma5N}) becomes asymptotically larger with increasing $N$, due to the factor $\varepsilon^{-N}$. Then surely the main asymptotic behavior cannot come from the range $N=O(1)$, and we thus consider the sum
\begin{equation}\label{3_sum}
\sum_{N=0}^{\lfloor-\phi_{**}/\varepsilon\rfloor}\sqrt{\varepsilon}g(\varepsilon N)e^{f(\varepsilon N)/\varepsilon},
\end{equation}
for $t>X+X_0$, since then $-\phi_{**}>0$. This is a Laplace type sum and the main contribution will come from the global maxima of the function $f(z)$ over the range $z\in(0,-\phi_{**})$. Note that $f(z)$ is smooth for $z\in(0,1/4)$, and $-\phi_{**}<1/4$ for any fixed $t$ $(>X+X_0)$, in view of (\ref{3_phistar}). We have
\begin{equation}\label{3_f'z}
f'(z)=-t+\frac{X+X_0}{\sqrt{1-4z}}+2\log\bigg(\frac{1+\sqrt{1-4z}}{2\sqrt{z}}\bigg)
\end{equation}
and 
\begin{equation}\label{3_f''z}
f''(z)=\frac{2(X+X_0)}{(1-4z)^{3/2}}-\frac{1}{z\sqrt{1-4z}}.
\end{equation}
Thus $f'(z)\to +\infty$ as either $z\to 0^+$ or $(1/4)^-$, and also
\begin{equation*}\label{3_f'phi}
f'(-\phi_{**})=2\log\bigg(\frac{t+X+X_0}{\sqrt{t^2-(X+X_0)^2}}\bigg)>0.
\end{equation*}
Thus $f(z)$ always has a local maximum at $z=-\phi_{**}$, which corresponds to the upper limit on the sum in (\ref{3_sum}). For $t$ sufficiently small we have $f'(z)>0$ for $z\in(0,1/4)$ in view of (\ref{3_f'z}). However for $t$ large $f'(z)=0$ will have solution(s). Solving $f''(z)=0$ leads to 
$z=\widetilde{z}\equiv 1/[2(X+X_0)+4]$
and thus the minimum value of $f'(z)$ occurs at $\widetilde{z}$ and the minimum value is
\begin{equation}\label{3_minf'}
\min_{0<z<1/4}f'(z)=-t+\sqrt{(X+X_0)(X+X_0+2)}+2\log\bigg(\frac{\sqrt{X+X_0+2}+\sqrt{X+X_0}}{\sqrt{2}}\bigg).
\end{equation}
Now, $\widetilde{z}$ may or may not lie within the interval $(0,-\phi_{**})$, but it will lie in this range for $t$ sufficiently large, as $-\phi_{**}\to 1/4$ as $t\to \infty$. By equating $\widetilde{z}$ to $-\phi_{**}$ we see that $\widetilde{z}\in(0,-\phi_{**})$ precisely when
\begin{equation}\label{3_t>}
t>\sqrt{(X+X_0)(X+X_0+2)}.
\end{equation}
If (\ref{3_t>}) holds then certainly $t>X+X_0$, which is needed for the sum in (\ref{3_p2lemma3}) to come into play. It follows that for $t>\sqrt{(X+X_0)(X+X_0+2)}$, $f'$ will have a unique minimum in the range $z\in(0,-\phi_{**})$, and if $f'(\widetilde{z})<0$ then $f'=0$ will have exactly two roots in this range. In view of (\ref{3_minf'}) let us define
\begin{equation}\label{3_t+}
t_+=t_+(X;X_0)=\sqrt{(X+X_0)(X+X_0+2)}+2\log\bigg(\frac{\sqrt{X+X_0+2}+\sqrt{X+X_0}}{\sqrt{2}}\bigg).
\end{equation}
Then for $t>t_+$ the function $f'$ has two zeros, call these $0<z_*<z_{**}<|\phi_{**}|$, and $z_*$ is a local maximum of $f$ (since $f''(z_*)<0$) while $z_{**}$ is a local minimum.

We next compute the contributions to (\ref{3_sum}) from $z=-\phi_{**}$ and then from $z=z_*$. By the Laplace method, this requires expanding the summand in (\ref{3_sum}) about a maximum of $f$. Setting
\begin{equation}\label{3_Ntilde}
\widetilde{N}=\widetilde{N}(X,t)=\bigg\lfloor\frac{t^2-(X+X_0)^2}{4\varepsilon t^2}\bigg\rfloor=\bigg\lfloor-\frac{\phi_{**}}{\varepsilon}\bigg\rfloor
\end{equation}
the contribution from $z=-\phi_{**}$ is 
\begin{equation}\label{3_redz**}
\sqrt{\varepsilon}\,g(\varepsilon \widetilde{N})\,e^{f(\varepsilon \widetilde{N})/\varepsilon}\sum_{m=0}^\infty e^{-mf'(\varepsilon \widetilde{N})}\,\big[1+O(\varepsilon)\big]\sim\frac{\sqrt{\varepsilon}\,g(\varepsilon \widetilde{N})}{1-e^{-f'(\varepsilon\widetilde{N})}}\,e^{f(\varepsilon \widetilde{N})/\varepsilon},
\end{equation}
where we set $N=\widetilde{N}-m$ in the sum over $N$ in (\ref{3_p2lemma3}). The contribution from $z=z_*$ (when $t>t_+$) is given by, using the Laplace method, 
\begin{eqnarray*}\label{3_redz*}
&&\sqrt{\varepsilon}\,g(z_*)\sum_{N}e^{f(z_*)/\varepsilon}\,\exp\bigg[\frac{1}{2\varepsilon}f''(z_*)(z-z_*)^2\bigg]\nonumber\\
&\sim &\frac{1}{\sqrt{\varepsilon}}\,g(z_*)\,e^{f(z_*)/\varepsilon}\int_{-\infty}^\infty\exp\bigg(-\frac{|f''(z_*)|}{2\varepsilon}u^2\bigg)\,du=\sqrt{2\pi}\frac{g(z_*)}{\sqrt{-f''(z_*)}}\,e^{f(z_*)/\varepsilon}.
\end{eqnarray*}
We summarize the asymptotics of the sum in (\ref{3_p2lemma3}) below.
\begin{lemma}\label{lemma6}
For $t>X+X_0$ the sum in (\ref{3_p2lemma3}) has the asymptotic expansions (with $t_+$ as in (\ref{3_t+}) and $\widetilde{N}$ in (\ref{3_Ntilde})):
\begin{enumerate}
\item[(i)] $X+X_0<t<t_+(X;X_0)$
\begin{equation}\label{3_le6i}
\mathrm{SUM}\sim\sqrt{\varepsilon}\frac{t}{4\pi}\frac{(t+X+X_0)^{3/2}}{\sqrt{t-X-X_0}}\frac{1}{(X+X_0)^2}\,e^{f(\varepsilon \widetilde{N})/\varepsilon}\equiv \mathrm{ENDPOINT},
\end{equation}
$$f(\varepsilon\widetilde{N})=-\varepsilon \,t\widetilde{N}-\frac{X+X_0+1}{2}\sqrt{1-4\varepsilon\widetilde{N}}+\frac{X_0-X}{2}+2\varepsilon\widetilde{N}\log\bigg(\frac{1+\sqrt{1-4\varepsilon\widetilde{N}}}{2\sqrt{\varepsilon\widetilde{N}}}\bigg).$$
\item[(ii)] $t>t_+(X;X_0)$
\begin{equation}\label{3_le6ii}
\mathrm{SUM}\sim \mathrm{ENDPOINT}+\frac{\sqrt{2\pi}\,g(z_*)}{\sqrt{-f''(z_*)}}\,e^{f(z_*)/\varepsilon},
\end{equation}
where $z_*$ is the unique solution to
$$t=\frac{X+X_0}{\sqrt{1-4z_*}}+2\log\bigg(\frac{1+\sqrt{1-4z_*}}{2\sqrt{z_*}}\bigg).$$
\item[(iii)] $t\approx X+X_0$ with $t=X+X_0+\varepsilon\,\widetilde{t}$, $\widetilde{t}=O(1)$, and $M=M(X,\widetilde{t})=\lfloor \widetilde{t}/(2(X+X_0))\rfloor$,
\begin{equation}\label{3_le6iii}
\mathrm{SUM}\sim\frac{\varepsilon^{-M}}{\sqrt{2\pi}\,M!}\exp\bigg[-\frac{1}{\varepsilon}\Big(X+\frac{1}{2}\Big)\bigg].
\end{equation}
\end{enumerate}
\end{lemma}

When $t-(X+X_0)=\varepsilon\,\widetilde{t}=O(\varepsilon)$, $\lfloor -\phi_{**}/\varepsilon\rfloor\sim M=O(1)$ and only a few terms in (\ref{3_p2lemma3}) (up to $N=M$) are present, and the last term dominates the others, leading to (\ref{3_le6iii}). This case is not important to the asymptotics of $p(x,t)$, or even $p_2(x,t)$. The right side of (\ref{3_le6i}) is the same as (\ref{3_redz**}), after we evaluate $f'$ and $g$ at $\varepsilon\, \widetilde{N}\sim-\phi_{**}$. In (\ref{3_le6ii}) the last term is the contribution from the local maximum of $f$ at $z=z_*$, and we note that $\sqrt{2\pi}\,g(z_*)\,[-f''(z_*)]^{-1/2}$ is precisely the algebraic factor that multiplies $e^{f(z_*)/\varepsilon}$ in (\ref{2_th2_2}), in view of (\ref{3_f''z}) and the definition of $g(z)$ below (\ref{3_lemma5}). Also, $f(z_*)$ in (\ref{2_th2_2f}) is equivalent to that computed from the expression below (\ref{3_lemma5}), when we use the fact that $z_*$ satisfies (\ref{2_th2_2z*}). The leading term in (\ref{3_le6ii}) can come from either the endpoint term, if $f(-\phi_{**})>f(z_*)$, or from $z_*$, if $f(-\phi_{**})<f(z_*)$. Along $t=t_+$ the endpoint contribution will be dominant, but for $t$ sufficiently large $f(-\phi_{**})\to -\infty$ and $f(z_*)\to -X-1/2$ so the $z_*$ contribution will dominate. However, we are only interested in estimating $p_2$, and then the full density $p$.

To obtain the leading term for $p_2$ $(=p-p_{_{BM}})$ we must add (\ref{3_INTapp}) to the results in Lemma \ref{lemma6}. For $t<X+X_0$, only the integral is present. For $X+X_0<t<t_+$, $p_2$ will be asymptotic to the sum of (\ref{3_INTapp}) and (\ref{3_le6i}), but the former dominates since $-(t-X-X_0)^2/(4t)>f(-\phi_{**})$ is equivalent to 
\begin{equation}\label{3_equiv}
X+\frac{X+X_0}{2t}>\frac{1}{2}\bigg[1-\frac{(X+X_0)^2}{t^2}\bigg]\log\bigg[\frac{t+X+X_0}{\sqrt{t^2-(X+X_0)^2}}\bigg],
\end{equation}
and (\ref{3_equiv}) is obviously true in view of the inequality $\xi^{2m+1}>\xi^{2m+1}/(2m+1)$ $(m=1,2,3,\cdots)$, which when summed over $m$ yields
$$\frac{2\xi}{1-\xi^2}>\sum_{m=0}^\infty\frac{2\xi^{2m+1}}{2m+1}=\log\bigg(\frac{1+\xi}{1-\xi}\bigg),\quad 0<\xi<1.$$
Setting $\xi=(X+X_0)/t$ we see that the left side of (\ref{3_equiv}) dominates the right side, even without the $X$ term in the former. Hence $p_2$ can be asymptotic to either (\ref{3_INTapp}), which must be true for $t<X+X_0$, or to the $z_*$ contribution in (\ref{3_le6ii}), which must occur for very large times $t$. The two contributions are comparable when $-X-(t-X-X_0)^2/(4t)=f(z_*)$, which can only occur for $t>t_+$. In view of (\ref{2_th2_2f}) this equation is equivalent to
\begin{equation}\label{3_equiv2}
\frac{t}{4}+\frac{(X+X_0)^2}{4t}=\frac{X+X_0+1}{2}\sqrt{1-4z_*}+\frac{X+X_0}{\sqrt{1-4z_*}}\,z_*.
\end{equation}
Setting $U=\sqrt{1-4z_*}$ so that $z_*=(1-U^2)/4$, (\ref{3_equiv2}) becomes the quadratic equation
\begin{equation*}\label{3_quadra}
(X+X_0+2)U+\frac{X+X_0}{U}=t+\frac{(X+X_0)^2}{t}
\end{equation*}
whose solution is 
\begin{equation}\label{3_U}
U=\frac{1}{2(X+X_0+2)}\bigg[t+\frac{(X+X_0)^2}{t}+\sqrt{\mathcal{B}}\bigg],
\end{equation}
where $\mathcal{B}$ is defined below (\ref{2_th2_2t1}). Then using (\ref{3_U}) in (\ref{2_th2_2z*}) yields the equation in (\ref{2_th2_2t1}) for the space time curve $t_1(X;X_0)$, which is the lower curve in Figure \ref{figure1}. We have thus established items (i) and (ii) in Theorem \ref{theorem2}. To establish item (iii) we define a second transition curve $t_2(X;X_0)$, along which $-(t+X-X_0)^2/(4t)=f(z_*)$, and this leads to (\ref{2_th2_2t2}). Then for $t>t_2$, $p(x,t)$ becomes asymptotic to the term in (\ref{3_le6ii}) arising from $z_*$. In section 5 we give more geometric interpretation to the various terms in the expansion of $p(x,t)$, and of the curve $t_+$. This completes the analysis for fixed $t=O(1)$. 

Now suppose we consider large time scales $t\approx\log(1/\varepsilon)$, setting $T=t-\log(1/\varepsilon)=O(1)$, as in (\ref{2_Tscale}). Now $\lfloor-\phi_{**}/\varepsilon\rfloor\sim 1/4$ and we have $\varepsilon^{-N}e^{-Nt}=e^{-NT}=O(1)$.
Then the asymptotics of the sum in (\ref{3_p2lemma3}) may be obtained by summing the $N=O(1)$ result in (\ref{3_lemma5N}), thus obtaining
\begin{equation}\label{3_Tscale}
p(x,t)\sim\frac{1}{\sqrt{2\pi}}\exp\bigg[-\frac{1}{\varepsilon}\bigg(X+\frac{1}{2}\bigg)\bigg]\sum_{N=0}^\infty\frac{e^{-N(T-X-X_0)}}{N!}.
\end{equation}
Evaluating the series in (\ref{3_Tscale}) gives the double exponential (in $T$) approximation in (\ref{2_pT}), and for $T\to\infty$ ($t-\log(1/\varepsilon)\to\infty$) this reduces to the steady state formula for $X$ (or $x$) $>0$ in (\ref{2_steady_appro}). On the $T$-scale the sum in (\ref{3_p2lemma3}) remains a discrete sum and every pole contributes roughly equally to its asymptotic behavior. Though the argument leading to the evaluation of (\ref{3_p2lemma3}) is different for the $t$ and $T$ time scales, the result in (\ref{2_th2_3}) actually reduces to (\ref{2_pT}), for large times $t$. This is easily seen by expanding $z_*$ and $f(z_*)$ for $t\to\infty$, which leads to (\ref{2_z*sim}) and (\ref{2_fz*sim}). This completes the derivation of Theorem \ref{theorem2}.

\section{Analysis for $x<0$}

We shall obtain Theorem \ref{theorem3} by using (\ref{2_th1_2}) in (\ref{2_invlap}) and again expanding the integral for $\beta\to +\infty$ with $x_0=\beta X_0=O(\beta)$ (and $X_0>0$). The analysis for $x<0$ is somewhat more complicated than for $x>0$, and (\ref{2_th1_2}) involves the special functions $D_{-\theta}(\cdot)$, which we have to approximate in various ranges of $\theta$, $\beta$, and $x$.

First we let, as in section 3, $\theta_{sa}$ have any real value in the range $(-\beta^2/4,\infty)$, and obtain from (\ref{2_invlap}) the following alternate representation.
\begin{lemma}\label{lemma7}
For $x<0$, $p(x,t)$ is given by
\begin{eqnarray}\label{4_lemma7}
p(x,t)&=&\sum_{N=0}^{\lfloor -\theta_{sa}\rfloor}e^{\theta_N t}\,e^{-x^2/4}\,e^{\beta(x_0-x)/2}\,e^{-x_0\sqrt{\theta_N+\beta^2/4}}\,\frac{D_{-\theta_N}(-\beta-x)}{D_{-\theta_N}(-\beta)}\,h_N\nonumber\\
&&+\frac{1}{2\pi i}\int_{Br(\theta_{sa})}e^{\theta\,t}\,e^{-x^2/4}\,e^{\beta(x_0-x)/2}\,\frac{e^{-x_0\sqrt{\theta+\beta^2/4}}}{\sqrt{\theta+\beta^2/4}-R_\beta(\theta)}\frac{D_{-\theta}(-\beta-x)}{D_{-\theta}(-\beta)}\,d\theta,
\end{eqnarray}
where $h_N$ is as in (\ref{3_hN}), $\theta_{sa}$ is real, and $Br(\theta_{sa})$ is a vertical Bromwich contour on which $\Re(\theta)=\theta_{sa}\in(-\beta^2/4,\infty)$. If $\theta_{sa}>0$ the sum in (\ref{4_lemma7}) is absent, and we assume that $\theta_{sa}\ne \theta_N$ for any $N$.
\end{lemma}

Next we give some asymptotic results for the parabolic cylinder functions. These may be found, for example, in \cite{tem}. We shall scale $\theta=\phi/\varepsilon=O(\varepsilon^{-1})$ and $x=\beta X=O(\beta)=O(1/\sqrt{\varepsilon})$, and we note that the function $D_{-\phi/\varepsilon}(Z/\sqrt{\varepsilon})$ has the integral representation
\begin{equation}\label{4_Dint}
D_{-\phi/\varepsilon}\bigg(\frac{Z}{\sqrt{\varepsilon}}\bigg)=\frac{\varepsilon^{\phi/(2\varepsilon)}}{i\sqrt{2\pi\varepsilon}}\,\exp\bigg(\frac{Z^2}{4\varepsilon}\bigg)\int_{Br}\exp\bigg[\frac{1}{\varepsilon}\widehat{F}(U;\phi,Z)\bigg]\,dU
\end{equation}
where 
\begin{equation}\label{4_Fhat}
\widehat{F}=-\phi\log U-ZU+\frac{1}{2}U^2,
\end{equation}
$\Re(U)>0$ on the vertical contour $Br$, and $\log U$ is defined to be real for $U$ real and positive. From (\ref{4_Fhat}) we see that (\ref{4_Dint}) has saddle points where $\widehat{F}'(U)=0$, or
\begin{equation}\label{4_upm}
U_\pm=\frac{1}{2}\Big[Z\pm\sqrt{Z^2+4\phi}\Big],
\end{equation}
and also $\widehat{F}_{UU}=1+\phi/U^2$. Then using the steepest descent method we find that the major contribution comes from the saddle $U_+$ in (\ref{4_upm}) and obtain, for $\phi>0$, 
\begin{eqnarray}\label{4_Dappx}
D_{-\phi/\varepsilon}\bigg(\frac{Z}{\sqrt{\varepsilon}}\bigg)&\sim&\varepsilon^{\phi/(2\varepsilon)}\,\frac{1}{(Z^2+4\phi)^{1/4}}\bigg(\frac{2\phi}{\sqrt{Z^2+4\phi}-Z}\bigg)^{1/2}\nonumber\\
&&\times\exp\bigg\{\frac{1}{\varepsilon}\bigg[\frac{\phi}{2}-\frac{Z}{4}\sqrt{Z^2+4\phi}-\phi\log\bigg(\frac{Z+\sqrt{Z^2+4\phi}}{2}\bigg)\bigg]\bigg\}.
\end{eqnarray}
We shall also need the expansion of $D_{-\phi/\varepsilon}(\cdot)$ for $\phi<0$, and in this range the function may have zeros. The expression in (\ref{4_Dappx}) still holds if $\phi<0$ and $Z>0$, as long as $Z^2+4\phi>0$. However (\ref{4_Dappx}) does not apply for $\phi<0$ and $Z<0$, which is what we need to approximate $D_{-\theta}(-\beta)$ in (\ref{4_lemma7}) if $\theta_{sa}<0$. But, we can use the symmetry relation in (\ref{3_Dsymm}) for the case $\phi<0$ and $Z<0$ (setting $-\beta=Z$). Then (\ref{4_Dappx}) may be used to estimate both terms in the right side of (\ref{3_Dsymm}) (to get $D_{\phi/\varepsilon-1}(-iZ/\sqrt{\varepsilon})$ we just replace $(\phi,Z)$ by $(\varepsilon-\phi, -iZ)$ and simplify the expression). As long as $\Gamma^{-1}(\theta)\ne 0$, or $\phi$ remains bounded away from $-\varepsilon N;\;N=0,1,2,\cdots$, the second term in (\ref{3_Dsymm}) is exponentially larger than the first and we thus obtain, for $\phi=\varepsilon\theta<0$, $Z<0$ and $Z^2+4\phi>0$, 
\begin{eqnarray}\label{4_Dappx2}
D_{-\phi/\varepsilon}\bigg(\frac{Z}{\sqrt{\varepsilon}}\bigg)&\sim & 2\sin(\pi\theta)\,\varepsilon^{\phi/(2\varepsilon)}\frac{1}{(Z^2+4\phi)^{1/4}}\bigg(\frac{-2\phi}{\sqrt{Z^2+4\phi}-Z}\bigg)^{1/2}\nonumber\\
&&\times \exp\bigg\{\frac{1}{\varepsilon}\bigg[\frac{\phi}{2}-\frac{Z}{4}\sqrt{Z^2+4\phi}-\phi\log\bigg(\frac{-Z-\sqrt{Z^2+4\phi}}{2}\bigg)\bigg]\bigg\}.
\end{eqnarray}
This applies as long as $\theta\ne 0,-1,-2,-3,\cdots$, more precisely as long as we are away from exponentially small neighborhoods of $\theta=-N$, which is where the zeros of $D_{-\phi/\varepsilon}(Z/\sqrt{\varepsilon})$ are located for $\phi,\,Z<0$. We can estimate these zeros precisely using (\ref{3_Dsymm}), similarly as we estimated the roots of (\ref{3_Rbeta}), but the present set of zeros will not be required for the analysis. Setting $Z=-1-X$ and then $Z=-1$ and taking the ratio of the two expressions, we thus have the following.
\begin{lemma}\label{lemma8}
For $\varepsilon\to 0^+$ and $\phi>0$ we have the asymptotic result
\begin{eqnarray}\label{4_le8}
\frac{D_{-\phi/\varepsilon}(-(1+X)/\sqrt{\varepsilon})}{D_{-\phi/\varepsilon}(-1/\sqrt{\varepsilon})}&\sim& \bigg[\frac{1+4\phi}{(1+X)^2+4\phi}\bigg]^{1/4}\bigg[\frac{\sqrt{1+4\phi}+1}{\sqrt{(1+X)^2+4\phi}+1+X}\bigg]^{1/2}\nonumber\\
&&\times\exp\bigg[\frac{\phi}{\varepsilon}\log\bigg(\frac{\sqrt{1+4\phi}-1}{\sqrt{(1+X)^2+4\phi}-1-X}\bigg)\bigg]\nonumber\\
&&\times\exp\bigg\{\frac{1}{4\varepsilon}\Big[-\sqrt{1+4\phi}+(1+X)\sqrt{(1+X)^2+4\phi}\Big]\bigg\}.
\end{eqnarray}
For $-1<X<0$ (thus $1+X>0$) and $-(1+X)^2/4<\phi<0$, (\ref{4_le8}) holds as long as $-\phi/\varepsilon=-\theta$ remains bounded away from the non-negative integers $0,1,2,\cdots$.
\end{lemma}

Note that in (\ref{4_Dappx}) the most rapidly varying factor in $\varepsilon$ is $\varepsilon^{\phi/(2\varepsilon)}=O[\exp(\varepsilon^{-1}\log\varepsilon)]$, but this disappears upon taking the ratio. In going from $\phi>0$ to $\phi<0$ in (\ref{4_le8}) we must only replace
$$\log\bigg[\frac{\sqrt{1+4\phi}-1}{\sqrt{(1+X)^2+4\phi}-1-X}\bigg]\quad\textrm{by}\quad \log\bigg[\frac{1-\sqrt{1+4\phi}}{1+X-\sqrt{(1+X)^2+4\phi}}\bigg].$$
Also, the factor $\sin(\pi\theta)$ in (\ref{4_Dappx2}) disappears upon taking the ratio of (\ref{4_Dappx2}) with the two different $Z$ values. But, (\ref{4_le8}) does become invalid when $\theta\approx -\varepsilon N$ for integer values of $N$. Lemma \ref{lemma8} does not apply for $\phi<0$ and $X<-1$, but this range will not be needed in our analysis.

We will need to approximate the parabolic cylinder functions when $\phi<0$, and also when $\phi\approx-(1+X)^2/4$. Then we need the following result.
\begin{lemma}\label{lemma9}
For $\phi<0$ and $\varepsilon\to 0$, with $\phi\approx-(1+X)^2/4$ and 
\begin{equation}\label{4_le9con}
\phi=-\frac{1}{4}(1+X)^2+\Big(\frac{1+X}{2}\Big)^{2/3}\,\varepsilon^{2/3}\,\delta,\quad\delta=O(1)
\end{equation}
we have
\begin{eqnarray}\label{4_Dappx3}
D_{-\phi/\varepsilon}\bigg(-\frac{1+X}{\sqrt{\varepsilon}}\bigg)&\sim&\varepsilon^{\phi/(2\varepsilon)}\frac{\sqrt{2\pi}\,(1+X)^{1/3}}{2^{1/3}\,\varepsilon^{1/6}}\,\exp\bigg\{\frac{1}{\varepsilon}\bigg[-\frac{(1+X)^2}{8}-\phi\log\Big(\frac{1+X}{2}\Big)\bigg]\bigg\}\nonumber\\
&&\times\Big[\sin(\pi\theta)\,\mathrm{Bi}(\delta)+\cos(\pi\theta)\,\mathrm{Ai}(\delta)\Big],
\end{eqnarray}
where $\mathrm{Ai}(\cdot)$ and $\mathrm{Bi}(\cdot)$ are the Airy functions. The expression in (\ref{4_Dappx3}) holds away from zeros of $\big[\sin(\pi\theta)\,\mathrm{Bi}(\delta)+\cos(\pi\theta)\,\mathrm{Ai}(\delta)\big]$.
\end{lemma}

This result is obtained by using well-known \cite{olv} results for approximating parabolic cylinder functions by Airy functions. This corresponds to situations where the saddle points $U_\pm$ in (\ref{4_upm}) nearly coalesce. Lemmas \ref{lemma8} and \ref{lemma9} give the approximations to the integrand in (\ref{4_lemma7}) that we shall need. Now consider the sum in (\ref{4_lemma7}), which is needed for $\theta_{sa}<0$. We have already estimated $h_N$ in (\ref{3_reshNfinal}) and shown that $\theta_N\sim -N$ with an exponentially small error. We thus have
\begin{equation}\label{4_Dratio}
\frac{D_{-\theta_N}(-\beta-X)}{D_{-\theta_N}(-\beta)}\sim\frac{D_N(-(1+X)/\sqrt{\varepsilon})}{D_N(-1/\sqrt{\varepsilon})}=\exp\bigg[-\frac{1}{\varepsilon}\bigg(\frac{X^2}{4}+\frac{X}{2}\bigg)\bigg]\,\frac{He_N((1+X)/\sqrt{\varepsilon})}{He_N(1/\sqrt{\varepsilon})},
\end{equation}
where $He_N(\cdot)$ is the $N^{\textrm{th}}$ Hermite polynomial (thus $He_0(Z)=1$, $He_1(Z)=Z$, etc.), which may be obtained, for example, from the integral
\begin{equation}\label{4_Hermite}
He_N(Z)=\frac{1}{\sqrt{2\pi}}\int_{-\infty}^\infty (Z+iv)^N\,e^{-v^2/2}\,dv.
\end{equation}
We shall need (\ref{4_Dratio}) only for $-1<X<0$ so the argument(s) of $He_N(\cdot)$ will be always positive. Using the saddle point method we can expand (\ref{4_Hermite}) for $N\to\infty$, $Z\to\infty$ and a natural scaling has $N=O(Z^2)$. Then different expansions are needed for $4N/Z^2<1$, and for $4N-Z^2=O(N^{1/3})$. Below we summarize only the main results.
\begin{lemma}\label{lemma10}
For $\varepsilon\to 0$, $N\to \infty$, and $X\in(-1,0]$ we have 
\begin{eqnarray}\label{4_le10}
He_N\bigg(\frac{1+X}{\sqrt{\varepsilon}}\bigg)&\sim&\sqrt{\frac{N}{2}}\,\bigg[\frac{(1+X)^2}{4\varepsilon}-N\bigg]^{-1/4}\bigg[\frac{1+X}{2}-\sqrt{\frac{(1+X)^2}{4\varepsilon}-N}\,\bigg]^{-1/2}\nonumber\\
&&\times\exp\bigg\{\frac{(1+X)^2}{4\varepsilon}-\frac{N}{2}-\frac{1+X}{\sqrt{\varepsilon}}\sqrt{\frac{(1+X)^2}{4\varepsilon}-N}\bigg\}\nonumber\\
&&\times\exp\bigg\{N\log\bigg[\frac{1+X}{2\sqrt{\varepsilon}}+\sqrt{\frac{(1+X)^2}{4\varepsilon}-N}\bigg]\bigg\}
\end{eqnarray}
where $(1+X)^2/(4\varepsilon N)>1$. For $(1+X)^2/(4\varepsilon N)\approx 1$, with the scaling
\begin{equation}\label{4_le10scale}
\varepsilon N=\frac{1}{4}(1+X)^2-\bigg(\frac{1+X}{2}\bigg)^{2/3}\,\varepsilon^{2/3}\,W,\quad W=O(1),
\end{equation}
we have
\begin{equation}\label{4_le10_2}
He_N\bigg(\frac{1+X}{\sqrt{\varepsilon}}\bigg)\sim 2^{1/6-N}\,\sqrt{\pi}\,\bigg(\frac{1+X}{\sqrt{\varepsilon}}\bigg)^{N+1/3}\,\exp\bigg[\frac{(1+X)^2}{8\varepsilon}\bigg]\,\mathrm{Ai}(W),
\end{equation}
where $\mathrm{Ai}(\cdot)$ is the Airy function.
\end{lemma}

We now asymptotically evaluate the expression in Lemma \ref{lemma7}, starting with the integral. We scale $\theta=\phi/\varepsilon$, use (\ref{4_le8}) to approximate the ratio of parabolic cylinder functions, and note that as long as we remain bounded away from the poles $\theta_N$, (\ref{3_D'/D}) may be used to approximate $R_\beta(\theta)$, and then
\begin{equation*}\label{4_Rbeta}
\frac{d\theta}{\sqrt{\theta+\beta^2/4}-R_\beta(\theta)}\sim\frac{d\phi}{\sqrt{\varepsilon}\sqrt{1+4\phi}}.
\end{equation*}
The integrand in (\ref{4_lemma7}) thus becomes 
\begin{eqnarray}\label{4_integrand}
&&\frac{1}{\sqrt{\varepsilon}}\,\frac{1}{(1+4\phi)^{1/4}}\,\frac{1}{\big[(1+X)^2+4\phi\big]^{1/4}}\,\bigg[\frac{\sqrt{1+4\phi}+1}{\sqrt{(1+X)^2+4\phi}+1+X}\bigg]^{1/2}\nonumber\\
&&\times\exp\bigg\{\frac{1}{\varepsilon}\bigg[-\frac{X^2}{4}-\frac{X}{2}+\frac{X_0}{2}+\phi\, t-\Big(\frac{X_0}{2}+\frac{1}{4}\Big)\sqrt{1+4\phi}+\frac{1+X}{4}\sqrt{(1+X)^2+4\phi}\bigg]\bigg\}\nonumber\\
&&\times\exp\bigg[\frac{\phi}{\varepsilon}\log\bigg(\frac{\sqrt{1+4\phi}-1}{\sqrt{(1+X)^2+4\phi}-1-X}\bigg)\bigg]\nonumber\\
&&\equiv\varepsilon^{-1/2}\mathcal{G}(\phi;t,X,X_0)\;\exp\big[\mathcal{F}(\phi;t,X,X_0)/\varepsilon\big],
\end{eqnarray}
and this we now integrate over the $\phi$ variable. So that there are saddle points where $\partial\mathcal{F}/\partial\phi=0$, or
\begin{equation}\label{4_saddlephi}
t-\frac{X_0}{\sqrt{1+4\phi}}+\log\bigg[\frac{\sqrt{1+4\phi}-1}{\sqrt{(1+X)^2+4\phi}-1-X}\bigg]=0.
\end{equation}

We write the solution of (\ref{4_saddlephi}) as $\phi_s=\phi_s(X,t;X_0)$. First observe that $\phi=0$ is a solution when $t-X_0+\log(1+X)=0$ (when $-1<X<0$), and this corresponds precisely to the fluid approximation $X=-1+e^{X_0-t}$ in (\ref{2_Yfl}), where $Y_{fl}(t)<0$ for $t>X_0$. The equation (\ref{4_saddlephi}) applies only for $\phi>-1/4$ and $\phi>-(1+X)^2/4$, as only then (\ref{4_le8}) can be used to approximate the integrand in (\ref{4_lemma7}). In Figure \ref{figure5} we sketch the ``surface" $\phi_s(X,t)$ over the $(X,t)$ plane for $X_0=1$. The curve $X+1=e^{X_0-t}$ divides the $(X,t)$ plane into regions where $\phi_s>0$ and $\phi_s<0$. We have $\phi_s>0$ when $X\le -1$, and also when $0<t<X_0-\log(1+X)$ with $X\in (-1,0]$. If $X=0$ we may solve (\ref{4_saddlephi}) explicitly to give $\phi_s(0,t;X_0)=(X_0^2-t^2)/(4t^2)$, and then
$$\frac{1}{\varepsilon}\bigg[\frac{X_0}{2}+t\,\phi_s(0,t;X_0)-\frac{X_0}{2}\sqrt{1+4\phi_s(0,t;X_0)}\bigg]=-\frac{1}{4\varepsilon\, t}(t-X_0)^2,$$
which is the exponential part in $p_{_{BM}}$ in (\ref{2_pbm}), and this we previously remarked holds also at $X=0$. For any $X\le 0$ we can obtain the behavior of $\phi_s$ at $t\to 0$, by expanding (\ref{4_saddlephi}) for $\phi\to\infty$, which leads to 
$\phi_s\sim{(X_0-X)^2}/(4t^2)$ as $t\to 0$.
For $t\to\infty$ and $X\le -1$ we have $\phi_s\to 0$ and from (\ref{4_saddlephi}) we can obtain the more precise estimate
\begin{equation}\label{4_phisest}
\phi_s\sim\left\{ \begin{array}{ll}
-(X+1)\,e^{X_0-t}, & X<-1\\
e^{2(X_0-t)}, & X=-1
\end{array} \right\},\quad t\to \infty.
\end{equation}
Thus for large times and $X\le -1$ the saddle point becomes exponentially close to the origin in the $\phi$-plane. For $X\le -1$ we have $\partial\phi_s/\partial t<0$ for all times, and $\phi_s$ decreases from $+\infty$ to $0$ as $t$ increases from $0$ to $+\infty$. In contrast, for $-1<X<0$ we again have $\partial\phi_s/\partial t<0$ but now $\phi_s$ decreases from $+\infty$ to $0$ as $t$ increases from $0$ to $X_0-\log(1+X)$ (which is $>0$), and $\phi_s$ further decreases from $0$ to $-(1+X)^2/4$ as $t$ increases from $X_0-\log(1+X)$ to $t_*$, where $t_*$ is obtained by setting $\phi=-(1+X)^2/4$ in (\ref{4_saddlephi}) and this leads to the expression in (\ref{2_t*}). For $\phi<-(1+X)^2/4$ the analysis leading to (\ref{4_saddlephi}) ceases to be valid, as then Lemma \ref{lemma8} no longer applies.

Using the standard saddle point estimate, we observe that $\partial^2\mathcal{F}/\partial \phi^2>0$ at $\phi=\phi_s$ and the steepest descent directions are $arg(\phi-\phi_s)=\pm \pi/2$, so we approximate the integral by 
\begin{equation}\label{4_intapprx}
\frac{1}{\sqrt{2\pi}}\,\mathcal{G}(\phi_s;t,X;X_0)\bigg[\frac{\partial^2\mathcal{F}}{\partial
\phi^2}(\phi_s;t,X;X_0)\bigg]^{-1/2}\,\exp\bigg[\frac{1}{\varepsilon}\mathcal{F}(\phi_s;t,X;X_0)\bigg].
\end{equation}
But, in view of (\ref{4_integrand}) and our definitions of $\mathcal{G}$ and $\mathcal{F}$ above (\ref{4_saddlephi}), (\ref{4_intapprx}) is the same as (\ref{2_th3_1}), with (\ref{2_th3_1F}) and (\ref{2_th3_1G}).

When $\phi_s>0$ the sum in (\ref{4_lemma7}) is absent, so we have established item (i) in Theorem \ref{theorem3}, at least for the range $\{X\le -1\}$ and for $t<X_0-\log(1+X)$ with $X\in (-1,0)$, which correspond to portions of Regions I and III. For the portions of these regions where $\phi_s<0$ we must also estimate the sum in (\ref{4_lemma7}), which we shall do shortly.

For $X<-1$ and large times $t$ we can use the estimate in (\ref{4_phisest}) to simplify (\ref{2_th3_1}), and this leads to (\ref{2_th3_1T}), on the $T$ time scale. Along the fluid approximation $\phi_s=0$, and by expanding (\ref{4_saddlephi}) for small $\phi$ we obtain the estimate
\begin{equation}\label{4_phisappx}
\phi_s\sim\frac{t-X_0+\log(1+X)}{1-2X_0-(1+X)^{-2}},\quad t\to X_0-\log(1+X).
\end{equation}
Using (\ref{4_phisappx}) in (\ref{2_th3_1})-(\ref{2_th3_1G}) gives the Gaussian approximation in (\ref{2_th3_1gau}), which applies only on the $\Delta$-scale, and thus for $X-(-1+e^{X_0-t})=O(\sqrt{\varepsilon})$.

We now take $-1<X<0$ and consider ranges where $-(1+X)^2/4<\phi_s<0$, and also the case(s) $t>t_*$, where $\phi_s$ fails to exist. By the saddle point method, setting $\theta_{sa}=\phi_s/\varepsilon$ allows us to obtain the estimate in (\ref{4_intapprx}) for the integral. The sum in (\ref{4_lemma7}) we then estimate by
\begin{equation}\label{4_sumappx}
\sum_{N=0}^{\lfloor-\phi_s/\varepsilon\rfloor}\sqrt{\varepsilon}\,\widetilde{g}(\varepsilon N)\,\exp\bigg[\frac{1}{\varepsilon}\,\widetilde{f}(\varepsilon N)\bigg]
\end{equation}
where
\begin{eqnarray}\label{4_sumappxf}
\widetilde{f}(z)&=&-\frac{X^2}{4}-\frac{X}{2}+\frac{X_0}{2}-zt-\bigg(\frac{X_0}{2}+\frac{1}{4}\bigg)\sqrt{1-4z}-\frac{1+X}{4}\sqrt{(1+X)^2-4z}\nonumber\\
&&+z\,\log\bigg[\frac{\big(1+X+\sqrt{(1+X)^2-4z}\,\big)\big(1+\sqrt{1-4z}\,\big)}{4z}\bigg]
\end{eqnarray}
and
\begin{equation}\label{4_sumappxg}
\widetilde{g}(z)=\frac{1}{2\pi\,(1-4z)^{1/4}\,\big[(1+X)^2-4z\big]^{1/4}}\bigg[\frac{1+\sqrt{1-4z}}{1+X-\sqrt{(1+X)^2-4z}}\bigg]^{1/2}.
\end{equation}
To obtain (\ref{4_sumappx})-(\ref{4_sumappxg}) we used (\ref{3_reshNfinal}) to estimate $h_N$ (also approximating $N!$ by Stirling's formula), (\ref{4_Dratio}) and (\ref{4_le10}) to estimate the parabolic cylinder functions in (\ref{4_lemma7}), and the facts that $\beta^2=\varepsilon^{-1}$ and $\theta_N\sim -N$. In (\ref{4_sumappxf}) and (\ref{4_sumappxg}) we emphasize the dependence of $z$, but of course $\widetilde{f}$ depends also on $t$ and $X$, and $\widetilde{g}$ depends also on $X$.

Note also that if we identify $z$ with $-\phi$, apart from a factor of $\varepsilon$, the summand in (\ref{4_sumappx}) is nearly identical to (\ref{4_integrand}), which corresponds to the integrand in (\ref{4_lemma7}), except that the term $\sqrt{(1+X)^2-4z}$ in the summand is replaced by $-\sqrt{(1+X)^2+4\phi}$ in (\ref{4_integrand}). Thus (\ref{4_integrand}) and the summand in (\ref{4_sumappx}) may be viewed as different branches of the same function, and in section 5 we explore further this connection. We now analyze Region IV and those portions of Regions I and III where $\phi_s<0$. We evaluate (\ref{4_sumappx}) by the Laplace method, which requires that we find the maximum of $\widetilde{f}(z)$ over the range $z\in[0,-\phi_s]$.

From (\ref{4_sumappxf}) we have
\begin{equation}\label{4_tildef'}
\widetilde{f}'(z)=-t+\frac{X_0}{\sqrt{1-4z}}+\log\bigg[\frac{1+X+\sqrt{(1+X)^2-4z}}{1-\sqrt{1-4z}}\bigg]
\end{equation}
and
\begin{equation}\label{4_tildef''}
\widetilde{f}''(z)=\frac{2X_0}{(1-4z)^{3/2}}-\frac{1}{2z\sqrt{1-4z}}-\frac{1+X}{2z\sqrt{(1+X)^2-4z}}.
\end{equation}
Also, $\widetilde{f}(0)=-(1+X)^2/2$, $\widetilde{f}'(0)=+\infty$ and 
\begin{equation}\label{4_tildef'phis}
\widetilde{f}'(-\phi_s)=\log\bigg[\frac{1+X+\sqrt{(1+X)^2+4\phi_s}}{1+X-\sqrt{(1+X)^2+4\phi_s}}\bigg]>0,
\end{equation}
where we used the fact that $\phi_s$ satisfies (\ref{4_saddlephi}). After some calculation we can explicitly solve $\widetilde{f}''(z)=0$ to find that it has two roots, at $z=z_\pm$, where
\begin{equation*}\label{4_zpm}
z_\pm=\frac{(X_0+1)^2(X+1)^2+2(X_0+1)-3(X+1)^2\pm X_0^{3/2}(X+1)\sqrt{(X_0+4)(X+1)^2-4}}{8(X_0+X+2)(X_0-X)}.
\end{equation*}
For $X\in(-1,0)$ the roots $z_\pm$ are real only when $X\in(X_{cusp},0)$, so we see the cusp point in (\ref{2_Xcusp}) entering the analysis. At the
cusp point we have
\begin{equation}\label{4_zatcusp}
z_+=z_-=\frac{3}{4(X_0+3)}>0,\quad X=X_{cusp}=-1+\frac{2}{\sqrt{X_0+4}}.
\end{equation}
At the points $z=z_\pm$ we need to consider the sign(s) of $\widetilde{f}'(z)$, and we note that, using (\ref{4_saddlephi}),
$\widetilde{f}'(z_-)=0$ implies $t=t_d(X;X_0)$, and $\widetilde{f}'(z_+)=0$ implies $t=t_c(X;X_0)$,
where $t_c$ and $t_d$ are the curves in (\ref{2_tc})-(\ref{2_alpha_d}).

In Region I, $X<X_{cusp}$ and $\widetilde{f}''(z)<0$, and thus $\widetilde{f}'$ decreases from $\infty$ at $z=0$, to the value in (\ref{4_tildef'phis}) when $z=-\phi_s$. Hence $\widetilde{f}'>0$ and the maximum of $\widetilde{f}$ occurs at the upper limit of the sum in (\ref{4_sumappx}). By expanding in Taylor series and noting that $\varepsilon\,\lfloor-\phi_s/\varepsilon\rfloor=-\phi_s+O(\varepsilon)$, we estimate (\ref{4_sumappx}) by
\begin{equation}\label{4_sumapprx2}
\sqrt{\varepsilon}\,\widetilde{g}(-\phi_s)\,\exp\bigg[\frac{1}{\varepsilon}\,\widetilde{f}\bigg(\varepsilon\,\bigg\lfloor-\frac{\phi_s}{\varepsilon}\bigg\rfloor\bigg)\bigg]\sum_{m=0}^{\infty}e^{-m\widetilde{f}'(-\phi_s)}.
\end{equation}
Evaluating the geometric sum in (\ref{4_sumapprx2}), and using (\ref{4_saddlephi}), (\ref{4_sumappxf}) and (\ref{4_tildef'}), we ultimately obtain
\begin{equation}\label{4_sumapprxult}
\frac{\sqrt{\varepsilon}}{8\pi\sqrt{-\phi_s}}\frac{\big(1+\sqrt{1+4\phi_s}\big)^{1/2}}{\big(1+4\phi_s\big)^{1/4}}\frac{\big(1+X+\sqrt{(1+X)^2+4\phi_s}\,\big)^{3/2}}{\big[(1+X)^2+4\phi_s\big]^{3/4}}\,\exp\bigg[\frac{1}{\varepsilon}\widetilde{f}\bigg(\varepsilon\bigg\lfloor-\frac{\phi_s}{\varepsilon}\bigg\rfloor\bigg)\bigg].
\end{equation}
But comparing the magnitude of (\ref{4_sumapprxult}), which is roughly of order $O(e^{\widetilde{f}(-\phi_s)/\varepsilon})$ to the saddle point contribution in (\ref{2_th3_1}) (or (\ref{4_integrand})), which is $O(e^{F/\varepsilon})$, we see that the latter dominates, as $F(X,t)>\widetilde{f}(-\phi_s)$. Hence we have now established (\ref{2_th3_1}) in all of Region I.

In the part of Region III above the fluid approximation we have $\phi_s<0$. Now $t<t_*$ if $X\in(X_{cusp},X_*)$ and $t<t_d$ if $X\in(X_*,0)$, where we recall that $(X_*,t_{**})$ is the intersection of the curves $t_*(X;X_0)$ and $t_d(X;X_0)$, so the intersection point depends only on $X_0$. If we view $\widetilde{f}(z)$ over $z\in(0,(1+X)^2/4)$ then $\widetilde{f}''(z)=0$ at $z=z_\pm$. We may have $-\phi_s<z_-$, $-\phi_s\in(z_-,z_+)$, or $-\phi_s>z_+$, depending on the values of $(X,t)$. However, since $t<t_d$ we have $\widetilde{f}'(z_+)>0$ and $\widetilde{f}'(z_-)>0$ and thus again $\widetilde{f}'(z)>0$ for $z\in[0,-\phi_s]$ and the estimate in (\ref{4_sumapprxult}) holds for the sum in (\ref{4_lemma7}), but the saddle contribution in (\ref{4_integrand}) dominates the sum. We have thus established (\ref{2_th3_1}) for the part of Region III where $\phi_s<0$.

Now consider Region IV, where $t\in(t_d,t_*)$ and $X\in(X_*,0)$. Then $\widetilde{f}'(z_-)<0$ and $\widetilde{f}'(z_+)>0$, so $\widetilde{f}'(z)$ has two zeros, say at $z_1$ and $z_2$, which satisfy $0<z_1<z_-<z_2<z_+$. Also, $z_1$ will be a local maximum of $\widetilde{f}$ while $z_2$ will be a local minimum. The global maximum of $\widetilde{f}(z)$ will come from either $z=z_1$, or the endpoint $z=-\phi_s$. But the endpoint contribution will be dominated again by the saddle point contribution to the integral in (\ref{4_lemma7}). The contribution to the sum from $z_1$ can be obtained by the Laplace estimate, as
\begin{equation}\label{4_reIV}
\sqrt{2\pi}\,\widetilde{g}(z_1)\,\big[-\widetilde{f}''(z_1)\big]^{-1/2}\,\exp\bigg[\frac{1}{\varepsilon}\widetilde{f}(z_1)\bigg].
\end{equation}
But $\sqrt{2\pi}\,\widetilde{g}(z_1)\,\big[-\widetilde{f}''(z_1)\big]^{-1/2}$ is the same as $\widetilde{G}(z_1)$ in (\ref{2_th3_5}). Then in Region IV the asymptotics are governed by either $\varepsilon\,N=z_1$ from the sum or $\phi=\phi_s$ in the integral, and either of these may dominate. This establishes item (v) in Theorem \ref{theorem3}.

In Regions II and V-VII we have $t>t_*$. Now we set $\varepsilon\theta_{sa}=-(1+X)^2/4$ in (\ref{4_lemma7}) and the integrand has no saddle points, as (\ref{4_saddlephi}) has no real solutions. Later we shall examine the case $t\approx t_*$ to see what happens to the solution of (\ref{4_saddlephi}) as $t$ increases past $t_*$. Now we must estimate the sum
\begin{equation}\label{4_sum_regII}
\sum_{N=0}^{\lfloor(1+X)^2/(4\varepsilon)\rfloor}\sqrt{\varepsilon}\,\widetilde{g}(\varepsilon N)\,\exp\bigg[\frac{1}{\varepsilon}\widetilde{f}(\varepsilon N)\bigg].
\end{equation}
Our approximation to the summand fails near the upper limit, where $z=(1+X)^2/4$. But this will be important asymptotically only for $t\approx t_*$, and this case we shall analyze separately. We note that
\begin{equation}\label{4_regIIf'}
\widetilde{f}'\bigg(\frac{(1+X)^2}{4}\bigg)=t_*-t<0
\end{equation}
so now $\widetilde{f}$ will no longer have a local maximum at the upper limit in (\ref{4_sum_regII}). 

In Region V we have $t_d<t<t_c$ if $X\in (X_{cusp},X_*)$ and $t_*<t<t_c$ if $X\in (X_*,0)$. In either case $\widetilde{f}'(z_-)<0$ and $\widetilde{f}'(z_+)>0$, so that $\widetilde{f}'(z)=0$ has three roots, at $z_1<z_2<z_3$, with $z_1$ and $z_3$ corresponding to local maxima of $\widetilde{f}$ and $z_2$ a local minimum. Then the contribution from the maximum at $z_1$ is again given by (\ref{4_reIV}). Adding to this the contribution from $z_3$ leads to (\ref{2_th3_4}) in Theorem \ref{theorem3}. Either contribution may dominate, and the two are comparable along the curve $t_\Gamma$, as previously discussed.

In Region VII we have $X\in(X_{cusp},X_*)$ and $t_*<t<t_d$. Now $\widetilde{f}'(z_-)>0$, $\widetilde{f}'(z_+)>0$ which together with (\ref{4_regIIf'}) shows that $ \widetilde{f}'(z)=0$ has a unique root, call it $z_3$ ($>z_+$), and this leads to (\ref{2_th3_3}) in Theorem \ref{theorem3}.

For Region VI we have $\widetilde{f}'(z_-)<0$ and $\widetilde{f}'(z_+)<0$, and, as always, $\widetilde{f}'(0)=+\infty$. Thus there is again a unique root, which we call $z_1$, and we obtain (\ref{2_th3_2}). 

For Region II, $\widetilde{f}''(z)=0$ has no solutions in $z\in(0,(1+X)^2/4)$, since $X<X_{cusp}$. Then $\widetilde{f}'(0)=+\infty$ and $\widetilde{f}'<0$ at $z=(1+X)^2/4$ so that $\widetilde{f}'$ has a unique zero and we again obtain (\ref{2_th3_2}).

It remains to consider the case $t\approx t_*$, where the approximations leading to the summand in (\ref{4_sumappx}) and (\ref{4_sum_regII}) cease to be valid. Also, we must separately analyze when $(X,t)=(X_{cusp},t_{cusp})$, where the local maxima $z_1$ and $z_3$ of $\widetilde{f}(z)$ coalesce, and the standard Laplace estimate in (\ref{4_reIV}) no longer holds (at the cusp point $z_1=z_3$ and $\widetilde{f}''(z_1)=0$).

We introduce the scaling
\begin{equation}\label{4_Lambda}
t=t_*(X;X_0)+\varepsilon^{1/3}\,\Lambda,\quad\Lambda=O(1)
\end{equation}
as in (\ref{2_th3_1lambda}), so that $t-t_*=O(\varepsilon^{1/3})$. Furthermore we shift the $Br$ contour in (\ref{4_lemma7}) to $\Re(\phi)=-(1+X)^2/4$. Then scaling $\phi$ as in (\ref{4_le9con}), we use (\ref{4_Dappx3}) to approximate $D_{-\phi/\varepsilon}(-(1+X)/\sqrt{\varepsilon})$ and (\ref{4_Dappx2}) to approximate $D_{-\phi/\varepsilon}(-1/\sqrt{\varepsilon})$. Then the integrand in (\ref{4_lemma7}) becomes, in this asymptotic limit,
\begin{equation}\label{4_t*scale}
\frac{\sqrt{\pi}}{2}\,\frac{\big[\sqrt{1-(1+X)^2}+1\big]^{1/2}}{\big[1-(1+X)^2\big]^{1/4}}\,\Big[\mathrm{Bi}(\delta)+\cot(\pi\theta)\mathrm{Ai}(\delta)\Big]\,e^{\mathcal{M}/\varepsilon}
\end{equation}
where
\begin{equation}\label{4_calM}
\mathcal{M}=-\frac{(X+1)^2}{8}-\frac{X^2}{4}-\frac{X}{2}+\frac{X_0}{2}+\phi\,t-\frac{\phi}{2}-\bigg(\frac{X_0}{2}+\frac{1}{4}\bigg)\sqrt{1+4\phi}+\phi\log\bigg(\frac{1-\sqrt{1+4\phi}}{1+X}\bigg),
\end{equation}
and (\ref{4_t*scale}) is to be integrated with respect to $\delta$, since $d\phi=(\frac{1+X}{2})^{2/3}\,\varepsilon^{2/3}\,d\delta$. Using the scaling in (\ref{4_le9con}) we can rewrite $\mathcal{M}$ in (\ref{4_calM}) as
\begin{equation}\label{4_calMFc}
\mathcal{M}=F_c-\frac{1}{4}(1+X)^2\,\varepsilon^{1/3}\,\Lambda+\varepsilon\,\Lambda_1\,\delta+o(\varepsilon),
\end{equation}
as this simply amounts to a Taylor expansion about $\phi=-(1+X)^2/4$. In (\ref{4_calMFc}), $F_c$ is obtained by setting $4\phi=-(1+X)^2$ and $t=t_*$ in (\ref{4_calM}), and then $F_c$ is as in (\ref{2_th3_1Fc}). Also, in (\ref{4_calMFc}) 
$\Lambda_1=[(1+X)/2]^{2/3}\Lambda$.
Using (\ref{4_calMFc}) in (\ref{4_t*scale}) we must evaluate the integral
\begin{equation}\label{4_t*int}
\frac{1}{2\pi i}\int_{Br}\Big[\mathrm{Bi}(\delta)+\cot(\pi\theta)\,\mathrm{Ai}(\delta)\Big]\,e^{\Lambda_1\delta}\,d\delta,
\end{equation}
where $Br$ has been shifted so that $\Re(\delta)=0$ on this contour. Now
$$\cot(\pi\theta)=i\,\frac{e^{i\pi\theta}+e^{-i\pi\theta}}{e^{i\pi\theta}-e^{-i\pi\theta}}\quad\textrm{and}\quad\theta=\frac{\phi}{\varepsilon}=-\frac{1}{4\varepsilon}(1+X)^2+\varepsilon^{-1/3}\bigg(\frac{1+X}{2}\bigg)^{2/3}\,\delta.$$
For $\delta$ positive imaginary we have $\cot(\pi\theta)\sim -i$, while for $\delta$ negative imaginary we have $\cot(\pi\theta)\sim i$. Hence parametrizing the $Br$ contour in (\ref{4_t*int}) and separately calculating the contributions from $\Im(\delta)>0$ and $\Im(\delta)<0$ we are led to
\begin{equation}\label{4_t*intsum}
\frac{1}{2\pi i}\int_{0}^{i\infty}\Big[\mathrm{Bi}(\delta)-i\,\mathrm{Ai}(\delta)\Big]\,e^{\Lambda_1\delta}\,d\delta+\frac{1}{2\pi i}\int_{-i\infty}^0\Big[\mathrm{Bi}(\delta)+i\,\mathrm{Ai}(\delta)\Big]\,e^{\Lambda_1\delta}\,d\delta.
\end{equation}
Using the symmetry relations between Airy functions
$\mathrm{Ai}(z)\pm i\,\mathrm{Bi}(z)=2\,e^{\pm \pi i/3}\,\mathrm{Ai}(z\,e^{\mp 2\pi i/3}),$
performing some contour rotations in (\ref{4_t*intsum}), and then multiplying the result by the factors in (\ref{4_t*scale}) that are not contained in (\ref{4_t*int}), we ultimately obtain
\begin{equation*}\label{4_p+}
\mathcal{P}_+\sim\frac{\big[\sqrt{1-(1+X)^2}+1\big]^{1/2}}{2\sqrt{\pi}\,\big[1-(1+X)^2\big]^{1/4}}\exp\bigg[\frac{F_c}{\varepsilon}-\frac{(1+X)^2}{4\,\varepsilon^{2/3}}\,\Lambda\bigg]\int_0^\infty\mathrm{Ai}(\delta)\,\Big(e^{\Lambda_1\,\omega^2\,\delta}+e^{\Lambda_1\,\omega\,\delta}\Big)\,d\delta,
\end{equation*}
where $\omega=e^{2\pi i/3}$ is a cube root of unity. Thus, $\mathcal{P}_+$ represents the contribution to $p(x,t)$ from the integral in (\ref{4_lemma7}) (with $\theta_{sa}=-(1+X)^2/(4\varepsilon)$), for $t-t_*=O(\varepsilon^{1/3})$.

To compute the contribution from the residue sum in (\ref{4_lemma7}), which we denote by $\mathcal{P}_-$, we use (\ref{4_Dratio}), approximate $He_N((1+X)/\sqrt{\varepsilon})$ using (\ref{4_le10_2}), and approximate $He_N(1/\sqrt{\varepsilon})$ using (\ref{4_le10}) with $X=0$. Note that near $t=t_*$ we must carefully consider the contribution from the upper endpoint in the sum in (\ref{4_lemma7}), since now $\phi_s\sim -(1+X)^2/4$. We thus obtain
\begin{equation}\label{4_upperend}
\frac{\varepsilon^{1/3}}{2^{1/3}\sqrt{\pi}}\frac{\big[1+\sqrt{1-(1+X)^2}\big]^{1/2}}{(1+X)^{2/3}\,\big[1-(1+X)^2\big]^{1/4}}\,\mathrm{Ai}(W)\,e^{\mathcal{S}/\varepsilon},
\end{equation}
where
\begin{equation}\label{4_Scal}
\mathcal{S}=-\frac{(1+X)^2}{8}+\frac{X_0}{2}+\frac{z}{2}-zt-\bigg(\frac{X_0}{2}+\frac{1}{4}\bigg)\sqrt{1-4z}+z\log\bigg(\frac{1+X}{2}\bigg)+z\log\bigg(\frac{1+\sqrt{1-4z}}{2z}\bigg).
\end{equation}
Using the scaling in (\ref{4_Lambda}) and (\ref{4_le10scale}) and expanding (\ref{4_Scal}) about $t=t_*$ and $z=(1+X)^2/4$ we obtain
\begin{equation}\label{4_ScalW}
\mathcal{S}=F_c-\frac{(1+X)^2}{4}\,\varepsilon^{1/3}\,\Lambda+\varepsilon\,\Lambda_1\,W+o(\varepsilon).
\end{equation}
With (\ref{4_upperend}) and (\ref{4_ScalW}), summing (\ref{4_upperend}) over $N$ corresponds to $\varepsilon^{-1/3}$ times an integral over $W>0$, in view of the Euler-Maclaurin formula and (\ref{4_le10scale}), as $W>0$ corresponds to $\varepsilon\,N<(1+X)^2/4$. Then the contribution to (\ref{4_lemma7}) from the sum is, to leading order,
\begin{equation*}\label{4_p-}
\mathcal{P}_-\sim\frac{1}{2\sqrt{\pi}}\,\frac{\big[\sqrt{1-(1+X)^2}+1\big]^{1/2}}{\big[1-(1+X)^2\big]^{1/4}}\exp\bigg[\frac{F_c}{\varepsilon}-\frac{(1+X)^2}{4\,\varepsilon^{2/3}}\,\Lambda\bigg]\,\int_0^\infty\mathrm{Ai}(W)\,e^{\Lambda_1\,W}\,dW.
\end{equation*}

To obtain the expansion of $p$ for $t\approx t_*$ we simply add $\mathcal{P}_+$ to $\mathcal{P}_-$. While both $\mathcal{P}_+$ and $\mathcal{P}_-$ involve integrals of Airy functions, their sum is much simpler. To see this let
$\mathcal{H}(Z)=\int_0^\infty\mathrm{Ai}(W)\big(e^{ZW}+e^{\omega ZW}+e^{\omega^2 ZW}\big)\,dW.$
We have
$\mathcal{H}(0)=3\int_0^\infty\mathrm{Ai}(W)dW=1.$
Using the Airy equation $\mathrm{Ai}''(W)=W\,\mathrm{Ai}(W)$, some integration by parts, and the facts that $1+\omega+\omega^2=0$ and $\omega^3=1$, we obtain
{\setlength\arraycolsep{2pt}
\begin{eqnarray*}\label{4_Hcal'}
\mathcal{H}'(Z)&=&\int_0^\infty W\,\mathrm{Ai}(W)\big(e^{ZW}+\omega\,e^{\omega ZW}+\omega^2\,e^{\omega^2 ZW}\big)\,dW=\int_0^\infty \mathrm{Ai}''(W)\big(e^{ZW}+\omega\,e^{\omega ZW}+\omega^2\,e^{\omega^2 ZW}\big)\,dW\nonumber\\
&=&-\int_0^\infty Z^2\,\mathrm{Ai}(W)\big(e^{ZW}+e^{\omega ZW}+e^{\omega^2 ZW}\big)\,dW=-Z^2\,\mathcal{H}(Z).
\end{eqnarray*}}
\noindent It follows that $\mathcal{H}(Z)=e^{-Z^2/3}$, and then $\mathcal{P}_+$ and $\mathcal{P}_-$ sum to give the explicit formula in (\ref{2_th3_1t*}), since $\Lambda_1^3=(1+X)^2\,\Lambda^3/4$.

In crossing $t=t_*$ we go from Region I to II ($-1<X<X_{cusp}$), from Region III to VII ($X_{cusp}<X<X_*$), or from Region IV to V ($X_*<X<0$). In the first two crossings (\ref{2_th3_1t*}) gives the approximation to $p(x,t)$, but in going from Region IV to V we must also consider the contribution from $z_1$. Then (\ref{2_th3_1t*}) gives the smooth continuation of the term $Ge^{F/\varepsilon}$ in (\ref{2_th3_5}) and it becomes the $z_3$ contribution in (\ref{2_th3_4}).

Also, we can show directly that the approximation to $p(x,t)$ is smooth as $t$ increases past $t_*$. For example, consider the transition from Region I to II. Then for $t<t_*$ (\ref{4_saddlephi}) (with $\phi=\phi_s$) holds. For $\phi_s$ slightly greater than $-(1+X)^2/4$ we let $\phi_s=-(1+X)^2/4+\varepsilon_0$. Then expanding (\ref{4_saddlephi}) about $t=t_*$ and for small $\varepsilon_0$ yields
$$t-t_*+O(\varepsilon_0)=\log\bigg[\frac{1+X-2\sqrt{\varepsilon_0}}{1+X}\bigg]\sim-\frac{2\sqrt{\varepsilon_0}}{1+X}$$
and hence
\begin{equation}\label{4_epsilon0}
\varepsilon_0\sim\frac{1}{4}(1+X)^2\,\big[t-t_*(X;X_0)\big]^2.
\end{equation}

By setting $4\phi_s+(1+X)^2=U$ and expanding (\ref{4_saddlephi}) for small $U$ we obtain
\begin{equation}\label{4_t*expandU}
t-t_*=-\frac{\sqrt{U}}{1+X}+\frac{1-(1+X_0)(1+X)^2}{2(1+X)^2\big[1-(1+X)^2\big]^{3/2}}\,U-\frac{U^{3/2}}{3(1+X)^3}+O(U^2).
\end{equation}
Also for small $U$, $F(X,t)$ in (\ref{2_th3_1}) becomes
\begin{equation}\label{4_F(X,t)}
F=F_c+\frac{1+X}{4}\sqrt{U}+\frac{(1+X_0)(1+X)^2-1}{8\big[1-(1+X)^2\big]^{3/2}}\,U+O(U^2)=F_c+\frac{(1+X)^2}{4}(t_*-t)-\frac{U^{3/2}}{12(1+X)}+O(U^2),
\end{equation}
where we used (\ref{4_t*expandU}). But $\sqrt{U}\sim(1+X)(t_*-t)$ so that $U^{3/2}\sim (t_*-t)^3(1+X)^3$ and (\ref{4_F(X,t)}) becomes the exponential part of (\ref{2_th3_1t*}). The algebraic part of (\ref{2_th3_1t*}) corresponds simply to $G(X,t_*(X;X_0))$. Thus (\ref{2_th3_1t*}) is the smooth limiting form of the Region I approximation as $t\uparrow t_*$, and a completely analogous calculation shows that (\ref{2_th3_1t*}) is also the limiting form (\ref{2_th3_2}) as $t\downarrow t_*$. To go from (\ref{2_th3_1}) to (\ref{2_th3_2}) we can simply replace $\sqrt{(1+X)^2+4\phi_s}$ by $-\sqrt{(1+X)^2+4\phi_s}$, and $(1+X)^2+4\phi_s$ has, in view of (\ref{4_epsilon0}), a double zero at $t=t_*$. 

Finally, we consider the vicinity of the cusp point $(X_{cusp},t_{cusp})$, defined by (\ref{2_Xcusp}) and (\ref{2_tcusp}). Then certainly $t>t_*$, and at the cusp the Regions II, V, VI and VII meet. We evaluate the sum in (\ref{4_sum_regII}), with (\ref{4_sumappxf}) and (\ref{4_sumappxg}). Now the maxima of $\widetilde{f}$ at $z_1$ and $z_3$ are very close, and the estimate in (\ref{4_reIV}) no longer holds. Note that at the cusp point $z_+=z_-$ and both $\widetilde{f}'$ and $\widetilde{f}''$ vanish.

Again applying the Laplace method to (\ref{4_sum_regII}), we expand the summand about
\begin{equation}\label{4_zcusp}
z=z_{cusp}\equiv\frac{3}{4(X_0+3)},
\end{equation}
and we note that (\ref{4_zcusp}) agrees with (\ref{4_zatcusp}). Let us now write $\widetilde{f}(z;X,t)$ and $\widetilde{g}(z;X)$, to indicate the dependence on all three variables in (\ref{4_sumappxf}) and (\ref{4_sumappxg}). Then, after some calculation, 
\begin{equation}\label{4_gtildecusp}
\widetilde{g}(z_{cusp};X_{cusp})=\frac{1}{2\pi}\sqrt{\frac{X_0+3}{X_0}}\sqrt{\sqrt{X_0(X_0+3)}+X_0+2}
\end{equation}
and
\begin{equation*}\label{4_ftildecusp}
\widetilde{f}(z_{cusp};X_{cusp},t_{cusp})=\frac{\sqrt{X_0}}{4(X_0+4)}\Big[-2(X_0+3)^{3/2}+\sqrt{X_0}\,(2X_0+9)\Big].
\end{equation*}
We shall need some of the higher order terms in the Taylor expansion of  $\widetilde{f}(z;X,t)$ about $(z,X,t)=(z_{cusp},X_{cusp},t_{cusp})$.

From (\ref{4_tildef'}) we find that
$\widetilde{f}'(z_{cusp})\sim t_{cusp}-t+\sqrt{{(X_0+4)(X_0+3)}/{X_0}}(X-X_{cusp})$,
where we first set $z=z_{cusp}$ in (\ref{4_tildef'}) and then expanded the result for $(X,t)\to (X_{cusp},t_{cusp})$. Also, from (\ref{4_tildef''}) we have
$\widetilde{f}''(z_{cusp})\sim 2\big[{(X_0+4)(X_0+3)}/{X_0}\big]^{3/2}\,(X-X_{cusp})$,
$\widetilde{f}'''(z_{cusp})=O(X-X_{cusp})$ as $X\to X_{cusp}$, and
$\widetilde{f}^{(iv)}(z_{cusp})\sim -{16(X_0+3)^{11/2}}/{X_0^{5/2}}$, as $X\to X_{cusp}.$
If we introduce the scaled variables $(\xi,\eta)$ in (\ref{2_th3_6XT}) and also set
$z-z_{cusp}=\varepsilon^{1/4}\,w=O(\varepsilon^{1/4})$,
then
$$\frac{\widetilde{f}^{(iv)}(z_{cusp})}{4!\varepsilon}\big(z-z_{cusp}\big)^4\sim-\frac{2(X_0+3)^{11/2}}{3X_0^{5/2}}\,w^4=O(1),\quad\frac{\widetilde{f}'''(z_{cusp})}{3!\varepsilon}\big(z-z_{cusp}\big)^3=O(\varepsilon^{1/4}),$$
and
$$\frac{\widetilde{f}''(z_{cusp})}{2!\varepsilon}\big(z-z_{cusp}\big)^2\sim\bigg[\frac{(X_0+4)(X_0+3)}{X_0}\bigg]^{3/2}\,\xi^2,\quad\frac{\widetilde{f}'(z_{cusp})}{\varepsilon}\big(z-z_{cusp}\big)\sim -\eta\,w.$$
We thus have the approximation
\begin{equation}\label{4_ftilde/ep}
\frac{\widetilde{f}(z)}{\varepsilon}=\frac{\widetilde{f}(z;X,t)}{\varepsilon}=\frac{\widetilde{f}_0(X,t)}{\varepsilon}-\eta\,w+\bigg[\frac{(X_0+4)(X_0+3)}{X_0}\bigg]^{3/2}\,\xi\,w^2-\frac{2(X_0+3)^{11/2}}{3X_0^{5/2}}\,w^4+O(\varepsilon^{1/4}),
\end{equation}
where
$\widetilde{f}_0(X,t)=\widetilde{f}(z_{cusp};X,t)$
is obtained by setting $z=z_{cusp}$ in (\ref{4_sumappxf}). Then we can further expand $\widetilde{f}_0(X,t)$ about $(X,t)=(X_{cusp},t_{cusp})$ in a double Taylor series, and this ultimately leads to 
\begin{eqnarray}\label{4_ftilde0/ep}
\frac{\widetilde{f}_0(X,t)}{\varepsilon}&=&\frac{\widetilde{f}(z_{cusp};X_{cusp},t_{cusp})}{\varepsilon}-\frac{\big(\sqrt{X_0}+2\sqrt{X_0+3}\,\big)^2}{4\sqrt{X_0(X_0+3)(X_0+4)}}\frac{\xi}{\sqrt{\varepsilon}}\nonumber\\
&&-\frac{3}{4(X_0+3)}\,\frac{\eta}{\varepsilon^{1/4}}-\frac{\sqrt{X_0}+2\sqrt{X_0+3}}{4\sqrt{X_0}}\,{\xi^2}+o(1)
\end{eqnarray}
where we again used the scaled variables $(\xi,\eta)$ in (\ref{2_th3_6XT}). Using (\ref{4_ftilde0/ep}) in (\ref{4_ftilde/ep}), and (\ref{4_gtildecusp}), we obtain the leading order approximation to the summand in (\ref{4_sum_regII}) for $z-z_{cusp}=O(\varepsilon^{1/4})$ and $(X,t)$ near the cusp. The sum is asymptotically the same as $\varepsilon^{-1}$ times an integral over $z$, or $\varepsilon^{-3/4}$ times an integral over $w\in(-\infty,\infty)$. We thus obtain the approximation in (\ref{2_th3_6}) and (\ref{2_th3_6J}), that involves the integral $\mathcal{J}$. This completes the analysis of the range $X<0$.

\section{Alternate approach}

We discuss a different approach to the asymptotics, based on geometrical optics and singular perturbations, which will lead to results equivalent to those in Theorems \ref{theorem2} and \ref{theorem3}.

Scaling $x=\beta X$ in (\ref{2_pde}) and using $\varepsilon=\beta^{-2}$ we argue that an asymptotic solution to (\ref{2_pde}) is in the form
\begin{equation}\label{5_pinpde}
p(x,t)=e^{F(X,t)/\varepsilon}\,\Big[G(X,t)+\varepsilon\,G^{(1)}(X,t)+O(\varepsilon^2)\Big].
\end{equation}
Then we distinguish $X>0$ and $X<0$ by setting $(F,G)=(F^+,G^+)$ for $X>0$ and $(F,G)=(F^-,G^-)$ for $X<0$. From (\ref{5_pinpde}) and (\ref{2_pde}) we find that $F^\pm$ satisfy the ``eikonal" equation(s)
\begin{equation}\label{5_eiconalF+}
F_t^+=\big(F_X^+\big)^2+F_X^+;\quad X,\,t>0,
\end{equation}
\begin{equation}\label{5_eiconalF-}
F_t^-=\big(F_X^-\big)^2+(X+1)F_X^-;\quad X<0,\,t>0,
\end{equation}
while $G^\pm$ satisfy the ``transport" equation(s)
\begin{equation}\label{5_transG+}
G_t^+=\big(2F_X^++1\big)G_X^++F_{XX}^+\,G^+;\quad X,\,t>0,
\end{equation}
\begin{equation}\label{5_transG-}
G_t^-=\big(2F_X^-+X+1\big)G_X^-+(F_{XX}^-+1)\,G^-;\quad X<0,\,t>0.
\end{equation}
The initial condition for $p$ is $p(x,0)=\delta(x-x_0)=\sqrt{\varepsilon}\,\delta(X-X_0),$
and the interface conditions at $X=0$, in view of (\ref{5_pinpde}), imply asymptotically that
\begin{equation}\label{5_ICF}
F^+(0^+,t)=F^-(0^-,t),
\end{equation}
\begin{equation}\label{5_ICG}
G^+(0^+,t)=G^-(0^-,t).
\end{equation}
Then $F_X$ is continuous at $X=0$ automatically, which follows from (\ref{5_eiconalF+}), (\ref{5_eiconalF-}) and (\ref{5_ICF}).

A singular perturbation analysis would first consider the short time scale $t=\varepsilon\,\tau_0=O(\varepsilon)$ and a spatial scale near the initial condition, with $X-X_0=\varepsilon\,\xi=O(\varepsilon)$. Then in terms of $(\xi,\tau_0)$, $p$ would satisfy the PDE $p_{\tau_0}=p_{\xi\xi}+p_\xi$, for $\tau_0>0$ and $-\infty<\xi<\infty$, and the initial condition $p=\varepsilon^{-1/2}\,\delta(\xi)$ at $\tau_0=0$. But then $p\sim\varepsilon^{-1/2}\,(4\pi\tau_0)^{-1/2}\exp\big[-(\xi+\tau_0)^2/(4\tau_0)\big]$, so on the $(\xi,\tau_0)$ scale $p$ may be approximated by the free space Brownian motion with unit negative drift.

On the $(X,t)$ scale we solve the non-linear PDE (\ref{5_eiconalF+}) by the method of characteristics. The appropriate solution, which is necessary to asymptotically match to that on the $(\xi,\tau_0)$ scale, is a singular solution where all characteristics, also called ``rays", start from the point $(X,t)=(X_0,0)$. Since (\ref{5_eiconalF+}) is ``constant-coefficient", the rays are all straight lines, and every point in the quarter plane $\{(X,t):X>0,t>0\}$ is reached by a unique ray from $(X_0,0)$. The solution of (\ref{5_eiconalF+}) corresponding to this ray family has $F^+=-(X-X_0+t)^2/(4t)$. Then we can also solve (\ref{5_transG+}), which is a first order linear PDE, by the method of characteristics. The solution will involve an arbitrary function of $(X-X_0)/t$, which is a function of which ray we are on. By using asymptotic matching between the $(X,t)$ and $(\xi,\tau_0)$ scales we can uniquely determine this unknown function, and ultimately obtain $G^+=1/\sqrt{4\pi t}$. We thus conclude that the free space Brownian motion $p_{_{BM}}$ applies also on the $(X,t)$ scale. But it will not be the leading term for $p$ for all $X>0$, as we show below. This must be true, since, for example, $p_{_{BM}}\to 0$ as $t\to\infty$ while the density $p(x,t)$ approaches the steady state in (\ref{2_steady}) or (\ref{2_steady_appro}).

Next we consider the range $X<0$ and solve (\ref{5_eiconalF-}) and (\ref{5_transG-}) subject to the continuity conditions (\ref{5_ICF}) and (\ref{5_ICG}), the latter yielding 
\begin{equation}\label{5_F-0}
F^-(0,t)=-\frac{1}{4t}(t-X_0)^2,
\end{equation}
\begin{equation}\label{5_G-0}
G^-(0,t)=\frac{1}{2\sqrt{\pi t}}.
\end{equation}
We can study in more detail the interface at $X=0$, by considering the scale $X=O(\varepsilon)$. Then to leading order this will simply regain the conditions in (\ref{5_F-0}) and (\ref{5_G-0}). But the next term(s) in the expansion will show a small $O(\varepsilon)$ probability that a ``direct" ray, starting from $(X_0,0)$, will reflect in the $t$-axis ($X=0$). This small reflection can be used to construct a new ray family, which we call the ``reflected" rays, and this will ultimately yield the right side of (\ref{2_th2_1}) in Theorem \ref{theorem2}. However, as we discussed in section 2, the reflected rays can never be the dominant part of $p(x,t)$. Even along $X=0$ they are smaller than $p_{_{BM}}$ (the direct ray expansion) by a factor $O(\varepsilon)$, which is also evident in (\ref{2_th2_1}). For $X>0$ the reflected rays lead to a solution exponentially smaller than $p_{_{BM}}$, and thus we do not consider the reflected rays further here.

\begin{figure}[h]   
  \begin{minipage}[t]{0.5\linewidth}  
    \centering   
    \includegraphics[width=0.85\textwidth]{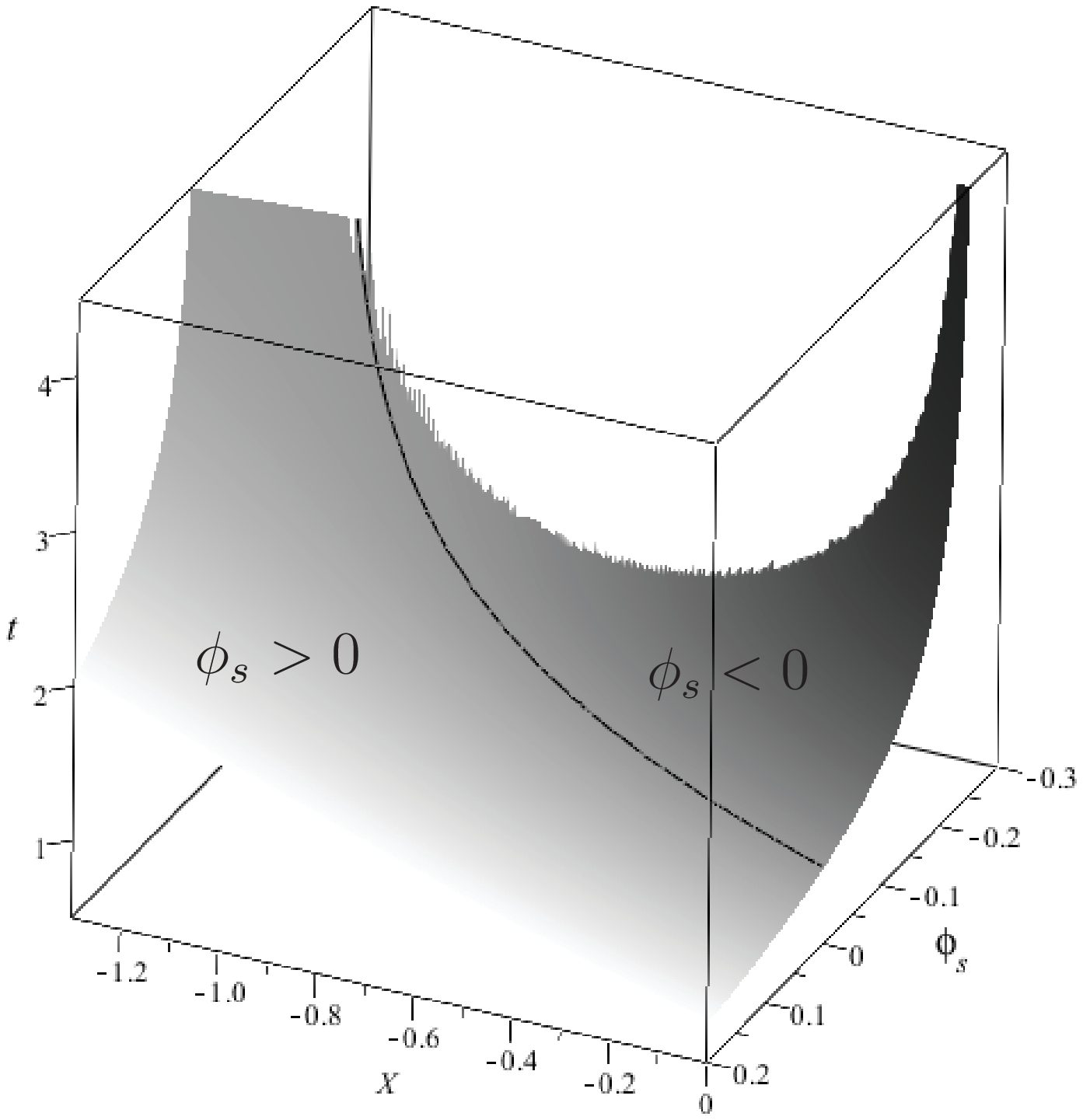}   
    \caption[]{The curve $X+1=e^{X_0-t}$ divides the\\surface $\phi_s(X,t)$ into regions where $\phi_s>0$ and\\$\phi_s<0$.}   
    \label{figure5}   
  \end{minipage}%
  \begin{minipage}[t]{0.5\linewidth}   
    \centering   
    \includegraphics[width=0.75\textwidth]{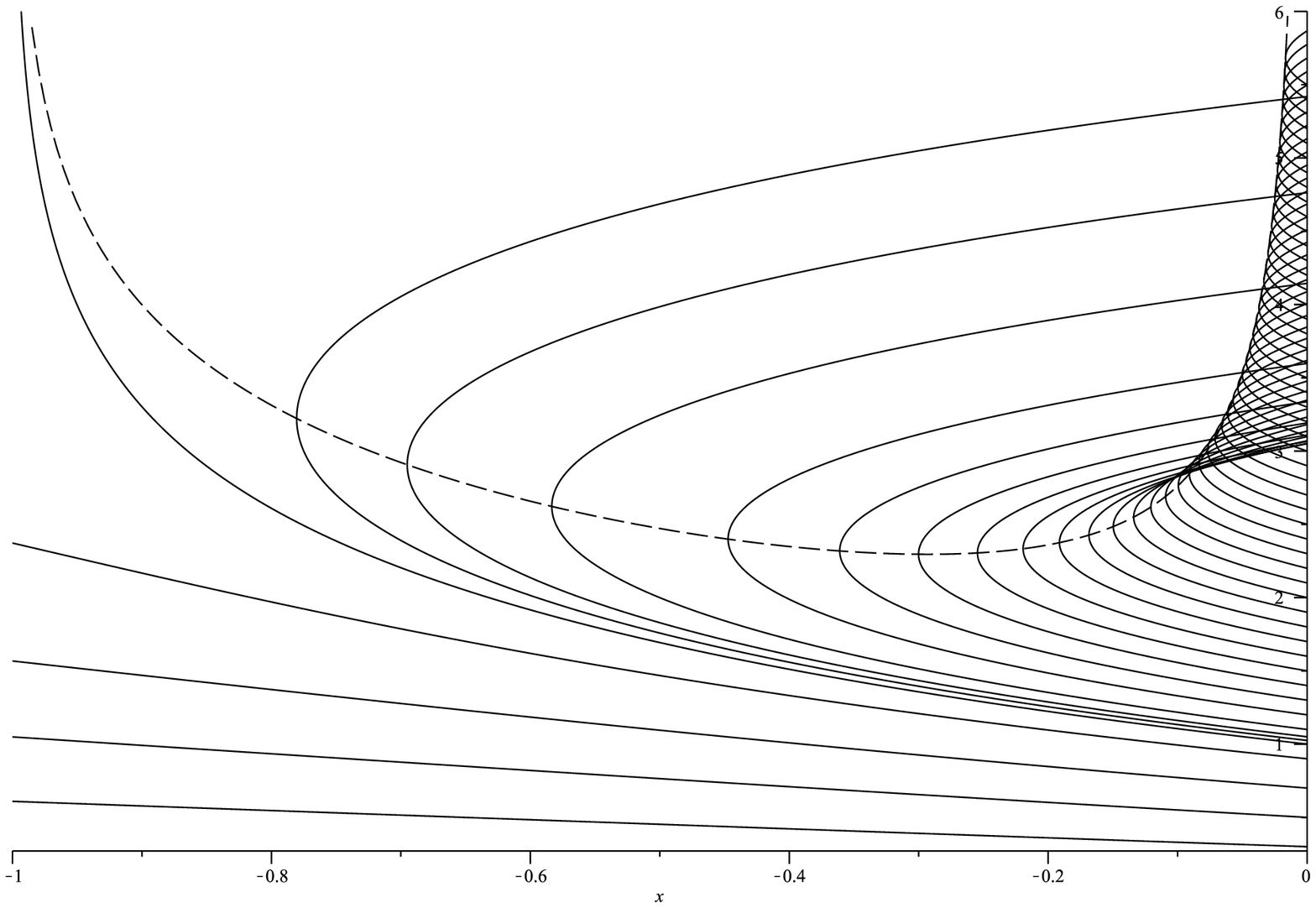}   
    \caption{Rays for $X<0$ and $X_0=1$ (the dashed curve is $t_*(X,X_0)$).}   
    \label{figure6}   
  \end{minipage}   
\end{figure}

We solve (\ref{5_eiconalF-}) subject to (\ref{5_F-0}), so that $X=0$ provides the ``initial manifold" for the first order PDE, and hence the solution for $X<0$ is not a singular solution. However, the solution will have certain ``singular" features, as we show below. Omitting the details, the solution of (\ref{5_eiconalF-}) and (\ref{5_F-0}) leads to the rays
\begin{equation}\label{5_rayt}
t=\alpha+\tau,
\end{equation}
\begin{equation}\label{5_rayX}
X=\bigg(\frac{1}{2}-\frac{X_0}{2\alpha}\bigg)\,e^{\tau}+\bigg(\frac{1}{2}+\frac{X_0}{2\alpha}\bigg)\,e^{-\tau}-1.
\end{equation}
Thus (\ref{5_rayt}) and (\ref{5_rayX}) give the mapping between the original $(X,t)$ variables and the $(\alpha,\tau)$ ray variables. When $\tau=0$, $X=0$ and thus $\alpha$ is the value of $t$ where the given ray hits the initial manifold. The solution $F^-$ can be expressed in terms of the ray variables as
\begin{equation}\label{5_rayF-}
F^-=-\frac{1}{2}\bigg(\frac{1}{2}-\frac{X_0}{2\alpha}\bigg)^2(e^{2\tau}-1)-\frac{\alpha}{4}+\frac{X_0}{2}-\frac{X_0^2}{4\alpha},
\end{equation}
and we note that when $\tau=0$, $\alpha=t$ and thus (\ref{5_F-0}) is satisfied. We can set $\tau=t-\alpha$ in (\ref{5_rayX}) and then this equation defines a one parameter family of curves in the $(X,t)$ plane, for $\alpha>0$, with $\alpha$ indexing the family. In Figure \ref{figure6} we sketch the rays in the space-time plane for $X<0$. The figure suggests that some rays maintain a negative slope ${dt}/{dX}$ for all $\tau>0$, while others change their slope and ultimately return to $X=0$ at some time $t>\alpha$. Furthermore, while some regions of the $(X,t)$ plane are reached by a single ray, in other regions the rays intersect and then $(X,t)$ corresponds to more than one value of $(\alpha,\tau)$ (the figure suggests that in the multi-valued range, exactly 3 rays reach a given $(X,t)$). Below we establish all these results analytically. 

Setting $\tau=t-\alpha$ in (\ref{5_rayX}) with $X=0$ and solving for $t$ leads to
\begin{equation}\label{5_t-alpha}
t-\alpha=\log\bigg(\frac{\alpha\pm X_0}{\alpha- X_0}\bigg).
\end{equation} 
Certainly $t=\alpha$ corresponds to $X=0$ for all $\alpha>0$, which corresponds to the $(-)$ sign in (\ref{5_t-alpha}), but the $(+)$ sign gives a second solution if $\alpha>X_0$. Thus for $\alpha<X_0$ a ray will not return to $X=0$, but if $\alpha>X_0$ the ray will return at time
\begin{equation}\label{5_t-alpha+}
t=\alpha+\log\bigg(\frac{\alpha+ X_0}{\alpha- X_0}\bigg),\quad \alpha>X_0.
\end{equation}
This return time becomes infinite as $\alpha\downarrow X_0$, but becomes very close to $\alpha$ if $\alpha\to +\infty$, which corresponds to an almost immediate return to $X=0$. 

When $\alpha=X_0$ we have $X=e^{-\tau}-1=e^{X_0-t}-1$, which is precisely the fluid approximation in the range $X<0$ (or $t>X_0$). Thus all rays that begin from $X=0$ at times before the fluid approximation never return to $X=0$ and maintain a negative slope, while all rays that begin at times $t>X_0$ return to $X=0$. Also, from (\ref{5_rayX}), if $\alpha<X_0$, a ray will have $X\to -\infty$ as $t\to \infty$, while if $\alpha=X_0$, the fluid approximation $\to -1$ as $t\to\infty$. Note also that $F^-=0$ if $\alpha=X_0$.

It is also instructive to find the minimum value of $t$ where the rays return to $X=0$. When $t=X_0^+$ the return value is very large. To find the minimum value we minimize the right side of (\ref{5_t-alpha+}), which yields $\alpha=\sqrt{X_0(X_0+2)}>0$ and thus the returned rays hit $X=0$ for times $t\ge t_{\min}(X_0)$, with 
\begin{equation}\label{5_tmin}
t_{\min}(X_0)=\sqrt{X_0(X_0+2)}+2\log\bigg(\frac{\sqrt{X_0+2}+\sqrt{X_0}}{\sqrt{2}}\bigg).
\end{equation}
Thus for $t<t_{\min}$, $X=0$ has one ray leaving and for $t>t_{\min}$ one ray leaves and at least one returns.

Next we consider the geometric envelope(s) of the rays, which Figure \ref{figure6} suggests will play a role in the analysis. Writing the rays as
\begin{equation}\label{5_rayXalpha}
X=\mathcal{F}(t,\alpha)=\bigg(\frac{1}{2}+\frac{X_0}{2\alpha}\bigg)\,e^{\alpha-t}+\bigg(\frac{1}{2}-\frac{X_0}{2\alpha}\bigg)\,e^{t-\alpha}-1,
\end{equation}
the envelope is obtained by setting $\partial\mathcal{F}/\partial\alpha=0$, which yields
\begin{equation}\label{5_tenve}
t-\alpha=\frac{1}{2}\log\bigg(\frac{\alpha^2+\alpha\,X_0-X_0}{\alpha^2-\alpha\,X_0-X_0}\bigg).
\end{equation}
Setting (\ref{5_tenve}) in (\ref{5_rayXalpha}) and performing some algebraic simplifications leads to the quartic equation
\begin{equation}\label{5_quartic}
\big[(X+1)^2-1\big]\alpha^4-\big[(2X_0+X_0^2)(X+1)^2-2X_0-2X_0^2\big]\alpha^2+X_0^2(X+X_0+2)(X-X_0)=0,
\end{equation}
which is also a quadratic equation for $\alpha^2$. Solving (\ref{5_quartic}) yields
\begin{equation}\label{5_alphapm}
\alpha^2_\pm=\frac{X_0}{2|X|(X+2)}\Big[(X_0+2)|X|(X+2)+X_0\pm\sqrt{X_0}(X+1)\sqrt{(X_0+4)(X+1)^2-4}\Big].
\end{equation}
The expression in (\ref{5_alphapm}) is real only for $X>-1+2/\sqrt{X_0+4}=X_{cusp}$. Then $\alpha_+$ corresponds to $\alpha_c$ in (\ref{2_alpha_c}), and $\alpha_-$ to $\alpha_d$ in (\ref{2_alpha_d}). The envelope(s) of the rays can then be obtained from (\ref{5_tenve}), with $\alpha_c$ corresponding to the ``upper caustic" $t=t_c(X;X_0)$ in (\ref{2_tc}) and $\alpha_d$ to the ``lower caustic" $t=t_d(X;X_0)$ in (\ref{2_td}). An alternate way of obtaining the caustic curves is to examine where the Jacobian between $(X,t)$ and $(\alpha,\tau)$ coordinates vanishes (or is infinite). Using (\ref{5_rayt}) and (\ref{5_rayX}) to compute the Jacobian $\partial(X,t)/\partial(\alpha,\tau)$ again leads to the two caustic curves $t_c$ and $t_d$, which are explicit functions of $X$ and $X_0$. When $X=X_{cusp}$, $\alpha_c=\alpha_d=\sqrt{X_0(X_0+3)}\equiv \alpha_{cusp}$ (which exceeds both $X_0$ and $\sqrt{X_0(X_0+2)}$, the minimum $\alpha$ needed for a ray to return to $X=0$) and then $t_d=t_c=t_{cusp}$, as in (\ref{2_tcusp}). At $X=X_{cusp}$ the two caustic curves have the same slope, ${dt}/{dX}=\sqrt{(X_0+3)(X_0+4)/X_0}$, and thus form a cusp. Also, $t_c\to +\infty$ as $X\to 0^-$, but $t_d$ approaches a finite value, with $t_d(0;X_0)=t_{\min}(X_0)$ as in (\ref{5_tmin}). When $X\to 0^-$, $\alpha_c\sim X_0/\sqrt{-2X}\to +\infty$ while $\alpha_d\to\sqrt{X_0(X_0+2)}$. 

We are now in a position to relate the present ray expansion results to those in Theorem \ref{theorem3}. First consider $\alpha <X_0$. Then we write $t$ in terms of $X$ and $\alpha$ as
\begin{equation}\label{5_talphaTh3}
t=\alpha+\log\bigg[\frac{\sqrt{(X+1)^2-1+X_0^2/\alpha^2}-X-1}{X_0/\alpha-1}\bigg]
\end{equation}
and then $F^-$ in (\ref{5_rayF-}) becomes
\begin{equation}\label{5_F-th3}
F^-=-\frac{1}{4}\bigg[(X+1)^2+\frac{X_0}{\alpha}-1-(X+1)\sqrt{(X+1)^2-1+\frac{X_0^2}{\alpha^2}}\;\bigg]-\frac{X_0^2}{4\alpha}-\frac{\alpha}{4}+\frac{X_0}{2}.
\end{equation}
Comparing (\ref{5_F-th3}) with $F(X,t)$ in (\ref{2_th3_1F}) we see they agree precisely if $\alpha$ and $\phi_s$ (the saddle location) are related by
$\alpha={X_0}/{\sqrt{1+4\phi_s}}$.
Then also (\ref{5_talphaTh3}) is equivalent to (\ref{2_lm1}). We have thus established the equivalence of the ray and saddle point approaches, at least for the exponential parts of the approximation, and for $\alpha<X_0$, which corresponds to $0<t<X_0-\log(1+X)$ if $X\in(-1,0)$, and to all $t$ if $X\le -1$.

Now take $\alpha>X_0$. Then we solve again (\ref{5_rayXalpha}) for $t$ in terms of $X$ and $\alpha$, but now must consider both branches of the quadratic, hence
\begin{equation}\label{5_talpha-}
t=\alpha+\log\bigg[\frac{X+1-\sqrt{(X+1)^2-1+X_0^2/\alpha^2}}{1-X_0/\alpha}\bigg]
\end{equation}
or
\begin{equation}\label{5_talpha+}
t=\alpha+\log\bigg[\frac{X+1+\sqrt{(X+1)^2-1+X_0^2/\alpha^2}}{1-X_0/\alpha}\bigg].
\end{equation}
When (\ref{5_talpha-}) holds, we differentiate along a ray (thus with $\alpha$ fixed) to get
\begin{equation}\label{5_slope<0}
\frac{dt}{dX}=-\bigg[(X+1)^2-1+\frac{X_0^2}{\alpha^2}\bigg]^{-1/2}<0.
\end{equation}
If (\ref{5_talpha+}) holds then $dt/dX$ is the negative of (\ref{5_slope<0}) and thus $dt/dX>0$. Then along a ray $dt/dX$ changes sign (or $dX/dt=0$) when $\alpha=X_0/\sqrt{1-(1+X)^2}$. But then (\ref{5_talpha-}) and (\ref{5_talpha+}) lead to precisely the curve $t_*=t_*(X;X_0)$ in (\ref{2_t*}). This curve is not itself a ray, but represents the locus of the minimal values of $X$ (or maximal $|X|$) that the rays attain, before changing slope and returning to $X=0$ at the time in (\ref{5_t-alpha+}). This is also illustrated in Figure \ref{figure6} (the dashed curve).

Now assume that $X_0-\log(1+X)<t<t_*$, with $X\in(-1,0)$, and that we are in the range where $dt/dX<0$. We define
\begin{equation}\label{5_R1}
R_1(\alpha)=t-\alpha-\log\bigg[\frac{X+1-\sqrt{(X+1)^2-1+X_0^2/\alpha^2}}{1-X_0/\alpha}\bigg]
\end{equation}
and show that $R_1(\alpha)$ has a unique zero. We have, from (\ref{5_R1}),
\begin{equation*}\label{5_R1'}
R_1'(\alpha)=-1+\frac{X_0}{\alpha^2-X_0^2}-\frac{X_0^2(X+1)}{\alpha(\alpha^2-X_0^2)\sqrt{(X+1)^2-1+X_0^2/\alpha^2}}
\end{equation*}
and we shall show that $R_1'(\alpha)<0$, or
\begin{equation}\label{5_R1'<0}
\frac{X_0^2}{\alpha}(X+1)>(X_0^2+X_0-\alpha^2)\sqrt{(X+1)^2-1+\frac{X_0^2}{\alpha^2}}.
\end{equation}
If $\alpha>\sqrt{X_0(X_0+1)}$ then (\ref{5_R1'<0}) is obviously true. If $X_0<\alpha<\sqrt{X_0(X_0+1)}$ we square both sides of (\ref{5_R1'<0}) and after some algebraic manipulation obtain the equivalent inequality
\begin{equation}\label{5_R1'<0eqiv}
\alpha^3-(X_0^2+X_0)\alpha+X_0^2>0.
\end{equation}
But (\ref{5_R1'<0eqiv}) has roots at $\alpha=X_0$ and $\alpha=\big[-X_0\pm\sqrt{X_0(X_0+4)}\big]/2<X_0$. Hence (\ref{5_R1'<0eqiv}) is true for all $\alpha>X_0$ and thus $R_1'(\alpha)<0$. We also have $R_1(\alpha)\to t-X_0+\log(X+1)>0$ as $\alpha\to X_0^+$
and 
\begin{equation}\label{5_R1end}
R_1\bigg(\frac{X_0}{\sqrt{1-(1+X)^2}}\bigg)=t-t_*<0,
\end{equation}
and thus $R_1(\alpha)$ has a unique zero. Then we again obtain (\ref{5_F-th3}), which is equivalent to the saddle point approximation in (\ref{2_th3_1}). We have now covered Regions I and III via the ray approach.

We shall next analyze Regions II, VI and VII, where $dt/dX>0$ and $t>t_*$. Then $R_1(\alpha)=0$ has no solutions, and we need to analyze the root(s) of 
\begin{equation*}\label{5_R2}
R_2(\alpha)=t-\alpha-\log\bigg[\frac{X+1+\sqrt{(X+1)^2-1+X_0^2/\alpha^2}}{1-X_0/\alpha}\bigg],
\end{equation*}
whose derivative is
\begin{equation*}\label{5_R2'}
R_2'(\alpha)=-1+\frac{X_0}{\alpha^2-X_0^2}+\frac{X_0^2(X+1)}{\alpha(\alpha^2-X_0^2)\sqrt{(X+1)^2-1+X_0^2/\alpha^2}}.
\end{equation*}

For Region II, we shall show that $R_2'(\alpha)>0$, or
\begin{equation}\label{5_R2'>0}
\frac{X_0^2}{\alpha}(X+1)>(\alpha^2-X_0^2-X_0)\sqrt{(X+1)^2-1+\frac{X_0^2}{\alpha^2}}.
\end{equation}
If $X_0<\alpha\le\sqrt{X_0(X_0+1)}$, (\ref{5_R2'>0}) is obviously true. If $\alpha>\sqrt{X_0(X_0+1)}$ we square both sides of (\ref{5_R2'>0}) and obtain the equivalent inequality
\begin{equation}\label{5_R2'>0sq}
\bigg[\frac{X_0^4}{\alpha^2}-(\alpha^2-X_0^2-X_0)^2\bigg](X+1)^2>\bigg(\frac{X_0^2}{\alpha^2}-1\bigg)(\alpha^2-X_0^2-X_0)^2.
\end{equation}
Clearly $X_0^2/\alpha^2-1<0$ so if $X_0^4/\alpha^2-(\alpha^2-X_0^2-X_0)^2\ge 0$, which occurs if $X_0<\alpha\le \big[X_0+\sqrt{X_0(X_0+4)}\big]/2$, then (\ref{5_R2'>0sq}) is true for any $X\in (-1,0)$. If, on the other hand, $X_0^4/\alpha^2-(\alpha^2-X_0^2-X_0)^2< 0$, then we multiply both sides of (\ref{5_R2'>0sq}) by $-1$, changing $>$ to $<$, and obtain
\begin{equation}\label{5_R2X<}
X+1<\sqrt{\frac{(\alpha^2-X_0^2)(\alpha^2-X_0^2-X_0)^2}{\alpha^2(\alpha^2-X_0^2-X_0)^2-X_0^4}}.
\end{equation}
But the inequality in (\ref{5_R2X<}) implies, in view of (\ref{5_quartic}), that $X<X_{cusp}$. Thus $R_2'(\alpha)>0$ in Region II. Then $R_2(\alpha)\to -\infty$ as $\alpha\to X_0^+$ and 
\begin{equation}\label{5_R2end}
R_2\bigg(\frac{X_0}{\sqrt{1-(1+X)^2}}\bigg)=t-t_*>0,
\end{equation}
so that $R_2(\alpha)$ has a unique root for $\alpha>X_0$. 

We next consider Regions VI and VII. We note that for $X\in(X_{cusp},0)$ the equation $R_2'(\alpha)=0$ has roots precisely at $\alpha_c$ and $\alpha_d$ (or $\alpha_\pm$ in (\ref{5_alphapm})). 

For Region VI, where $t>t_c$ and $X\in(X_{cusp},0)$, we have $R_2(\alpha_d)=t-t_d>0$, $R_2(\alpha_c)=t-t_c>0$, $R_2(X_0^+)=-\infty$ and the equality in (\ref{5_R2end}) holds. Thus $R_2(\alpha)$ has a unique root for $\alpha>X_0$.

For Region VII, where $t_*<t<t_d$ and $X\in(X_{cusp},X_*)$, we have  $R_2(X_0^+)=-\infty$ and the equality in (\ref{5_R2end}) holds as usual, but now $R_2(\alpha_d)=t-t_d<0$ and $R_2(\alpha_c)=t-t_c<0$. Again, $R_2(\alpha)$ has a unique root for $\alpha>X_0$.

We thus have shown that in Regions II, VI and VII, $R_2(\alpha)$ has a unique root for $\alpha>X_0$. If we identify
$\alpha={X_0}/{\sqrt{1-4z}}$,
and note that, using (\ref{5_talpha+}) and (\ref{5_rayt})-(\ref{5_rayF-}),
\begin{equation}\label{5_R2F-}
F^-=-\frac{1}{4}\bigg[(X+1)^2+\frac{X_0}{\alpha}-1+(X+1)\sqrt{(X+1)^2-1+\frac{X_0^2}{\alpha^2}}\;\bigg]-\frac{X_0^2}{4\alpha}-\frac{\alpha}{4}+\frac{X_0}{2},
\end{equation}
then $F^-$ agrees with $\widetilde{f}(z_1)$ (cf. (\ref{2_th3_2})) or $\widetilde{f}(z_3)$ (cf. (\ref{2_th3_3})). Here we use the expression below (\ref{2_th3_2}) for $\widetilde{f}(z)$ and the fact that $z_1$ or $z_3$ is the unique solution to (\ref{2_lm2}), which is equivalent to (\ref{5_talpha+}).

Next we examine the interior of the caustics, where $t_d<t<t_c$, so certainly $X\in(X_{cusp},0)$. This corresponds to Regions IV and V. The equation $R_1(\alpha)=0$ has a unique zero in Region IV and no solution in Region V. Note that (\ref{5_R1end}) holds in Region IV. However, $R_2(\alpha)$ will now have multiple roots, and also
$R_2(\alpha_c)=t-t_c<0$, $R_2(\alpha_d)=t-t_d>0$,
$R_2(X_0^+)=-\infty$ and the equality in (\ref{5_R2end}) holds. When $t>t_*$, $t>t_d$ and $t<t_c$ (Region V) we have three sign changes in $R_2(\alpha)$ over the interval $\alpha\in (X_0,X_0/\sqrt{1-(1+X)^2}\;)$, and we denote the three roots by $\alpha_1<\alpha_2<\alpha_3$, where $\alpha_j=\alpha_j(X,t)$. When $t_d<t<t_*$ (Region IV) there are two sign changes and we call the roots $\alpha_1<\alpha_2$, and then let $\alpha_3$ be the unique root of $R_1(\alpha)=0$. For both Regions IV and V, a given $(X,t)$ corresponds to three distinct values of $(\alpha,\tau)$, so that the mapping from ray to $(X,t)$ coordinates in 3-to-1. In the five other regions this mapping is 1-to-1 and only a single ray reaches a point $(X,t)$ that is within Regions I-III, VI or VII.

Omitting the details, given $F^-$ in the single valued case we can easily solve the linear PDE in (\ref{5_transG-}) subject to (\ref{5_G-0}), with the result
\begin{equation}\label{5_G-}
G^-=\sqrt{\frac{X_0}{2\pi}}\;\Big|(\alpha^2+\alpha\,X_0-X_0)\,e^{-2(t-\alpha)}-(\alpha^2-\alpha\,X_0-X_0)\Big|^{-1/2}.
\end{equation}
We can show that (\ref{5_G-}) agrees with $G$ in (\ref{2_th3_1}), and also with $\widetilde{G}$ in (\ref{2_th3_2}) and (\ref{2_th3_3}). Then with (\ref{5_G-}) and (\ref{5_rayF-}) we have the leading order ray approximation $p\sim G^-e^{F^-/\varepsilon}$, and this applies in the union of Regions I-III, VI and VII, including the curves that separate them. We observed already in section 4 that $Ge^{F/\varepsilon}$ in (\ref{2_th3_1}) may be smoothly continued from Regions I and III into II, VI and VII, and in particular is smooth along $t=t_*$. It is only upon trying to eliminate $t=\alpha+\tau$ in (\ref{5_rayF-}) that leads to the two different expressions in (\ref{5_F-th3}) and (\ref{5_R2F-}). Outside the caustic range, the rays with $dt/dX<0$ correspond to the saddle point $\phi_s$ in the integral in (\ref{4_lemma7}), while rays with $dt/dX>0$ correspond to the unique interior maximum of $\widetilde{f}$ in (\ref{4_sumappx}) or (\ref{4_sum_regII}), which arises from the residue sum in (\ref{4_lemma7}).

In Regions IV and V, three rays from $X=0$ reach a given $(X,t)$, so we write 
\begin{equation}\label{5_G123}
p(x,t)\sim G_1^-\,e^{F_1^-/\varepsilon}+G_2^-\,e^{F_2^-/\varepsilon}+G_3^-\,e^{F_3^-/\varepsilon},
\end{equation}
where $F_j^-$ corresponds to $(\alpha_j,\tau_j)$ for $j=1,2,3$. The dominant contribution to $p$ will come from the maximal $F_j^-$. In Region IV, $F_3^-$ arises from rays with $dt/dX<0$ and this term corresponds to $Ge^{F/\varepsilon}$ in (\ref{2_th3_5}). But (\ref{2_th3_5}) has only two terms while (\ref{5_G123}) has three. We recall that in Region IV our analysis of the residue sum (\ref{4_lemma7}) led to the equation $\widetilde{f}'(z)=0$ which had two roots $z_1<z_2$, with $z_1$ (resp. $z_2$) being a maximum (resp. minimum) of $\widetilde{f}(z)$. But the Laplace method would never consider minima, and $\widetilde{f}(z_1)$ in (\ref{2_th3_5}) corresponds to $F_1^-$ in (\ref{5_G123}), while $F_2^-$ would correspond to the minimum of $\widetilde{f}$ at $z=z_2$, since it is also a root of $\widetilde{f}'(z)=0$. Geometrically, for every $(X,t)$ in Region IV, three rays meet, one with $dt/dX<0$ ($F_3^-$) and two with $dt/dX>0$ ($F_1^-$ and $F_2^-$). But $F_1^->F_2^-$ so the middle term in (\ref{5_G123}) can never be dominant asymptotically. Also, the term with $F_2^-$ corresponds to a ray that has been tangent to one of the caustic curves before reaching $(X,t)$, while the $F_1^-$ ray has not been tangent to a caustic. For Region V comparing (\ref{5_G123}) to (\ref{2_th3_4}) we now have the correspondence $(\widetilde{f}(z_1),\widetilde{f}(z_3))\leftrightarrow (F_1^-,F_3^-)$ and the middle term in (\ref{5_G123}) is again negligible, corresponding to the local minimum of $\widetilde{f}$ at $z=z_2$, and the $F_2^-$ ray has been tangent to a caustic before reaching $(X,t)$. In Region V all three rays that reach $(X,t)$ have $dt/dX>0$.

Thus we may drop the middle term in (\ref{5_G123}) as it is exponentially smaller than one of the other two terms. Along the caustic $t=t_d$ we have $F_1^-=F_2^-$ but then the $F_3^-$ term dominates. Along $t=t_c$ we have $F_2^-=F_3^-$ but then $F_1^-$ dominates. We also have $G_1^-$ and $G_2^-$, obtained from (\ref{5_G-}), becoming singular along $t=t_d$, while $G_2^-$ and $G_3^-$ become singular along the caustic $t=t_c$. But the $G_j^-$ corresponding to the dominant term in (\ref{5_G123}) will not be singular along the caustic. However, at the cusp point $(X,t)=(X_{cusp},t_{cusp})$, $F_1^-$, $F_2^-$ and $F_3^-$ all agree, and the corresponding $G_j^-$ all become singular. Thus here a genuine non-uniformity develops, which we analyzed using the Laplace method and obtained (\ref{2_th3_6}). It should also be possible to start directly from the PDE and analyze a vicinity of the cusp, say using the local $(\xi, \eta)$ variables in (\ref{2_th3_6XT}), and thus obtain (\ref{2_th3_6}) more directly. However, this would require significant additional work which we do not attempt here. Deciding whether $F_1^-$ or $F_3^-$ in (\ref{5_G123}) dominates amounts to solving $F_1^-=F_3^-$, and this leads to the implicit curve $t_\Gamma$ which we already discussed, and evaluated numerically. Certainly $t_d<t_\Gamma<t_c$ and $t_\Gamma$ hits $X=0$ and a finite point.

Now we return to the range where $X>0$. We must consider the effects of the rays that returned to $X=0$ after traveling in the range $X<0$, on the range $X>0$. We recall that these rays hit $X=0$ at all times $t\ge t_{\min}$, cf. (\ref{5_tmin}). The rays that return to $X=0$ continue into the range $X>0$ and maintains the same slope $dt/dX>0$ at $X=0^+$ and $X=0^-$. The continuity of ray slopes follows from the continuity of $F_X$ at $X=0$, and this follows from (\ref{5_eiconalF+}) and (\ref{5_eiconalF-}). Let us denote this new ray family for $X>0$ by $F^R$, so $F^R$ will again satisfy the PDE in (\ref{5_eiconalF+}). But now the initial condition for $F^R$ at $X=0$ must be obtained from the value of $F^-(0^-,t)$ for a ray that returned. By differentiating (\ref{5_rayXalpha}) implicitly along a ray, we find that $dt/dX=\alpha/X_0$, at the time in (\ref{5_t-alpha+}) when the ray returns to $X=0$. Then we find that the ``returned" rays correspond to 
\begin{equation}\label{5_returned}
X=\frac{X_0}{\alpha}\bigg[t-\alpha-\log\bigg(\frac{\alpha+X_0}{\alpha-X_0}\bigg)\bigg]\equiv \frac{X_0}{\alpha}\,\tau,
\end{equation}
which is a family of straight lines. Also, (\ref{5_returned}) defines the $\tau$ variable for this new family, with $\tau>0$ and $\alpha>\sqrt{X_0(X_0+2)}$. To obtain $F^R$ we solve (\ref{5_eiconalF-}) subject to (\ref{5_R2F-}) at $X=0$, which yields
\begin{equation*}\label{5_FR0}
F^R(0^+,t)=\frac{X_0}{2}-\frac{\alpha}{4}-\frac{X_0^2}{4\alpha}-\frac{X_0}{2\alpha}.
\end{equation*}
We let $\alpha=X_0/\sqrt{1-4z_*}$ where now $z_*(X,t)$ satisfies
\begin{equation*}\label{5_tz*}
t=\frac{X+X_0}{\sqrt{1-4z_*}}+\log\bigg[\frac{1+\sqrt{1-4z_*}}{1-\sqrt{1-4z_*}}\bigg].
\end{equation*}
Then the solution $F^R$ is given by
\begin{equation*}\label{5_FR}
F^R(X,t)=-\bigg(\frac{1}{2}+\frac{X_0}{2\alpha}\bigg)^2\,\tau+\frac{X_0}{2}-\frac{\alpha}{4}-\frac{X_0^2}{4\alpha}-\frac{X_0}{2\alpha}=\frac{X_0-X}{2}-\frac{X+X_0+1}{2}\sqrt{1-4z_*}-\frac{z_*(X+X_0)}{\sqrt{1-4z_*}},
\end{equation*}
which agrees with $f(z_*)$ in (\ref{2_th2_2f}) (with (\ref{2_th2_2z*})). By solving for $G^R(X,t)$, using continuity with $G^-$ at $X=0$, we can similarly obtain the factor in the right side of (\ref{2_th2_2}) that multiplies $e^{f(z_*)/\varepsilon}$. 

The ray family in (\ref{5_returned}) also has a caustic curve, which can be obtained by setting the derivative with respect to $\alpha$ equal to zero, and then eliminating $\alpha$. This yields the curve
\begin{equation}\label{5_t+}
t=t_+(X;X_0)=\sqrt{(X+X_0)(X+X_0+2)}+2\log\bigg[\frac{\sqrt{X+X_0+2}+\sqrt{X+X_0}}{\sqrt{2}}\bigg]
\end{equation}
and this is the same as that in (\ref{3_t+}). Also, the returned rays fill only the domain $t>t_+$ (and $X>0$) and are absent for $t<t_+$. We have $t_+(0;X_0)=t_d(0;X_0)$ so we may view $t_+$ as a caustic in $X>0$ that was induced by the lower caustic $t_d$ in the range $X<0$. Now consider both the direct ray and returned ray expansions for $X>0$. For $t<t_+$ only the direct rays are present and then $p\sim p_{_{BM}}$. For $t>t_+$ every point $(X,t)$ is reached by three rays, the direct path from $(X_0,0)$, and two ``returned" rays from $X=0$. One returned ray reached $(X,t)$ by first being tangent to the caustic $t_+$, and the other returned ray reached $(X,t)$ before being tangent to the caustic. The former returned rays correspond to the minimum of $f(z)$ in (\ref{3_sum}) at $z_{**}$, which we discussed below (\ref{3_t+}). But again clearly $f(z_{**})<f(z_*)$ so we have, for $t>t_+$, 
\begin{equation}\label{5_pfort>t+}
p(x,t)\sim G^R\,e^{F^R/\varepsilon}+p_{_{BM}}.
\end{equation}
Along $t=t_+$, $p_{_{BM}}$ dominates, but for times sufficiently large the first term in (\ref{5_pfort>t+}) overtakes the direct ray solution. Here $F^R$ is understood to arise from the returned ray that was not tangent to $t_+$ before hitting $(X,t)$. The transition from $p_{_{BM}}$ to the returned ray result occurs precisely when $F^R=-(t+X-X_0)^2/(4t)$, and this corresponds to (\ref{2_t2eqn}) and the curve $t=t_2(X;X_0)$.

\begin{figure}[h]   
  \begin{minipage}[t]{0.5\linewidth}  
    \centering   
    \includegraphics[width=0.75\textwidth]{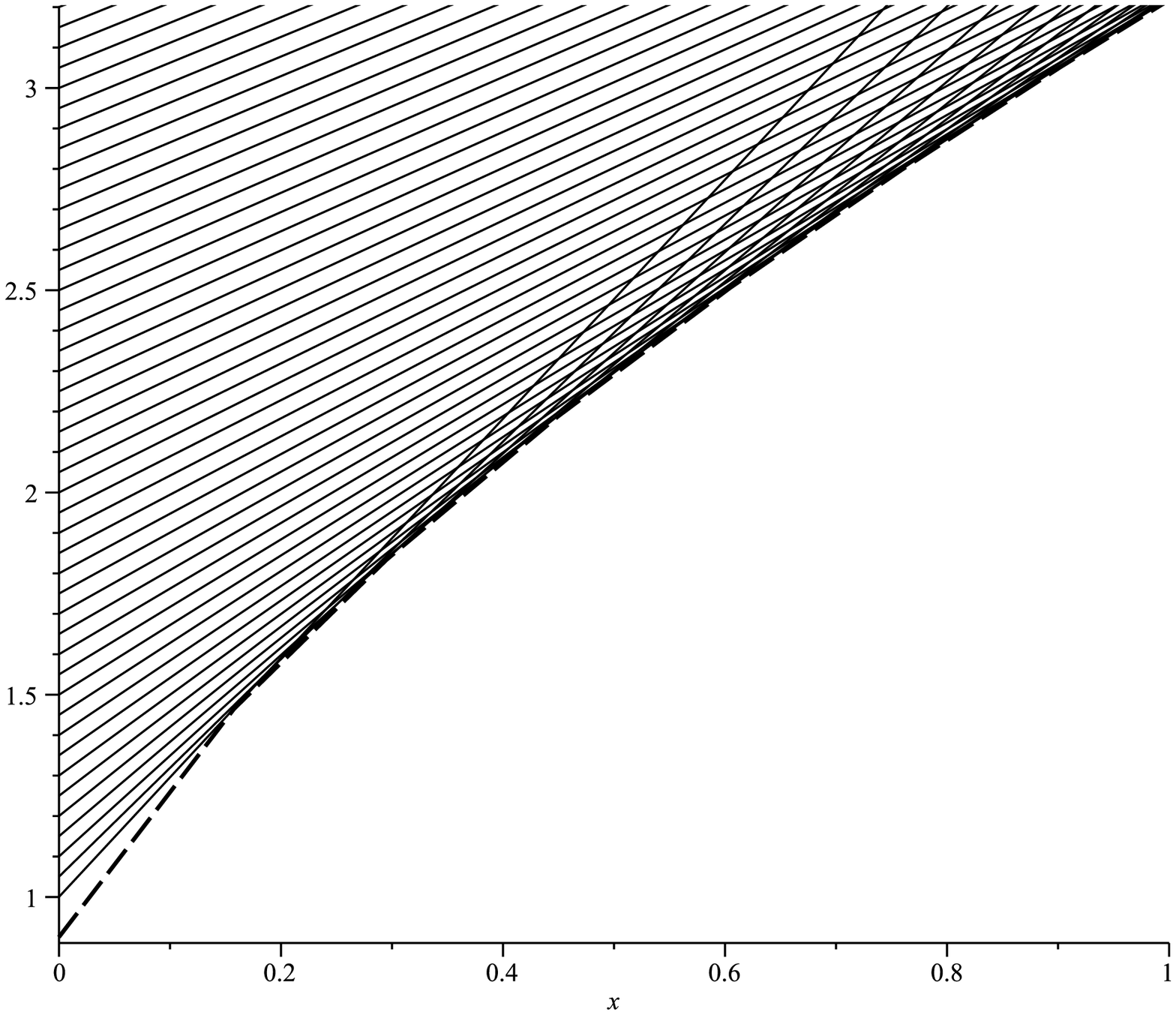}   
    \caption[]{Returned rays in the range $t>t_+$ above\\the caustic with $X_0=0.1$ (the dashed curve is $t_+$).}   
    \label{figure7}   
  \end{minipage}%
  \begin{minipage}[t]{0.5\linewidth}   
    \centering   
    \includegraphics[width=0.9\textwidth]{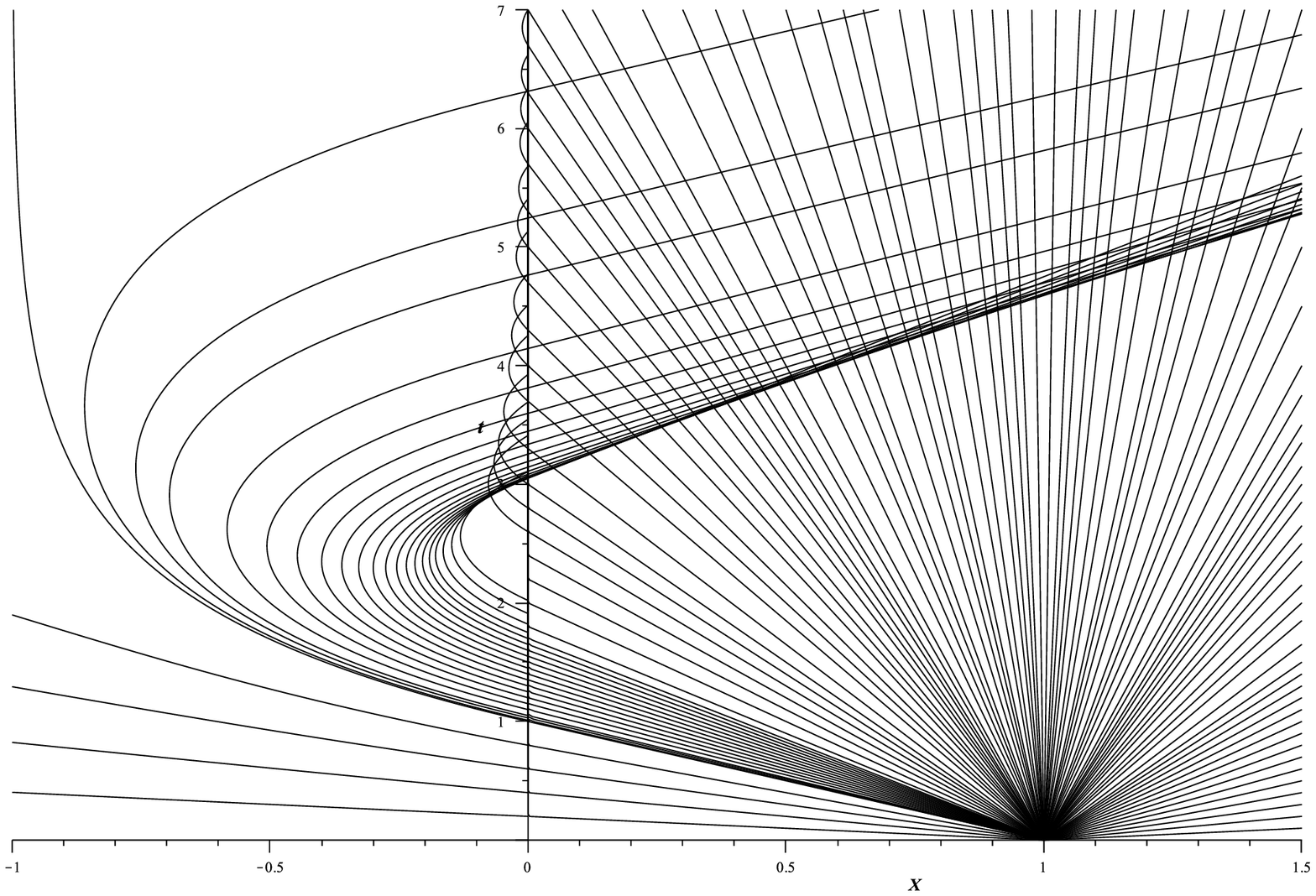}   
    \caption{All the rays for both $X<0$ and $X>0$.}   
    \label{figure8}   
  \end{minipage}   
\end{figure}

In Figure \ref{figure7} we sketch the returned rays in the range $t>t_+$ above the caustic for $X_0=0.1$. Note that for larger values of $X_0$, such as $X_0\ge 1$, (\ref{5_t+}) numerically resembles a straight line (and in fact $t_+\sim X$ as $X\to\infty$), and it is difficult to see the curvature of $t_+$ or the tangency of the rays to this curve.  In Figure \ref{figure8} we sketch all the rays for both $X<0$ and $X>0$, which illustrates the multi-valuedness. Note that numerically part of the curve $t=t_*$ where the rays change slope is hard to distinguish from the upper caustic $t=t_c$. Our analysis here shows that the ray approach is much different from the saddle/singularity analysis, but ultimately yields equivalent results. Either method may be used to obtain higher order terms in (\ref{5_pinpde}).

\newpage
\begin{center}
\begin{tabular} {ccc}
{\includegraphics[width=0.3\textwidth]{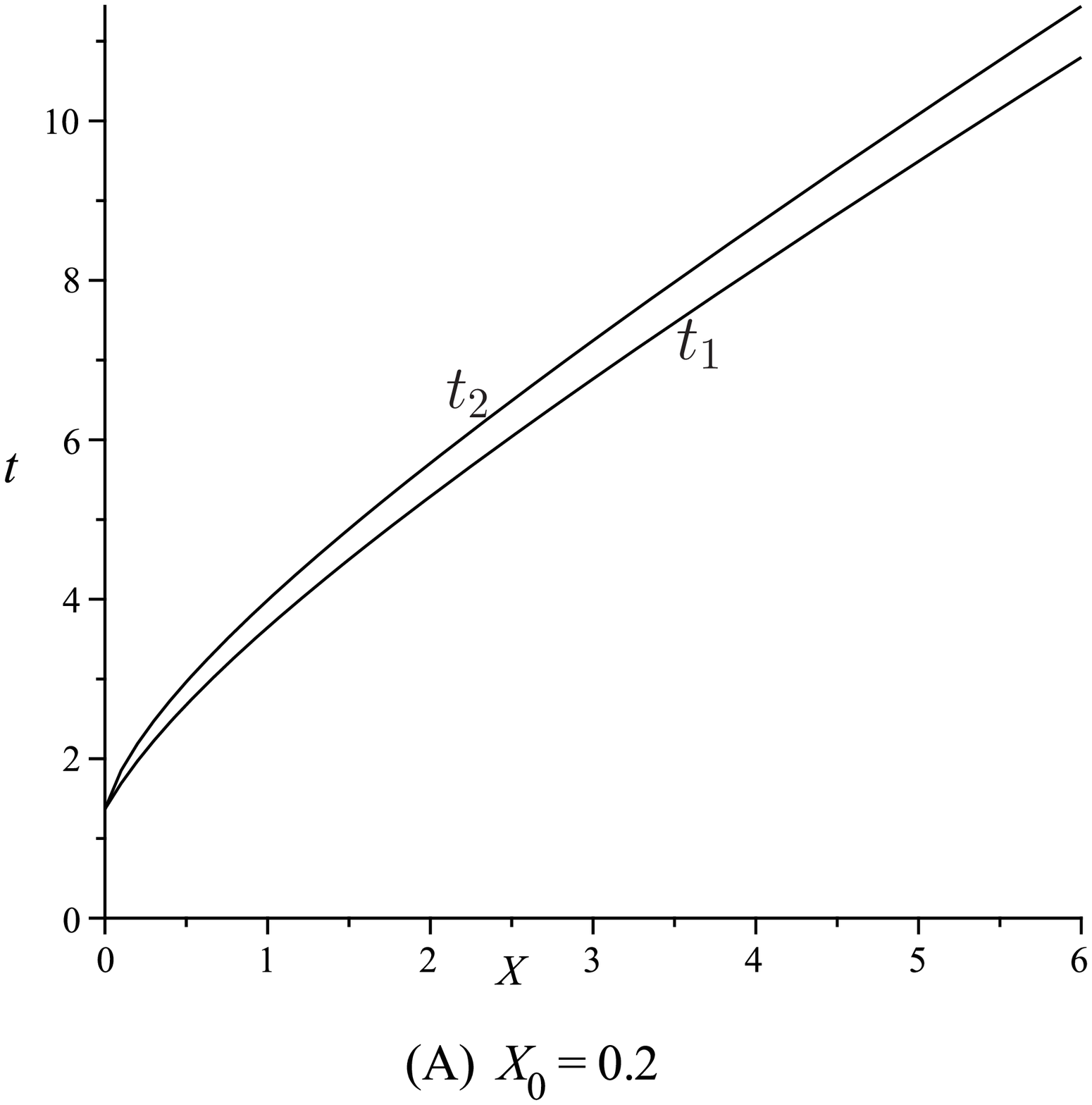}} & {\includegraphics[width=0.3\textwidth]{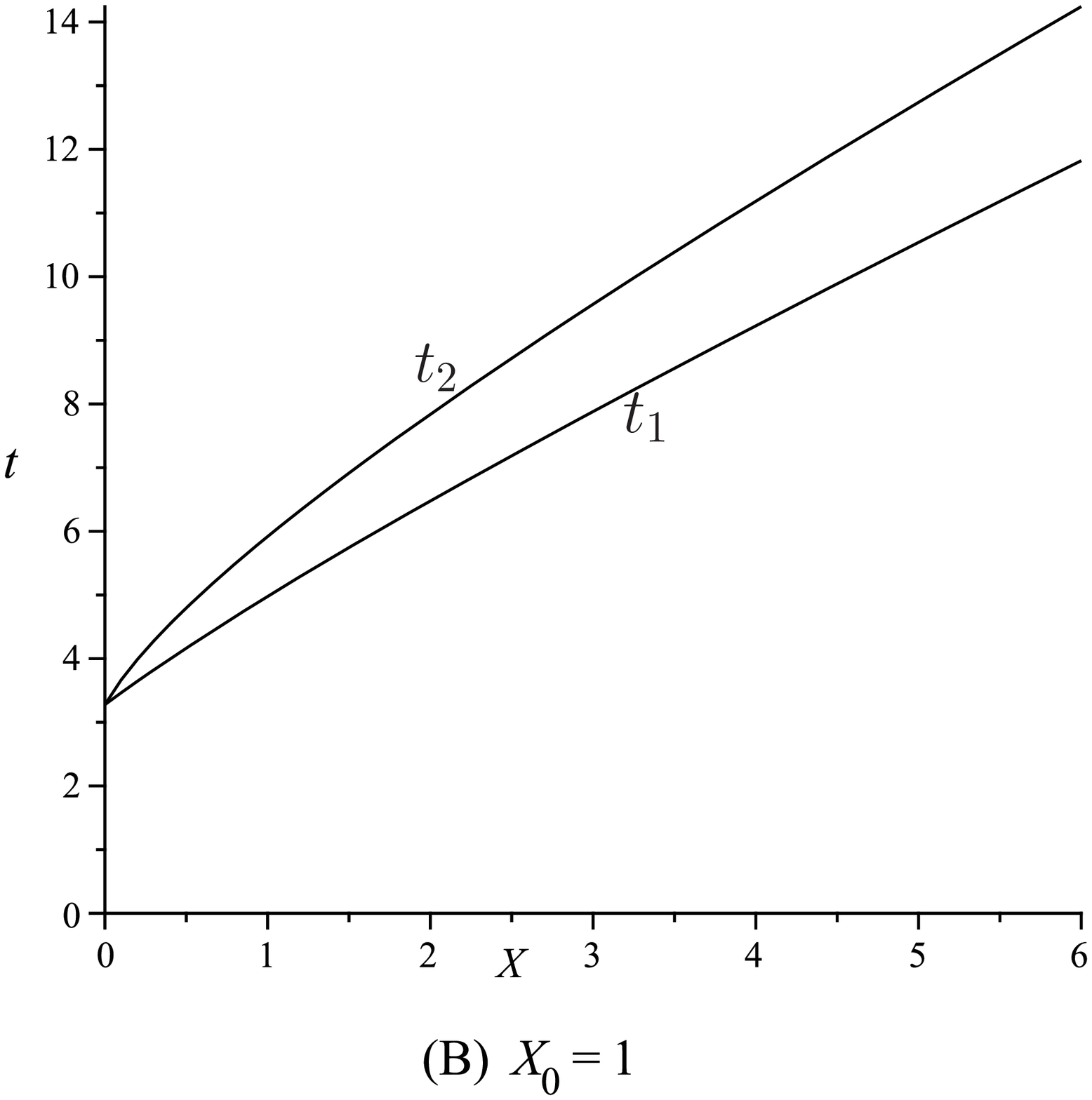}} & {\includegraphics[width=0.3\textwidth]{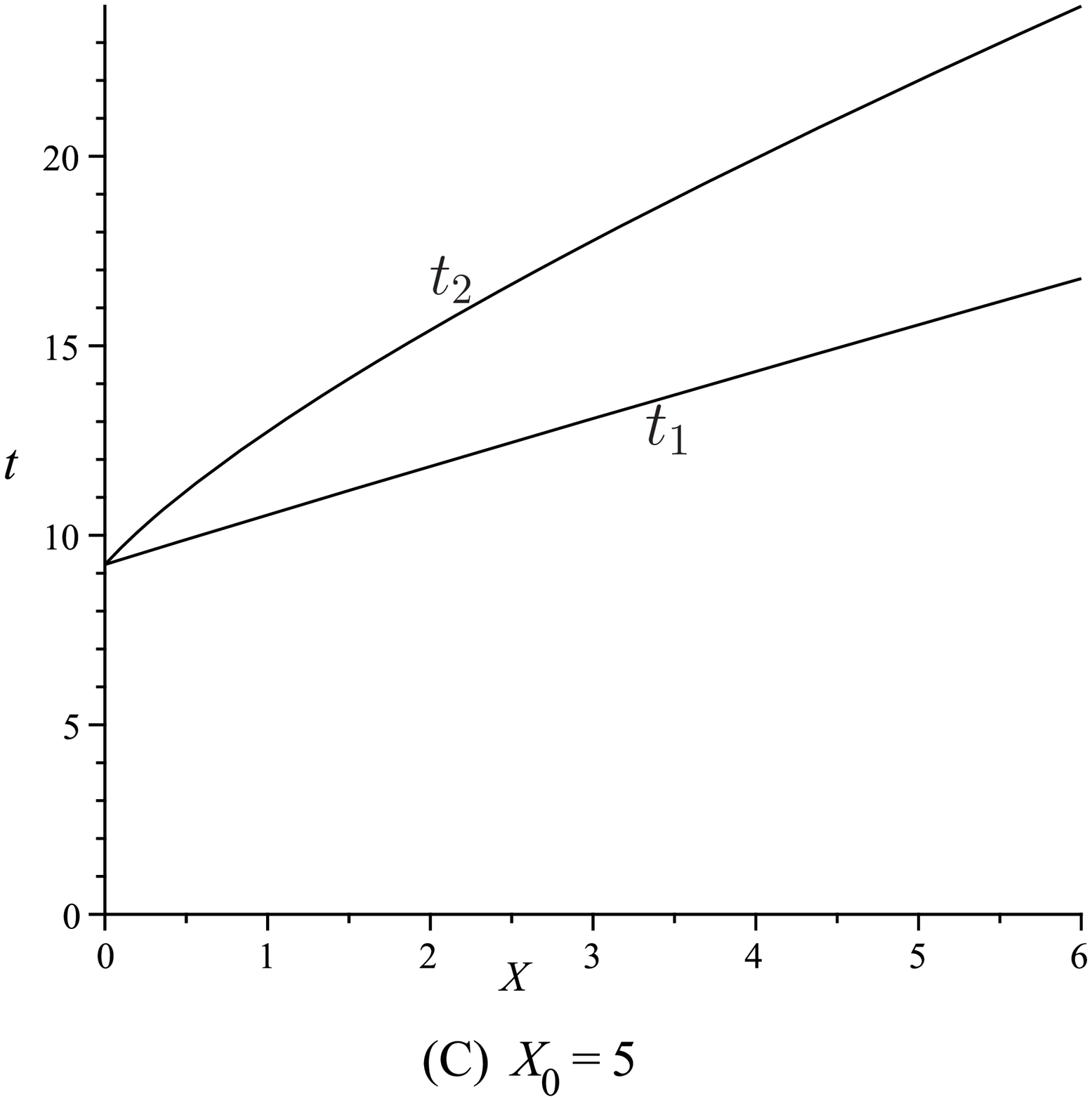}}\\
\end{tabular}
\caption{The transition curves $t_1$ and $t_2$ for $X_0=0.2,\;1,\;5$.}\label{figure1}
\end{center}

\begin{center}
\begin{tabular} {ccc}
{\includegraphics[width=0.3\textwidth]{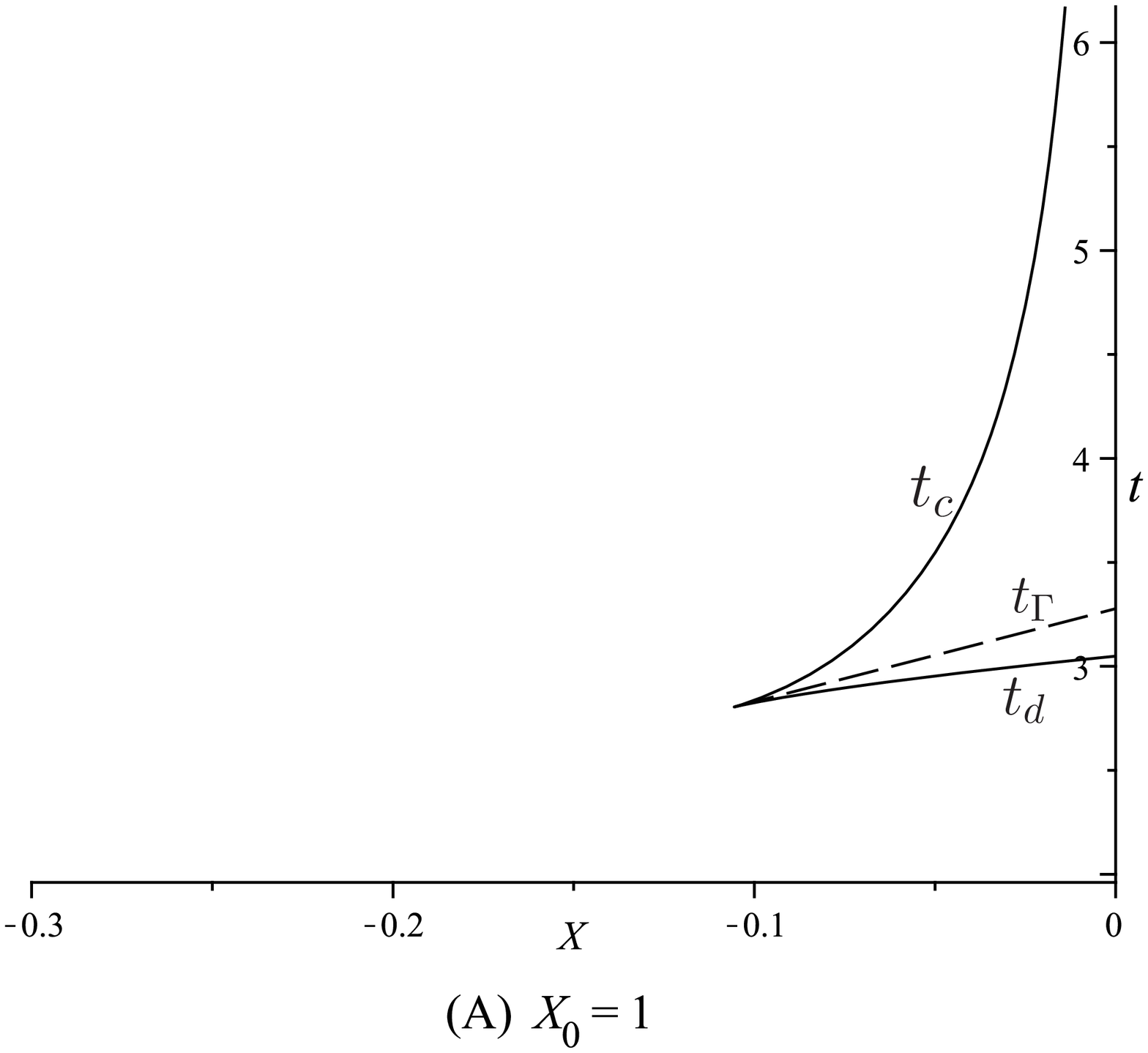}} & {\includegraphics[width=0.3\textwidth]{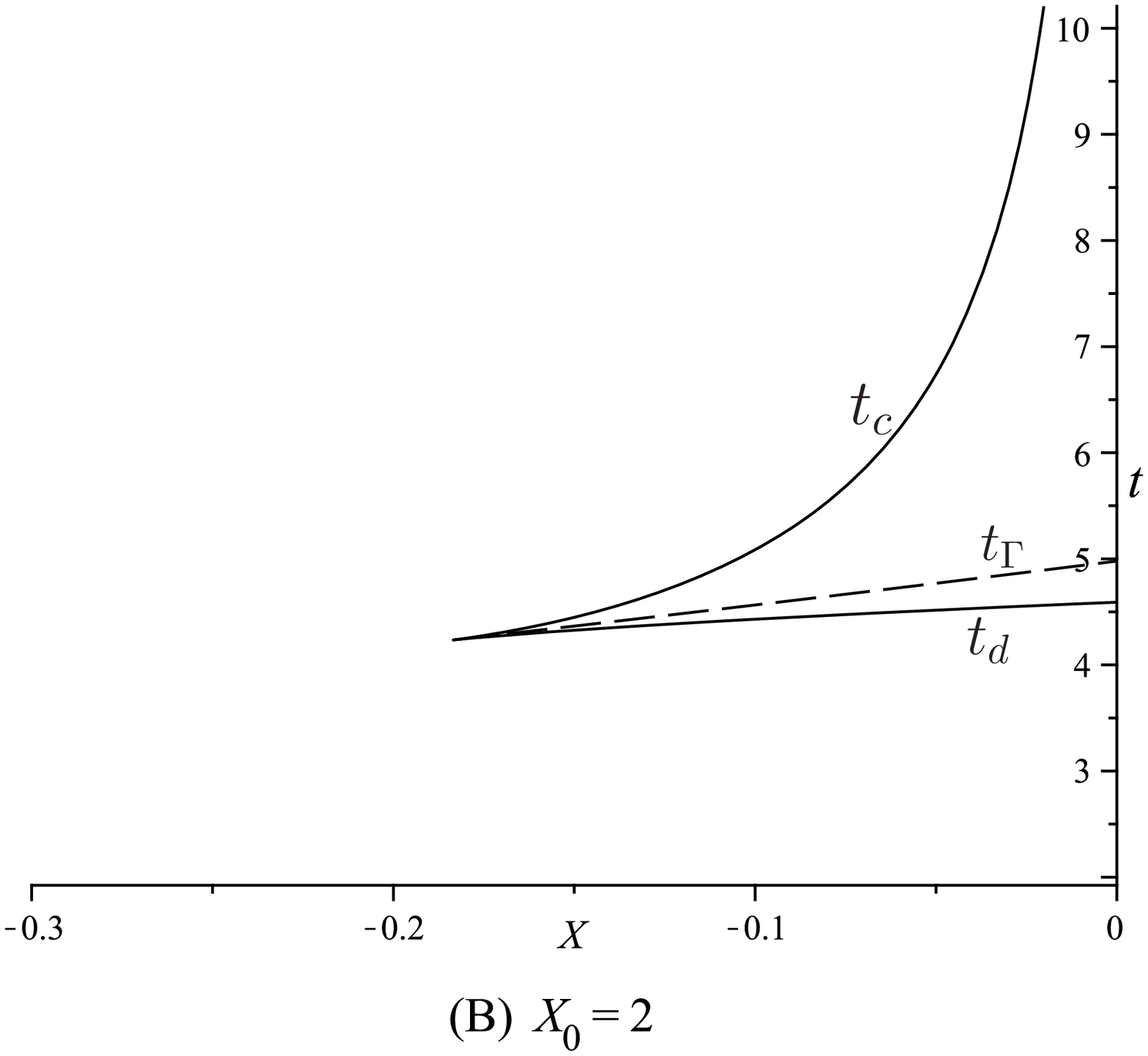}} & {\includegraphics[width=0.3\textwidth]{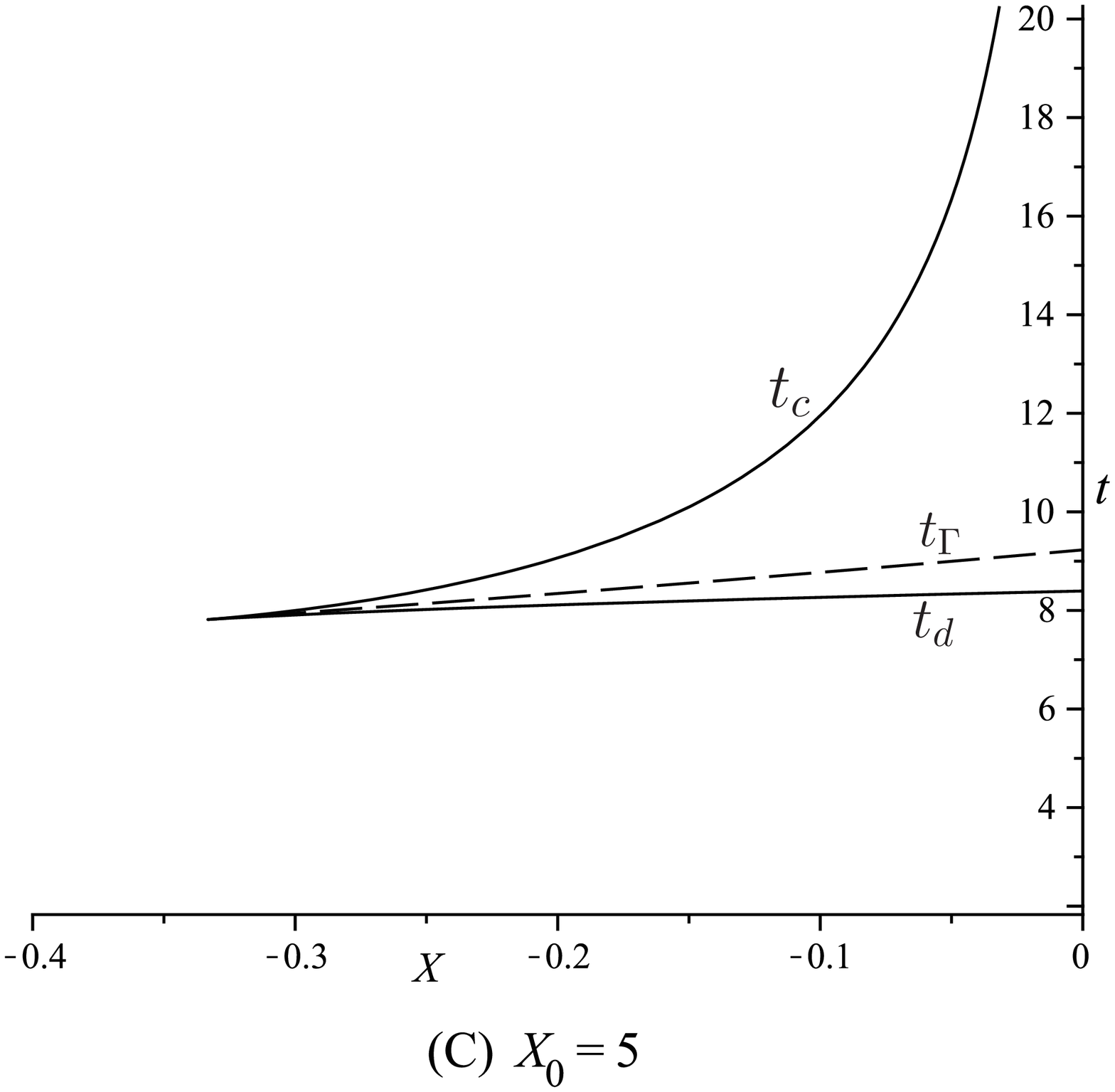}}\\
\end{tabular}
\caption{The curves $t_c$, $t_d$ and $t_\Gamma$ for $X_0=1,\;2,\;5$. }\label{figure4}
\end{center}

\newpage


\begin{thebibliography}{1}
\bibitem{hal} S. Halfin and W. Whitt, Heavy-traffic limits for queues with many exponential servers, Oper. Res. 29 (1981), 567-588.

\bibitem{bro} E. Brockmeyer, H. L. Halstr\o m, and A. Jensen, The Life and Works of A. K. Erlang, Trans. Danish Acad. Tech. Sci. 2 (1948), 277 pp.

\bibitem{pol} F. Pollaczek, \"Uber zwei Formeln aus der Theorie des Wartens vor Schaltergruppen, Elektr. Nachr. 8 (1931), 256-268.

\bibitem{lee}J. S. H. van Leewaarden and C. Knessl, Transient behavior of the Halfin-Whitt diffusion, Stochastic Process. Appl. 121 (2011), 1524-1545.

\bibitem{gam} D. Gamarnik and D. A. Goldberg, On the rate of convergence to stationarity of the $M/M/n$ queue in the Halfin-Whitt regime, 2011, Preprint.

\bibitem{bor} S. Borst, A. Mandelbaum, and M. Reiman, Dimensioning large call centers, Oper. Res. 52 (2004), 17-34.

\bibitem{gan} N. Gans, G. Koole, and A. Mandelbaum, Telephone call centers: tutorial, review and research prospects, Manuf. Serv. Oper. Manag. 5 (2003), 79-141.

\bibitem{jan} A. J. E. M. Janssen, J. S. H. van Leeuwaarden and B. Zwart, Refining square-root safety staffing by expanding Erlang C, Oper. Res. 59 (2011), 1512-1522.

\bibitem{man} A. Mandelbaum and P. Mom\v cilovi\'c, Queues with many servers: the virtual waiting-time process in the QED regime, Math. Oper. Res. 33 (2008), 561-586.

\bibitem{gar} O. Garnett, A. Mandelbaum, and M. Reiman, Designing a call center with impatient customers, Manuf. Serv. Oper. Manag. 4 (2002), 208-227.

\bibitem{jel} P. Jelenkovi\'c, A. Mandelbaum, and P. Mom\v cilovi\'c, Heavy traffic limits for queues with many deterministic servers, Queueing Syst. 47 (2004), 53-69.


\bibitem{ree} J. Reed, The $G/GI/N$ queue in the Halfin-Whitt regime, Ann. Appl. Probab. 19 (2009), 2211-2269.

\bibitem{mag} C. Maglaras and A. Zeevi, Diffusion approximations for a multiclass Markovian service system with ``guaranteed" and ``best-effort" service levels, Math. Oper. Res. 29 (2004) 786-813.

\bibitem{won} R. Wong, Asymptotic Approximation of Integrals, SIAM, Philadelphia, PA, 2001.

\bibitem{ble} N. Bleistein and R. A. Handelsman, Asymptotic Expansions of Integrals, Dover, New York, 1986.

\bibitem{ben} C. M. Bender and S. A. Orszag, Advanced Mathematical Methods for Scientists and Engineers, McGraw-Hill, New York, 1978.

\bibitem{fla} P. Flajolet and R. Sedgcwick, Analytic Combinatorics, Cambridge University Press, Cambridge, 2009.

\bibitem{szp} W. Szpankowski, Average Case Analysis of Algorithms on Sequences, Wiley-Interscience, New York, 2001. 

\bibitem{abr} M. Abramowitz and I. A. Stegun, Handbook of Mathematical Functions with Formulas, Graphs, and Mathematical Tables (Tenth Printing), Dover, New York, 1972.

\bibitem{gra} I. S. Gradshteyn and I. M. Ryzhik, Table of Integrals, Series and Products, Seventh Edition, Elsevier/Academic Press, Amsterdam, 2007.

\bibitem{magnus} W. Magnus, F. Oberhettinger, and R. P. Soni, Formulas and Theorems for the Special Functions of Mathematical Physics, Springer-Verlag, New York, 1966.

\bibitem{tem} N. M. Temme, Parabolic Cylinder Function, NIST Handbook of Mathematical Functions, 303-319, U. S. Dept. Commerce, Washington, DC, 2010. 

\bibitem{olv} F. W. J. Olver, Uniform asymptotic expansions for Weber parabolic cylinder functions of large orders, J. Res.
Nat. Bur. Standards Sect. B 63B (1959), 131-169.


\end{thebibliography}
\end{document}